\documentclass[10pt,leqno]{amsart}
\usepackage{amsmath,amsthm,amsfonts}
\usepackage {latexsym}
\usepackage{amssymb}
\usepackage[dvips]{epsfig}

\renewcommand{\thesection}{\arabic{section}}

\newtheorem{theorem}{Theorem}[section]
\newtheorem{lemma}[theorem]{Lemma}
\newtheorem{prop}[theorem]{Proposition}
\newtheorem{defi}[theorem]{Definition}
\newtheorem{corollary}[theorem]{Corollary}
\newtheorem{remark}[theorem]{Remark}

\newtheoremstyle{prime}
{5pt}   
{5pt}   
{}
{}
{}
{}
{ }
{\thmname{{\bfseries#1}}\thmnumber{ \textbf{#2${\boldsymbol'}$.}}\thmnote{ #3}}

\theoremstyle{prime}

\newtheoremstyle{letter}
{5pt}   
{5pt}   
{}
{}
{}
{}
{ }
{\thmname{{\bfseries#1}}\thmnumber{}\thmnote{ \textbf{#3.}}}

\theoremstyle{letter}

\newtheoremstyle{tilde}
{5pt}   
{5pt}   
{}
{}
{}
{}
{ }
{\thmname{{\bfseries#1}}\thmnumber{}\thmnote{ $\boldsymbol{\tilde{\mbox{\textbf{#3}}}}${\bf.}}}

\theoremstyle{tilde}

\newtheoremstyle{letterprime}
{5pt}   
{5pt}   
{}
{}
{}
{}
{ }
{\thmname{{\bfseries#1}}\thmnumber{}\thmnote{ \textbf{#3$\boldsymbol{'}$.}}}

\theoremstyle{letterprime}

\renewcommand{\theequation}{\thesection .\arabic{equation}}
\let\sect\section
\renewcommand\section{\setcounter{equation}{0}
\gdef\theequation{\thesection .\arabic{equation}}\sect}

\newcommand{\cB}{{\mathcal{B}}}

\newcommand{\cE}{{\mathcal{E}}}

\newcommand{\cJ}{{\mathcal{J}}}

\newcommand{\cR}{{\mathcal{R}}}

\newcommand{\cM}{{\mathcal{M}}}
\newcommand{\cT}{{\mathcal{T}}}

\newcommand{\cF}{{\mathcal{F}}}
\newcommand{\cK}{{\mathcal{K}}}
\newcommand{\cL}{{\mathcal{L}}}
\newcommand{\cP}{{\mathcal{P}}}

\newcommand{\IC}{{\mathbb{C}}}

\newcommand{\IN}{\mathbb{N}}
\newcommand{\IR}{{\mathbb{R}}}
\newcommand{\TT}{{\mathbb{T}}}

\newcommand{\tor}{\TT}
\newcommand{\ZZ}{{\mathbb{Z}}}
\newcommand{\IZ}{{\mathbb{Z}}}

\newcommand{\oN}{{\overline{N}}}

\newcommand{\be}{\begin{eqnarray}}
\newcommand{\ee}{\end{eqnarray}}

\newcommand{\diag}{\mathop{\rm{diag}}}

\newcommand{\dist}{\mathop{\rm{dist}}}

\newcommand{\mes}{\mathop{\rm{mes}\, }}

\renewcommand{\mod}{{\rm{mod}\, }}

\newcommand{\spec}{\mathop{\rm{spec}}}

\newcommand{\w}{\omega}

\newcommand{\ve}{{\varepsilon}}
\newcommand{\vp}{{\varphi}}

\newcommand{\nn}{\nonumber}

\newcommand{\la}{\langle}
\newcommand{\ra}{\rangle}

\newcommand{\op}{\operatorname}

\def\beeq{\begin{equation}}
\def\eneq{\end{equation}}
\def\eps{\varepsilon}
\def\les{\lesssim}

\def\cZ{{\mathcal Z}}

\def\bm{\begin{matrix}}
\def\endm{\end{matrix}}

\def\oN{\bar{N}}
\def\ooN{\Bar{\oN}}
\def\tom{\widetilde\omega}
\def\tcE{\widetilde\cE}
\def\less{\lesssim}

\newcommand{\mybigcup}{\mathop{\textstyle{\bigcup}}\limits}
\newcommand{\mapt}{T_\omega^{(m-1)\ell}x}

\begin{document}

\title[non-perturbative localization]{On non-perturbative Anderson
localization for $C^\alpha$ potentials generated by shifts and
skew-shifts}
\author{Jackson Chan, Michael Goldstein, Wilhelm Schlag}
\thanks{The second author was partially supported by NSERC, and the
third author was partially supported by the NSF}
\address{first and second authors: Department of Mathematics,
University of Toronto, Toronto, Ontario, CANADA}
\email{jchan@math.toronto.edu, gold@math.toronto.edu}

\address{third author: Department of Mathematics, University of
Chicago, Chicago, IL 60637, USA}

\email{schlag@math.uchicago.edu}

\maketitle

\begin{abstract} In this paper we address the question of
proving Anderson localization (AL) for the operator
\[
\bigl[H(x,\omega)\psi\bigr](n) := - \vp(n+1) - \vp(n-1) +
V\bigl(T^n_\omega x\bigr)\psi(n), \qquad n\in\mathbb Z
\]
where $T:\tor^2\to\tor^2$ is either the shift or the skew-shift and
$V$ is only $C^\alpha(\tor^2)$ for some $\alpha>0$. We show that
under the assumption of positive Lyapunov exponents, (AL) takes
place for a.e.~frequency, phase, and energy.
\end{abstract}

\section{Statement of the main results}

Consider the one-dimensional difference Schr\"odinger equation \be
\label{eq:1.eq} \bigl[H(x,\omega)\psi\bigr](n) := - \vp(n+1) -
\vp(n-1) + V\bigl(T^n_\omega x\bigr)\psi(n) = E\psi(n)\ ,\quad n \in
\IZ \ee where $V(x)$ is a real-valued H\"older continuous function
on the two-dimensional torus $\TT^2$, and $T_\omega: \TT^2 \to
\TT^2$ is an ergodic transformation which in this paper will be
either the shift $T_\omega(x,y) = (x,y) + \omega$, $\omega \in
\TT^2$, or the skew-shift $T_\omega(x,y) = (x+y, y +\omega)$,
$\omega\in\tor$. Let $M_{[a,b]}(x,\omega, E)$ be the monodromy
matrix of \eqref{eq:1.eq} on the interval $[a,b]$, i.e.
$$
M_{[a,b]} (x,\omega, E) = \prod^a_{m=b} \begin{bmatrix}
V\bigl(T^m_\omega x\bigr) - E & -1\\
1 & 0\end{bmatrix}
$$
Let $L(\omega, E)$ be the Lyapunov exponent of the cocycle $M_{[1,
N]}(x,\omega, E)$, $N > 0$, i.e.,
$$
L(\omega, E) = \lim_{N\to \infty}\ N^{-1} \int_{\TT^2} \log
\big\|M_N(x,\omega, E)\big\|\, dx
$$

\begin{theorem}\label{th:1.1} Let $V(x)$ be $C^\alpha$--smooth with some $\alpha> 0$.
Assume that $L(\omega, E) > 0$ for all $\omega$ and all $E \in \IR$.
Then there exists $Q$ with $\mes Q = 0$ such that for any $\omega
\notin Q$ there exists $\cE_\omega$ with $\mes\cE_\omega = 0$ such
that for a.a. $x \in \TT$ and all $E \notin \cE_\omega$ the
following assertion holds: if $\bigl[H(x,\omega)\psi\bigr](n) =
E\psi(n)$, $n \in \IZ$, for some polynomially bounded function
$\psi\not\equiv0$, then
$$
|\psi(n)| \le C(x,\omega, E)\exp\bigl(-L(\omega, E)|n|/2\bigr)
$$
for all $n \in \IZ$.
\end{theorem}

If we take the disorder to be large, then we arrive at the following
theorem ($\kappa=\kappa(\alpha)>0$ is a small constant):

\begin{theorem}\label{thm:1.2}
 Let $V(x)$ be $C^\alpha$--smooth with some $\alpha>0$.  Let $L(\omega,\lambda, E)$ be the
 Lyapunov exponent with potential $\lambda V(x)$, $\lambda \in \IR$.  There exists $\lambda_0 = \lambda_0(V)$
 such that for each $|\lambda| > \lambda_0$, there exists a set $\cE_\lambda \in \IR$ such that
 $\mes (\lambda^{-1}\cE_\lambda) < \lambda^{-\kappa}$, and $L(\omega,\lambda, E) > {1\over 2} \log |\lambda|$
 for all $\omega$ and all $E \notin \cE_\lambda$.
 Moreover, for each $|\lambda| > \lambda_0$ there exists $Q_\lambda$ with $\mes Q_\lambda = 0$
 such that for each $\omega \notin Q_\lambda$ there exists $\cE_{\lambda, \omega}$ with $\mes\cE_{\lambda,\omega} = 0$
 such that for a.a.~$x$ the following assertion holds:
for all $E\not\in \cE_\lambda\cup \cE_{\lambda, \omega}$, if
$\left[H(\lambda, x, \omega)\psi\right](n) = E\psi(n)$, $n \in \IZ$
for some polynomially bounded function $\psi(n)$, then
 $$
 |\psi(n)| \le C(\lambda, x, \omega, E) \exp\bigl(-L(\lambda, \omega, E)|n|/2\bigr), n \in \IZ
 $$
 Here $H(\lambda, x, \omega)$ stands for the Schr\"odinger operator \eqref{eq:1.eq} with $\lambda V(x)$ in the role of $V(x)$.
 \end{theorem}

The main novel feature in these theorems is the low regularity
of the potential function~$V$. While much is known about the
case of analytic $V$, see Bourgain-Goldstein\cite{BG} and
Bourgain-Goldstein-Schlag\cite{BGS} and for the case of almost
Mathieu, Jitmoriskaya\cite{Jit}, comparatively little is known about
to non-analytic category. Klein\cite{Kl} proved (AL) for positive
Lyapunov exponents, and potentials from the Gevrey class.
Bjerklov\cite{Bje} showed (AL) and positive Lyapunov exponents for
$C^1$ potentials, large disorder, and off a set of energies of
positive measure. Chan\cite{C} proved (AL) for large disorder, for
generic $C^3$ potentials in a suitable sense and for a.e.~energy.

The methods of this paper are in spirit related to those
of~\cite{BG} and \cite{BGS}. Thus, the main ingredients are large
deviation theorems and the elimination of resonant frequencies.
Analyticity has so far played a crucial role in obtaining these
tools. Hence, we need to take a very different route here. A basic
principle in this paper is to reduce matters to the study of the
eigenvalues as functions parameterized by the phase (the so-called
Rellich functions). Firstly, we note that by Weyl's comparison
theorem the determinant of the Hamiltonian $H_{[-N,N]}(x,\omega)-E$
of~\eqref{eq:1.eq} at energy $E$ restricted to the interval $[-N,N]$
is comparable to the product of the determinants of the Hamiltonians
at energy $E$ corresponding to $[n_j,n_{j+1}]$ where
$[-N,N]=\bigcup_{j=1}^J[n_j,n_{j+1}]$; albeit, their ratio can be
very large, namely $\eta^{-J}$ where $\eta $ is the distance of the
spectrum of $H_{[-N,N]}(x,\omega)$ to~$E$. This is one source of
energy removal: evidently we will need to control~$\eta$. Secondly,
first order eigenvalue perturbation shows that for $C^1$ potentials
these functions are again $C^1$ (in fact, they also inherit the
H\"older regularity of the potential). Thus, we can apply the
ergodic theorem to these individual functions and then sum up over
all of them to obtain a large deviation theorem for the sum of
shifts of logarithms of Dirichlet determinants. A second source of
energy removal arises at this stage: we will need to exclude those
energies that serve as critical values of the Rellich functions.

\medskip

We conclude this introduction with a heuristic discussion of the
wider context of our results. In particular, we shall mention some
natural ramifications that Theorems~\ref{th:1.1} and~\ref{thm:1.2}
above appear to possess.

\medskip

\begin{itemize}
\item[(1)]
It seems natural to combine the methods of this paper with those
based on subharmonicity (developed in \cite{BG}, \cite{BGS},
\cite{GS1}, \cite{GS2}) to show that for the case of analytic
potentials the elimination of resonances for the bulk of energies
(as in this work) suffices to obtain complete (AL) at all energies.
\item[(2)]
The crucial  component needed to make
 progress in (1) consists of a count of the number of the
Dirichlet eigenvalues of the operator $H_{[-N, N]} (x, \omega)$
which fall into the set $\cE_\omega$. Recall that the latter is the
set of ``forbidden" energies, which appear in the elimination of
resonances in this work.
\item[(3)]
It is not hard to prove that the Hausdorff dimension of the set
$\cE_\omega$ is equal to zero in case of analytic potentials (as
well as for $C^K$-smooth potentials with large $K$). It seems that
there is a possibility to use this fact and to modify the method
of~\cite{GS2}, which is based on the multi-scale (or avalanche
principle) expansion of the function $\log|f_{[-N, N]} (z, \omega ,
E)|$, to evaluate the averaged (in phases) number of eigenvalues
falling into $\cE_\omega$. Recall that the expansion itself is valid
for the shift and skew-shift, provided the Lyapunov exponent is
positive (see \cite{GS2}).
\item[(4)]
It is not clear to what extent the results of Theorems~\ref{th:1.1}
and~\ref{thm:1.2} are optimal for smooth potentials. For example, it
is unclear whether the Hausdorff dimension of the set $\cE_\omega$
also vanishes for H\"older continuous potentials. In this context,
it seems natural to ask the following question:

{\em Are there any smooth potentials $V(x),\ x\in \tor$, with $L(E,
\omega )>0$ for all $E, \omega$ and with purely singular continuous
spectrum for a.a.~$\omega$? In other words, the spectrum of $H(x,
\omega)$ for such a potential would fall into the set $\cE_\omega$
from Theorem~\ref{th:1.1} for a.a.~$\omega$.}

On the other hand, it looks promising to modify the technology of
\cite{C} to show that for "generic" smooth potentials the set
$\cE_\omega$ does not contribute anything to the spectrum, i.e.,
that  complete (AL) takes place in Theorem~\ref{th:1.1}.
\item[(5)] The large deviations estimates and the process of elimination of
resonances which are developed in Sections 2, 3 of this paper can
also be established for the case of quasi-periodic Schr\"odinger
operators on the lattice. Moreover, not only the multidimensional
Laplacian, but also its long-range versions can be treated in this
fashion. Taking this into account, it seems plausible to establish
(AL) for higher-dimensional quasi-periodic lattice models in the
case of analytic potentials at large values of $\lambda$. In
particular, this would prove the absolute continuity of the spectrum
of these models in the regime of small values of $\lambda$.
\end{itemize}
In the work \cite{GK} the ideas of this work are modified for the
analysis of localization at almost all energies in the case of
random potentials with fast decaying correlations and in particular
for the potentials generated by the doubling map on circle.

\section{Large deviation estimates for the averages of shifts and skew-shifts of logarithms of $C^1$-smooth functions}\label{sec:2}

In this section, we develop a general framework of averaging of
functions of the form $\log|f(x)-\xi|$ over orbits of the shift
and skew-shift. The reader will find the relevant quantitative
ergodic theorems in the appendices. In the process we shall need to
remove those values of $\xi$ for which the function~$\log|f(x)-\xi|$
becomes too singular. This is comparatively easy: it will only
require Fubini's theorem. Throughout, we assume that the potential
is $C^1$ for the sake of simplicity. The generalization to H\"older
classes is elementary, see Appendix C.

\begin{defi}\label{def:basic}
Suppose $f\in C^m\left(\TT^2\right)$.  If $\alpha + \beta\leq m$,
let $B_{\alpha,\beta}(f) := \underset{x\in\TT^2}{\max}\,
\left|\partial_{x_1}^{\alpha}\, \partial_{x_2}^{\beta}f(x)\right|$.
Also, if $k\leq m$, let $B_k(f):=\underset{\alpha + \beta\leq
k}{\max} B_{\alpha,\beta}(f)$.  In particular, $B_0(f) =
\underset{x\in\TT^2}{\max}\, |f(x)|$. Throughout this paper, we let
\[S_f(\xi,\delta) := \left\{ x\in\TT^2 : |f(x) - \xi| <
\delta\right\}\] denote level sets of $f$.
\end{defi}

 \noindent Let $f\in C^1 \left(\TT^{2}\right)$.
\ Our first goal is to estimate
\begin{equation}
\label{eq:orbits} \#\left\{ k\in\IN : 1\leq k\leq N,
\left|f\left(T^kx\right) - \xi\right| < \delta\right\} = \#\left\{
k\in\IN : 1\leq k\leq N, T^k x\in S_f(\xi,\delta)\right\}
\end{equation}
where $\xi\in\IR$ is a parameter, $0 < \delta < 1$, $T :
\TT^2\rightarrow\TT^2$ is the shift $T(x_1,x_2) = (x_1 + \omega_1,
x_2 + \omega_2)$, or the skew-shift $T(x_1, x_2) = (x_1 + x_2, x_2 +
\omega)$ (addition here is always $\mod\ZZ^2$). In order to answer
\eqref{eq:orbits}, we will need to use quantitative ergodic
properties of these transformations. As a preliminary step, we
introduce the following functions for the purpose of mollifying
given $C^1$ functions.

\begin{defi}
\label{def:moll} Given $\tau
> 0$, let $h_\tau\in C^4(\IR)$ be $1$-periodic such that
\begin{itemize}
\item $h_\tau\geq 0$
\item $\op{supp}h_\tau\subset\bigcup_{k\in\IZ}[k - \tau, k + \tau]$
\item $\int_0^1h_\tau(y)dy = 1$
\item
$\max_{y\in\IR}\left|\left(\frac{\displaystyle d}{\displaystyle
dy}\right)^{m}h_\tau(y)\right| \lesssim \tau^{-(m + 1)}\,
$ for $m\leq 4$.
\end{itemize}
Moreover, we set $ \tilde h_\tau(x_1, x_2)  = h_\tau(x_1)h_\tau(x_2)
$.
\end{defi}

The following lemma is a well-known quantitative statement
concerning the  mollifiers of a given $C^1$ function.

\begin{lemma}\label{lem:2.C4} Given $\varphi\in C^1(\TT^2)$ and $\tau>0$, define
\[
\psi(x_1, x_2) : = \int_{\TT^2}\varphi(y_1, y_2)\tilde h_\tau(x_1 -
y_1, x_2 - y_2)\, dy
\]
Then $\psi\in C^4(\TT^2)$ satisfies
\begin{itemize}
\item[(1)]$\max_{x\in\TT^2}|\varphi(x) - \psi(x)| \lesssim B_1(\varphi)\tau$
\item[(2)]$B_4(\psi) \lesssim B_0(\varphi)\tau^{-4}$
\end{itemize}
\end{lemma}

Turning to the dynamics, we will of course need a Diophantine
condition. Throughout this paper, constants will be allowed to
depend on the constants appearing in this definition.

\begin{defi}
\label{def:Dioph} Under a {\em Diophantine condition} on $\omega$ we
shall mean the following: if $T$ is the shift, then we will assume $\|
k_1\omega_1 + k_2\omega_2 \|
> c(|k_1| + |k_2|)^{-A}$ for $(k_1, k_2)\in\IZ^2\setminus\{ 0\}$
where $c > 0$, $A > 2$. If $T$ is the skew-shift, then we will assume $\|
k\frac{\omega}{2}\| > c|k|^{-(1 + \eps)}$ for $k\in\IZ\setminus\{
0\}$ where $0 < c < 1$ and $0 < \eps \ll 1$.
\end{defi}

\noindent The following proposition is a quantitative version of the
ergodic theorem for smooth (i.e., $C^4$) functions.

\begin{prop}
\label{prop:2.largen} For sufficiently large $N$, one has (with $T$
being either the shift or skew-shift)
\[
\left|\frac{1}{N}\sum_{m = 1}^N\psi\left(T^mx\right) - \la\psi\ra\right| \lesssim B_4(\psi)N^{-\sigma}
\]
for all $x\in\TT^2$. The constant $\sigma$ depends on the parameters
in Definition~\ref{def:Dioph}.
\end{prop}
\begin{proof}
The proof can be found in Appendix~A.
\end{proof}

We now turn to estimating \eqref{eq:orbits}. It will be convenient
to work with $C^1$ functions instead of indicators of level sets.
The following lemma introduces the standard transition between the
two.

\begin{lemma}\label{lem:2.list}
Given $\delta>0$ small, let $\chi_{\delta}\in C^1(\IR)$ be such that
\begin{itemize}
\item $0 \leq \chi_{\delta} \leq 1$
\item $\chi_{\delta}(y) = 1$ for $ y\in[-\delta,\delta]$
\item $\op{supp}\chi_{\delta}\subset[-2\delta,2\delta]$
\item $\max_{y\in\IR}\left|\chi_{\delta}'(y)\right|\lesssim\delta^{-1}$
\end{itemize}
Then the following holds:
\begin{enumerate}
\item $\#\left\{ 1\leq k\leq N : T^kx\in S_f(\xi, \delta)\right\} \leq \sum_{k=1}^N\chi_{\delta}\left(f\left(T^kx\right) - \xi\right)$
\item $\mes S_f(\xi, \delta)\leq\la\chi_{\delta}(f(\cdot) - \xi)\ra = \int_{\TT^2}\chi_{\delta}\bigl(f(x) - \xi\bigr)\,dx \leq \mes S_f(\xi, 2\delta)$
\end{enumerate}
for any real $\xi$ and positive integer $N$.
\end{lemma}

We can now apply Proposition~\ref{prop:2.largen} to deduce the
required bound on~\eqref{eq:orbits}.

\begin{corollary}
\label{cor:2.cprime} Let $f\in C^1(\TT^2)$ and suppose $\omega$ is Diophantine. Then for any $\xi\in\IR$,
$\delta > 0$, one has
\[
\frac{1}{N}\#\left\{ 1\leq k\leq N : T^kx\in S_f(\xi,\delta)\right\} \lesssim\mes S_f(\xi,2\delta) + \left(1 + B_1(f)\right)\delta^{\frac{1}{2}}
\]
for all $x\in\TT^2$ provided $N\geq\delta^{-\frac{10}{\sigma}}$. Here $\sigma>0$ is the small constant from the ergodic theorem,
see Proposition~\ref{prop:2.largen}.
\end{corollary}

\begin{proof}
Using the notations of Lemma~\ref{lem:2.list}, we have to estimate
$\frac{1}{N}\sum_{k = 1}^N\chi_{\delta}\left(f(T^kx) - \xi\right)$.
Note that $\varphi(x) = \chi_{\delta}\left(f(x) - \xi\right)$ is
$C^1$ on $\TT^2$, $B_1(\varphi) \lesssim B_1(f)\delta^{-1}$. By
Lemma~\ref{lem:2.C4}, given $\tau > 0$, there is $\psi\in
C^4\left(\TT^2\right)$ such that
\begin{itemize}
\item[(1)]$\max_{x\in\TT^2}\left|\varphi(x) - \psi(x)\right| \lesssim B_1(f)\delta^{-1}\tau$
\item[(2)]$B_4(\psi) \lesssim \tau^{-6}$
\end{itemize}
Due to (1), $\left|\la\varphi\ra - \la\psi\ra\right| \lesssim B_1(f)\delta^{-1}\tau$ and $\left|\frac{1}{N}\sum_{k=1}^N\varphi\left(T^kx\right) - \frac{1}{N}\sum_{k=1}^N\psi\left(T^kx\right)\right| \lesssim B_1(f)\delta^{-1}\tau$ for all $x\in\TT^2$. Therefore,
\begin{align*}
\frac{1}{N}\#\left\{ 1\leq k\leq N : T^kx\in S_f(\xi, \delta)\right\} & \leq \frac{1}{N}\sum_{k=1}^N\varphi\left(T^kx\right) \\
& \lesssim \la\varphi\ra + B_4(\psi)N^{-\sigma} + B_1(f)\delta^{-1}\tau \\
& \leq \mes S_f(\xi, 2\delta) + \tau^{-4}N^{-\sigma} + B_1(f)\delta^{-1}\tau.
\end{align*}
The assertion follows if we take $\tau = \delta^{\frac{3}{2}}$.
\end{proof}

Since we are not making a non-degeneracy assumption on $f$ (in
particular, $f$ may be constant) it will be necessary to remove
certain values of $\xi$ for which $S_f(\xi,\delta)$ is very large.
This can be done easily by means of Fubini's theorem.

\begin{corollary}
\label{cor:2.xiset} Let $\omega$ be Diophantine. Given $\delta > 0$, there exists a set $\cE_{\delta}\subset\IR$,
$\mes\cE_{\delta} \lesssim \delta^{\frac{1}{2}}$ such that for any
$\xi\notin\cE_{\delta}$, one has $\frac{1}{N}\#\left\{ 1\leq k\leq N
: T^kx\in S_f(\xi, \delta)\right\} \lesssim \left(1 +
B_1(f)\right)\delta^{\frac{1}{2}}$ for all $x\in\TT^2$ provided
$N\geq\delta^{-\frac{10}{\sigma}}$.
\end{corollary}
\begin{proof} Clearly,
$\mes\left\{(x,\xi)\in\TT^2\times\IR : x\in S_f(\xi,
2\delta)\right\} = 4\delta$. By Fubini's Theorem, there exists
$\cE_{\delta}\subset\IR$, $\mes\cE_{\delta}\leq
4\delta^{\frac{1}{2}}$, such that $\mes S_f(\xi,
2\delta)\leq\delta^{\frac{1}{2}}$ for $\xi\notin\cE_{\delta}$. The
assertion now follows from Corollary~\ref{cor:2.cprime}.
\end{proof}

\begin{remark}\label{rem:2.8FAT}
Let $\cE_\delta$ be as follows
\[
\cE_\delta = \{\xi\::\:\mes S_f(\xi,\delta)>\delta^{\frac12} \}
\]
Given an arbitrary subset $\cE\subset\IR$ and $r>0$ introduce
\[
[\cE](r)=\{\xi\in\IR\::\: \dist(\xi,\cE)<r\}
\]
Note that if $\xi\in\cE_\delta$ and $|\xi_1-\xi|<r$, then
\[ S_f(\xi_1,\delta)\subset S_f(\xi,\delta+r) \]
Define
\[ \cE_{\delta,\delta_1} = \{\xi\::\:\mes S_f(\xi,\delta)>\delta_1^{\frac12} \}\]
Then $\mes \cE_{2\delta,\delta}\lesssim \delta^{\frac12}$. On the other hand,
$[\cE_\delta](r)\subset \cE_{2\delta,\delta}$ for $r<\delta$. In particular,
$\mes [\cE_\delta](r)\lesssim \delta^{\frac12}$.
\end{remark}

Our next goal is to estimate
\begin{equation}
\label{eq:LDT} \mes\Big\{ x :
\Big|\frac{1}{N}\sum_{k=1}^N\log|f\left(T^kx\right) - \xi| -
\la\log|f(\cdot) - \xi|\ra\Big| > \delta\Big\}.
\end{equation}
This is of course motivated by the {\em large deviation theorems} in
\cite{BG}, \cite{GS1}, and \cite{GS2}. As a first step, note that
\begin{equation}
\label{eq:oldstar}
\begin{aligned}
\frac{1}{N}\hspace*{-12pt}\sum_{\substack{1\leq k\leq N\\ T^kx\notin S_f(\xi,\delta)}}\hspace*{-16pt}\log\left|f\left(T^kx\right) - \xi\right| & + \frac{1}{N}\#\left\{ 1\leq k \leq N : T^kx\in S_f(\xi, \delta)\right\}\cdot\log\left(\min_{1\leq k\leq N}\left|f\left(T^kx\right) - \xi\right|\right) \\
& \leq \frac{1}{N}\sum_{k=1}^N\log\left|f\left(T^kx\right) -
\xi\right|
 \leq \frac{1}{N}\hspace*{-12pt}\sum_{\substack{1\leq k\leq N\\
T^kx\notin S_f(\xi,\delta)}}\log\left|f\left(T^kx\right)-\xi\right|
\end{aligned}
\end{equation}
since $\delta < 1$. The following lemma shows that the term
involving the minimum in~\eqref{eq:oldstar} can be controlled at the
expense of removing an exponentially small set of~$x$. In what
follows, we use the notation
\begin{equation}
\label{eq:J_def} [-B_0(f) ,B_0(f)] = : \cJ_0(f)
\end{equation}
where $f$ will be a given $C^1$ function.

\begin{lemma}\label{lem:2.6}
Let $\kappa>0$ be arbitrary. There exists $\cE(N)\subset\IR$,
$\mes\cE(N) < \exp\left(-\frac{N^\kappa}{2}\right)$ such that for
any $\xi\notin\cE(N)$, one has
\begin{equation}
\label{eq:2.logmin}
\mes\Big\{ x\in\TT^2 : \log\Big(\min_{1\leq k\leq N}|f(T^kx) -
\xi|\Big) < -N^\kappa\Big\} \leq
\exp\left(-\frac{1}{4}N^\kappa\right).
\end{equation}
provided $N\ge N_0(\kappa)$.
\end{lemma}

\begin{proof}
This follows from Fubini's theorem. Indeed, with $\cJ_0(f)$ as
above,
\begin{align}
& \int_{\cJ_0(f)} \mes\left\{ x\in\TT^2 : \log\left(\min_{1\leq
k\leq N}\left|f\left(T^kx\right) - \xi\right|\right) <
-N^\kappa\right\}\, d\xi \label{eq:mes_int}
\\
&\leq \sum_{k=1}^N \int_{\cJ_0(f)}\mes\left\{ x\in\TT^2 :
\log\left|f\left(T^kx\right) - \xi\right|<
-N^\kappa\right\} \, d\xi \nn \\
&\leq 2N  e^{-N^\kappa}  \nn
\end{align}
Hence, we can remove a set of $\xi$ of measure not exceeding
$e^{-\frac12 N^\kappa}$ so that the integrand in~\eqref{eq:mes_int}
is at most $e^{-\frac14 N^\kappa}$, as claimed.
\end{proof}

\begin{remark}\label{rem:2.disect} One can see that the
following
version of Lemma~\ref{lem:2.6} holds: For any $x_2^{(0)} \in \TT$
there exists $\cE^{(1)}(x_2^{(0)},N)\subset \IR$ with measure $\le
\exp\bigl(-N^{\kappa}/2\bigr)$ such that for any $\xi \notin
\cE^{(1)}(x_2^{(0)}, N)$ one has
$$
\mes\left\{x_1 \in \TT: \log \left(\min_{1 < k \le
N}\big|f\bigl(T^k(x_1, x_2^{(0)}) - \xi\big|\right) < -
N^\kappa\right\} \le \left(-\tfrac14 N^\kappa\right)\ .
$$
\end{remark}

Combining this lemma with \eqref{eq:oldstar} we obtain the following:

\begin{corollary}
Let $\delta>0$ and $\kappa>0$ be small. Then for all
$\xi\notin\cE(N)$ there exists $B(\xi)\subset\TT^2$, $\mes
B(\xi)\leq\exp\left(-\frac{1}{4}N^\kappa\right)$ such that for any
$x\not\in B(\xi)$ one has
\begin{align*}
\frac{1}{N}\hspace*{-12pt}\sum_{\substack{1\leq k\leq N\\ T^kx\notin S_f(\xi,\delta)}}\hspace*{-16pt}\log\left|f\left(T^kx\right) - \xi\right| - \frac{1}{N^{1-\kappa}}\sup_{y\in\TT^2}\left(\#\left\{ 1\leq k\leq N : T^ky\in S_f(\xi,\delta)\right\}\right) & \leq \frac{1}{N}\sum_{k=1}^N\log\left|f\left(T^kx\right) - \xi\right| \\
& \leq \frac{1}{N}\hspace*{-12pt}\sum_{\substack{1\leq k\leq N\\ T^kx\notin S_f(\xi,\delta)}}\hspace*{-16pt}\log\left|f\left(T^kx\right) - \xi\right|
\end{align*}
Here $N\ge N_0(\kappa)$ is a positive integer.
\end{corollary}

In order to bound the averages on the left and right-hand sides here
we introduce the following auxiliary function.

\begin{defi}\label{def:rho}
Henceforth, given $\delta>0$ small we define $\rho=\rho_{\delta}\in
C^2(\IR)$ to be such that
\begin{itemize}
\item $\rho(y) = |y|$ if $|y| \geq\delta$ and $\rho(y)\ge|y|$ for all $y$
\item $\frac{\delta}{2}\leq\rho(y)\leq\delta$ if $y\in(-\delta,\delta)$
\item $\max_{y\in\IR}\left|\rho''(y)\right|\lesssim\delta^{-1}$
\end{itemize}
\end{defi}

The significance of this definition can be seen from the following
lemma.

\begin{lemma}
If $f\in C^1(\TT^2)$, and $0<\delta<1$, then for $\xi \notin \cE_{2\delta,\delta}$
\[
\frac{1}{N}\sum_{k=1}^N\log\rho\left(f\left(T^kx\right) - \xi\right) \leq \frac{1}{N}\hspace*{-12pt}
\sum_{\substack{1\leq k\leq N\\ T^kx\notin S_f(\xi,\delta)}}\hspace*{-16pt}\log\left|f\left(T^kx\right) - \xi\right| \leq \frac{1}{N}\sum_{k=1}^{N}\log\rho\left(f\left(T^kx\right) - \xi\right) + \left(1 + B_1(f)\right)\delta^{\frac13}
\]
\end{lemma}

\begin{proof}
The first inequality is clear since $\delta < 1$. For the second,
note that
\[
\frac{1}{N}\hspace*{-12pt}\sum_{\substack{1\leq k\leq N\\ T^kx\notin
S_f(\xi,\delta)}} \hspace*{-16pt}\log\left|f\left(T^kx\right) -
\xi\right| -\frac{1}{N}\#\left\{ 1\leq k\leq N : T^kx\in
S_f(\xi,\delta)\right\}|\log(\delta/2)|
\leq\frac{1}{N}\sum_{k=1}^N\log\rho\left(f\left(T^kx\right) -
\xi\right)
\]
By Corollary~\ref{cor:2.xiset},
\begin{align*}
\frac{1}{N}\#\left\{ 1\leq k\leq N : T^kx\in S_f(\xi,\delta)\right\}\cdot\Big|\log\frac{\delta}{2}
\Big| & \lesssim \left(1 + B_1(f)\right)\delta^{\frac{1}{2}}\Big|\log\frac{\delta}{2}
\Big|
\end{align*}
and the lemma follows.
\end{proof}

Note that we also need to exclude a set of $\xi$ in order to prevent
the averages $\la\log\left|f(\cdot) - \xi\right|\ra$ in
\eqref{eq:LDT} from being too large.

\begin{lemma}
\label{2:lem.capL} Given $R > 0$ there exists $\cL_R\subset
\cJ_0(f)$, $\mes\cL_R\lesssim B_0 R^{-1}$, such that for any
$\xi\notin\cL_R$, one has
\begin{itemize}
\item[(1)]$\int_{\TT^2}\bigl|\log\left|f(x) - \xi\right|\bigr|^2dx\leq\left(\log B_0\right)^2 R$
\item[(2)]
$\left|\int_{\TT^2}\log\left|f(x) - \xi\right|dx -
\int_{\TT^2\setminus S_f(\xi,\delta)}\log\left|f(x) -
\xi\right|dx\right|\leq\left(\log B_0\right)R^{\frac{1}{2}}
\left[\mes S_f(\xi,\delta)\right]^{\frac{1}{2}}$
\end{itemize}
\end{lemma}

\begin{proof}
Since $\bigl|\log\left|f(x) - \xi\right|\bigr|^2 >0$, by Fubini's
theorem,
\begin{align*}
\int_{\cJ_0(f)}\int_{\TT^2}\bigl|\log\left|f(x) - \xi\right|\bigr|^2dx\,d\xi & =
\int_{\TT^2}\int_{\cJ_0(f)}\bigl|\log\left|f(x) - \xi\right|\bigr|^2d\xi\, dx \\
& \lesssim B_0\left(\log B_0\right)^2
\end{align*}
Hence, there exists $\cL_R\subset \cJ_0(f)$, $\mes\cL_R\leq B_0
R^{-1}$ such that $\int_{\TT^2}\bigl|\log\left|f(x) -
\xi\right|\bigr|^2dx\leq\left(\log B_0\right)^2 R$ for
$\xi\notin\cL_R$. This proves (1). To prove (2), we use
Cauchy-Schwarz:
\begin{align*}
\left|\int_{\TT^2}\log\left|f(x) - \xi\right|dx - \int_{\TT^2\setminus S_f(\xi,\delta)}\log\left|f(x) - \xi\right|\,dx\right|
& \leq \int_{S_f(\xi,\delta)}\bigl|\log\left|f(x) - \xi\right|\bigr|\,dx \\
& \leq \left[\int_{\TT^2}\bigl|\log\left|f(x) - \xi\right|\bigr|^2dx\right]^{\frac{1}{2}}
\left[\mes S_f(\xi,\delta)\right]^{\frac{1}{2}} \\
& \leq \left(\log B_0\right)R^{\frac{1}{2}}\left[\mes
S_f(\xi,\delta)\right]^{\frac{1}{2}}
\end{align*}
and the lemma follows.
\end{proof}

Now the same for the regularized functions $\rho\left(f(x) -
\xi\right)$:

\begin{corollary}
\label{cor:2.11} There exists $\cM_{\delta}\subset \cJ_0(f)$, $\mes
\cM_{\delta}\lesssim B_0\delta^{\frac{1}{4}}$, such that for any
$\xi\notin \cM_{\delta}$ one has
\[
\left|\int_{\TT^2}\log\rho\left(f(x) - \xi\right)\,dx -
\int_{\TT^2}\log\left|f(x) - \xi\right|\,dx \right|
\lesssim\left(\log B_0\right)\delta^{\frac{1}{8}}.
\]
\end{corollary}

\begin{proof} Using the notations of the previous lemma,
suppose $\xi\notin\cL_R\cup\cE_\delta$ where $\cE_\delta$ is as in Remark~\ref{rem:2.8FAT}. Then
\begin{align*}
\left|\int_{\TT^2}\log\rho\left(f(x) - \xi\right)dx - \int_{\TT^2}\log\left|f(x) - \xi\right|dx\right| & \le \int_{S(\xi,\delta)}\bigl|\log\rho\left(f(x) - \xi\right)\bigr|dx + \int_{S(\xi,\delta)}\bigl|\log\left|f(x) - \xi\right|\bigr|dx \\
& \leq \left[\mes S(\xi,\delta)\right]\left|\log\frac{\delta}{2}\right| + \left(\log B_0\right)R^{\frac{1}{2}}\left[\mes S_f(\xi,\delta)\right]^{\frac{1}{2}} \\
& \leq \delta^{\frac{1}{2}}\left|\log\frac{\delta}{2}\right| +
\left(\log B_0\right)R^{\frac{1}{2}}\delta^{\frac{1}{4}}.
\end{align*}
Take $R = \delta^{-\frac{1}{4}}$, $\cM_\delta :=
\cL_{\delta^{-\frac{1}{4}}}\cup\cE_\delta$. Then
$\left|\int_{\TT^2}\log\rho\left(f(x) - \xi\right)dx -
\int_{\TT^2}\log\left|f(x) - \xi\right|dx\right| \lesssim \left(\log
B_0\right)\delta^{\frac{1}{8}}$ for any $\xi\notin \cM_\delta$.
Moreover, $\mes \cM_\delta\lesssim B_0\delta^{\frac{1}{4}}$.
\end{proof}

We are finally ready to state a large deviation theorem for averages
of $C^1$ functions, albeit at the expense of removing some dangerous
level sets (i.e., values of $\xi$).

\begin{theorem}\label{thm:A}
Let $f\in C^1\left(\TT^2\right)$ and suppose $\omega$ satisfies a
Diophantine condition. Then there is a sufficiently small $\kappa>0$
so that for all large  $N\ge N_0(\kappa)$  there exists
$\cT(N)\subset\cJ_0(f)$, $\mes\cT(N) < N^{-\kappa}$, such that for
any $\xi\in\cJ_0(f)\setminus\cT(N)$ one has
\begin{equation}\label{eq:largedt}
\mes\Big\{ x :
\Big|\frac{1}{N}\sum_{k=1}^N\log\left|f\left(T^kx\right) -
\xi\right| - \la\log|f(\cdot) - \xi|\ra\Big|
 > N^{-\kappa}\Big\}\leq\exp\left(-N^\kappa\right)
\end{equation}
Moreover,  one has
\begin{equation}\label{eq:uniform}
\sup_{\substack{\cB\subset[1,N]\\ \#\cB<N^{1-2\kappa}}}\sup_{x\in\TT^2}\frac{1}{N}\sum_{k\in[1,N]\setminus\cB}\log\left|f\left(T^kx\right) -
\xi\right| \le \la\log|f(\cdot) - \xi|\ra+ N^{-\kappa}
\end{equation}
for any $\xi\in\cJ_0(f)\setminus\cT(N)$.
\end{theorem}

\begin{proof} Let $\rho$ be as in Definition~\ref{def:rho} with $\delta$ to be specified later. Then
\begin{align*}
\left|\frac{1}{N}\sum_{k=1}^N\log\left|f\left(T^kx\right) - \xi\right| - \la\log\left|f(\cdot) - \xi\right|\ra\right| & \leq \left|\frac{1}{N}\sum_{k=1}^N\log\left|f\left(T^kx\right) - \xi\right| - \frac{1}{N}\sum_{k=1}^N\log\rho\left(f\left(T^kx\right) - \xi\right)\right| \\
& \hspace*{12pt} + \left|\frac{1}{N}\sum_{k=1}^N\log\rho\left(f\left(T^kx\right) - \xi\right) - \int_{\TT^2}\log\rho\left(f(y) - \xi\right)dy\right| \\
& \hspace*{12pt} + \left|\int_{\TT^2}\log\rho\left(f(y) - \xi\right)dy - \int_{\TT^2}\log\left|f(y) - \xi\right|dy\right|
\end{align*}
Let $B_{N,\xi} = \left\{ y\in\TT^2 : \min_{1\leq k\leq
N}\left|f\left(T^ky\right) - \xi\right| \leq e^{-N^\kappa}\right\}$, $\cE(N)$ be as in Lemma~\ref{lem:2.6} and $\cE_{2\delta,\delta}$ be as in Remark~\ref{rem:2.8FAT}.
If $\xi\notin\cE(N)\cup \cE_{2\delta,\delta}$, and $x\notin B_{N,\xi}$, then
\begin{align*}
&\left|\frac{1}{N}\sum_{k=1}^N\log\left|f\left(T^kx\right) - \xi\right| - \frac{1}{N}\sum_{k=1}^N\log\rho\left(f\left(T^kx\right) - \xi\right)\right| \\
& \leq \frac{1}{N}\#\left\{ 1\leq k\leq N : T^kx\in S_f(\xi,\delta)\right\}
 \left[\left|\log\min_{1\leq k\leq N}\left|f\left(T^kx\right) - \xi\right|\right| + \left|\log\frac{\delta}{2}\right|\right] \\
& \lesssim \left(1 + B_1(f)\right)\delta^{\frac{1}{2}}\left[N^\kappa + \left|\log\frac{\delta}{2}\right|\right]
\end{align*}
Moreover, $\mes B_{N,\xi} \leq \exp\left(-\frac{1}{4}N^\kappa\right)$.
Let $\varphi(y) = \log\rho\left(f(y) - \xi\right)$. By Lemma~\ref{lem:2.C4}, for any $\tau >
0$, there is $\psi\in C^4\left(\TT^2\right)$ such that
\begin{itemize}
\item[(1)]$\max_{y\in\TT^2}\left|\varphi(y) - \psi(y)\right|\lesssim B_1(\varphi)\tau\leq B_1(f)\delta^{-1}\tau$
\item[(2)]$B_4(\psi)\lesssim B_0(\varphi)\tau^{-4}\leq(|\log(\delta/2)|+1+B_0(f))\tau^{-4}$
\end{itemize}
Then with some $\sigma>0$,
\begin{equation}
\label{eq:smooth_ergodic}
\begin{aligned}
\left|\frac{1}{N}\sum_{k=1}^N\varphi\left(T^kx\right) - \la\varphi\ra\right| & \leq \left|\frac{1}{N}\sum_{k=1}^N\varphi\left(T^kx\right) - \frac{1}{N}\sum_{k=1}^N\psi\left(T^kx\right)\right| + \left|\frac{1}{N}\sum_{k=1}^N\psi\left(T^kx\right) - \la\psi\ra\right| + \left|\la\psi\ra - \la\varphi\ra\right| \\
& \lesssim B_1(f)\delta^{-1}\tau + (|\log\delta/2|+1+B_0(f))\tau^{-4}N^{-\sigma}
\end{aligned}
\end{equation}
provided $N$ is sufficiently large (see Proposition~\ref{prop:2.largen}).
By Corollary~\ref{cor:2.11} there exists $\cM_\delta$, with $\mes \cM_\delta\lesssim B_0\,\delta^{\frac14}$, such that for any
$\xi\notin \cM_{\delta}$
\[
\left|\int_{\TT^2}\log\rho\left(f(y) - \xi\right)dy - \int_{\TT^2}\log\left|f(y) - \xi\right|dy\right| \leq \left(\log B_0\right)\delta^{\frac{1}{8}}.
\]
Take $\delta^2 = \tau = N^{-\frac{\sigma}{10}}$, and let $\cT(N) := \cE(N)\cup \cE_{2\delta,\delta}\cup
\cM_\delta$. Then \[\mes\cT(N)\lesssim\exp\left(-\frac{N^\kappa}{2}\right)
+\left(N^{-\frac{\sigma}{10}}\right)^{\frac{1}{4}}
+ B_0\left(N^{-\frac{\sigma}{10}}\right)^{\frac{1}{4}} <
N^{-\frac{\sigma}{50}}.\]
Finally, we conclude from the preceding that   if $\xi\notin\cT(N)$ and $x\notin
B_{N,\xi}$ then
\[
\left|\frac{1}{N}\sum_{k=1}^N\log\left|f\left(T^kx\right) - \xi\right| - \la\log\left|f(\cdot) - \xi\right|\ra\right|
 < N^{-\kappa}
\]
provided $\kappa$ was chosen sufficiently small.

The uniform upper bound~\eqref{eq:uniform} is implicit in the preceding. Indeed, fixing $\cB\subset[1,N]$ with $\#\cB<N^{1-2\kappa}$,
we obtain as above
\begin{align}
&\frac{1}{N}\sum_{k\in[1,N]\setminus\cB}\log\left|f\left(T^kx\right)
- \xi\right| - \la\log\left|f(\cdot) - \xi\right|\ra \nn\\&
= \frac{1}{N}\sum_{k\in[1,N]\setminus\cB}\log\left|f\left(T^kx\right) - \xi\right| -
\frac{1}{N}\sum_{k\in[1,N]\setminus\cB}\log\rho\left(f\left(T^kx\right) - \xi\right) \label{eq:up0}\\
& \hspace*{12pt} + \frac{1}{N}\sum_{k\in[1,N]\setminus\cB}\log\rho\left(f\left(T^kx\right)
 - \xi\right) - \int_{\TT^2}\log\rho\left(f(y) - \xi\right)\,dy \label{eq:up1} \\
& \hspace*{12pt} + \int_{\TT^2}\log\rho\left(f(y) - \xi\right)dy - \int_{\TT^2}\log\left|f(y) - \xi\right|\,dy \label{eq:up2}
\end{align}
It was shown above that for all $\xi\not\in \cT(N)$ we have
$\eqref{eq:up2} \le N^{-\kappa}$
uniformly in $x$. Moreover, with $\phi$ and $\psi$ as in \eqref{eq:smooth_ergodic} above,
\begin{align*}
&\Big|\frac{1}{N}\sum_{k\in[1,N]\setminus\cB}\varphi\left(T^kx\right) - \la\varphi\ra\Big| \\
& \leq \Big|\frac{1}{N}\sum_{k\in[1,N]\setminus\cB}\varphi\left(T^kx\right) -
\frac{1}{N}\sum_{k\in[1,N]\setminus\cB}\psi\left(T^kx\right)\Big|
+ \Big|\frac{1}{N}\sum_{k\in[1,N]\setminus\cB}\psi\left(T^kx\right) - \la\psi\ra\Big|
+ |\la\psi\ra - \la\varphi\ra| \\
& \leq \frac{1}{N}\sum_{k\in[1,N]\setminus\cB} |\varphi(T^kx)-\psi(T^kx)| +
\Big|\frac{1}{N}\sum_{k=1}^N\psi\left(T^kx\right) - \la\psi\ra\Big| + \frac{1}{N}\sum_{k\in\cB}|\psi(T^k x)|
+ |\la\psi\ra - \la\varphi\ra| \\
& \lesssim B_1(f)\delta^{-1}\tau + (|\log\delta/2|+1+B_0(f))\tau^{-4}[N^{-\sigma}+N^{-2\kappa}]
\end{align*}
which implies that $|\eqref{eq:up1}|$ is controlled uniformly in $x$. Finally,
\[ \eqref{eq:up0} = \frac{1}{N} \sum_{k\in[1,N]\setminus\cB} \log \frac{|f(T^kx) - \xi|}{\rho(|f(T^kx) - \xi|)} \le 0\]
where the last inequality follows from the fact that $|y|\le\rho(y)$.
\end{proof}

\begin{remark}
\label{rem:Nga1ga2} In the previous proof we can relax the
Diophantine assumption on $\omega$. Indeed, in the case of shift $Tx
= x + \omega$,  it suffices to require that $\omega$ is
$(N,\gamma_1,\gamma_2)$-Diophantine for some $\gamma_1,\gamma_2>0$
and $N$ the same as in~\eqref{eq:largedt}. This follows from the
fact that the main ergodic theorem for the shift holds under this
weaker Diophantine assumption, see Remark~\ref{rem:weak_dioph}. For
the skew shift $T_\omega(x_1, x_2) = (x_1 + x_2, x_2 + \omega)$ it
suffices to require that $\omega \in \TT_{c,\ve_1, N}$, see
Remark~\ref{rem:A.skewn} in Appendix~A.
\end{remark}

\begin{remark}\label{rem:2.18}
Inspection of the proof of Theorem~\ref{thm:A} shows that the set $\cT(N)$ is a union of two sets $\cT(N)
=\cT'(N)\cup\cT''(N)$ with the following properties:
\begin{itemize}
\item $\mes\cT''(N) \less \exp(-N^{\kappa/2})$
\item $\mes\cT'(N)\less N^{-\kappa}$ and $\cT'(N)$ can be chosen the same for $N^{\frac12}\le N'\le N$. In
particular, the following version of \eqref{eq:uniform} holds:
\[\sup_{N^{\frac12}\le N'\le N}\sup_{\substack{\cB\subset[1,N']\\ \#\cB<(N')^{1-2\kappa}}}\sup_{x\in\TT^2}\frac{1}{N'}
\sum_{k\in[1,N']\setminus\cB}\log\left|f\left(T^kx\right) -
\xi\right| \le \la\log|f(\cdot) - \xi|\ra+ N^{-\kappa}
\]
for any
$\xi\in[-B_0(f),B_0(f)]\setminus\cT(N)$.
\end{itemize}
Moreover, invoking Remark~\ref{rem:2.8FAT} yields
\[
\mes[\cT(N)](r)\less N^{-\kappa/2}
\]
where $r=\exp(-N^\kappa)$.
\end{remark}

\begin{remark}\label{rem:2.ldtsect} The set of exceptional phases $x \in \TT^2$ in Theorem~\ref{thm:A}
 derives only from Lemma~\ref{lem:2.6}, since all other estimates are uniform in $x \in \TT^2$.
 Taking into account Remark~\ref{rem:2.disect} (for Lemma~\ref{lem:2.6}) one obtains the following version
 of  the first statement of Theorem~\ref{thm:A} (which we will use in Section~\ref{sec:3} for the case of the skew shift):
 for any $x_2\in\tor$ there exists $\cT^{(1)}(x_2, N)$ such that
\be\label{eq:2.ldtsect} \mes\left\{x_1 \in \TT: \Big|{1\over
N}\sum^N_{k=1} \log \big|f\bigl(T^k(x_1, x_2)\bigr) - \xi\big| -
\langle\log\big|f(\cdot)-\xi\big|\rangle\Big|> N^{-\kappa}\right\}
\le \exp(-N^\kappa) \ee provided $\xi \in
\cJ_0(f)\setminus\cT^{(1)}(x_2, N)$, where $\mes\cT^{(1)}(x_2, N)
\le N^{-\kappa}$.
\end{remark}

The method of proof of Theorem~\ref{thm:A} is quite robust and
applies to other dynamics as well.  For our applications of
Theorem~\ref{thm:A} to localization, we need the following
modifications involving functions that depend also on $\omega$. Let
$f\in C^1\left(\TT^2\times\TT^2\right)$ and write $T_{\omega} :
\TT^2 \rightarrow \TT^2$ to indicate the dependence on~$\omega$. As
before, we define $S_{f(\cdot,\omega)}(\xi,\delta) = \left\{
x\in\TT^2 : \left|f(x,\omega) - \xi\right| < \delta\right\}$. In
analogy with Corollary~\ref{cor:2.cprime} we now have the following
result.
\begin{corollary}
\label{cor:2'.torus} Let $f\in C^1\left(\TT^2\times\TT^2\right)$.
\begin{itemize}
\item Let $T_\omega: \TT^2 \to \TT^2$ be the shift.  Assume that
$\omega_0$ is $(N,\gamma_1,\gamma_2)$--Diophantine. Then for any
$\xi\in\IR$, any small $\delta > 0$, and
$N>\delta^{-\frac{20}{\sigma}}$, as well as
$|\omega-\omega_0|<N^{-1}$, one has \be\label{eq:2.count}
\frac{1}{N}\#\left\{ 1\leq k\leq N : T_{\omega}^kx\in
S_{f(\cdot,\omega)}(\xi,\delta)\right\}\lesssim\mes
S_{f(\cdot,\omega)}(\xi,2\delta) + \left(1 +
B_1(f)\right)\delta^{\frac{1}{2}} \ee \item Let $T_\omega: \TT^2 \to
\TT^2$ be the skew-shift. Assume $\omega_0 \in \TT_{c,\ve_1, N}$.
Then \eqref{eq:2.count} is valid, provided $N >
\delta^{-\frac{20}{\sigma}}$, $\delta$ is small, and $|\omega -
\omega_0| < N^{-3}$.
\end{itemize}
\end{corollary}

\begin{proof} The proof is basically the same as that of Corollary~\ref{cor:2.cprime}. More precisely,
with
$\chi_{\delta}$ as in Lemma~\ref{lem:2.list},
\begin{itemize}
\item $\#\left\{ 1\leq k\leq N : T_{\omega}^kx\in S_{f(\cdot,\omega)}(\xi,\delta)\right\}\leq\sum_{k=1}^N\chi_{\delta}\left(f\left(T_{\omega}^kx,\omega\right) - \xi\right)$
\item $0\leq\la\chi_{\delta}\left(f(\cdot,\omega) - \xi\right)\ra\leq\mes S_{f(\cdot,\omega)}(\xi,2\delta)$
\end{itemize}
Here $\delta>0$ is any small number, $\xi$ an arbitrary real number,
$\omega\in\TT^2$, and $N$ a positive integer. Now one proceeds as in
Corollary~\ref{cor:2.cprime} using the ergodic theorem, i.e.,
Propositions~\ref{prop:A.ineq} and~\ref{prop:A.skewineq} from
Appendix~A. The point to notice here is that the constants in the
ergodic theorem are uniform in $|\omega-\omega_0|< N^{-1}$.
\end{proof}

We can again remove a set of exceptional $\xi$ for which the measure $\mes S_{f(\cdot,\omega)}(\xi,2\delta)$
 is too large;  as in Corollary~\ref{cor:2.xiset} this is an easy consequence of Fubini's theorem with the added feature
that the set we remove can be chosen to be the same for all
$\omega$ close to a given $\omega_0$.

\begin{lemma}\label{lem:stableS}
Let $\omega_0\in\TT^2$, $\xi\in\IR$, $ \delta>0$. Then
\[ S_{f(\cdot,\omega)}(\xi,\delta)\subset S_{f(\cdot,\omega_0)}(\xi,2\delta) \]
for all $|\omega-\omega_0|<B_1(f)^{-1}\delta$. In particular, there exists
$\cE_{\delta,\omega_0}\subset\IR$,
$\mes\cE_{\delta,\omega_0}\lesssim\delta^{\frac{1}{2}}$ such that
 for any $\xi\notin{\cE}_{\delta,\omega_0}$, one has
 \[
\mes S_{f(\cdot,\omega)}(\xi,\delta) \lesssim \delta^{\frac12}
\]
provided $|\omega-\omega_0|<B_1(f)^{-1}\delta$.
\end{lemma}
\begin{proof}
Clearly,
$\left|f(x,\omega) - f(x,\omega_0)\right|\leq B_1(f)|\omega -
\omega_0| < \delta$.  Thus, if $|f(x,\omega) - \xi| < \delta$ then also
$|f(x,\omega_0) - \xi|\leq\left|f(x,\omega_0) - f(x,\omega)\right| +
\left|f(x,\omega) - \xi\right|  < 2\delta$
and the lemma follows. The second statement follows from
\[ \int \mes S_{f(\cdot,\omega)}(\xi,\delta)\, d\xi = 2\delta \]
and Fubini's theorem.
\end{proof}

Combining the previous two statements yields the following:

\begin{corollary}
\label{cor:2.DC} Given $\delta > 0$, let
$\cE_{\delta,\omega_0}\subset\IR$ be as in the previous lemma.
\begin{itemize}
\item Let  $T_\omega: \TT^2 \to \TT^2$ be the shift.  Assume $\omega_0$
is $(N,\gamma_1,\gamma_2)$--Diophantine. Then
 for any $\xi\notin{\cE}_{\delta,\omega_0}$, one has
\be\label{eq:2.newcount} \frac{1}{N}\#\left\{ 1\leq k\leq N :
T_{\omega}^kx\in
S_{f(\cdot,\omega)}(\xi,\delta)\right\}\lesssim\left(1 +
B_1(f)\right)\delta^{\frac{1}{2}} \ee for all $x\in\TT^2$ and
$\left|\omega - \omega_0\right| < ((1+B_1(f))N)^{-1}$ provided
$N\geq\delta^{-\frac{20}{\sigma}}$ with $\sigma>0$ a
sufficiently small constant depending on~$\omega_0$.

\item Let $T_\omega: \TT^2 \to \TT^2$ be the skew-shift.  Assume
$\omega_0 \in \TT_{c,\ve_1,N}$.  Then \eqref{eq:2.newcount} is valid
provided $|\omega- \omega_0| < (1 + B_1(f))^{-1}N^{-3}$, $ N >
\delta^{-{\frac{20}{\sigma}}}$ with some small $\sigma>0$.
\end{itemize}
\end{corollary}

\begin{proof}
We can apply Corollary~\ref{cor:2'.torus} for large $N$, since
$\left|\omega - \omega_0\right| \le \frac{c}{2} N^{-1}$ for the
shift, and $\left|\omega - \omega_0\right| \le \frac{c}{2}N^{-3}$
for the skew-shift. Furthermore, since
 $|\omega - \omega_0| < \delta/ B_1(f)$ we can apply the previous lemma to conclude that
\[
\mes S_{f(\cdot,\omega)}(\xi,3\delta) \lesssim \delta^{\frac12}
\]
for all $\xi\notin{\cE}_{\delta,\omega_0}$.
\end{proof}

\begin{corollary}
Using the notation of Corollary~\ref{cor:2.DC},  one has
\begin{align*}
\frac{1}{N}\sum_{k=1}^N\log\rho\left(f\left(T_{\omega}^kx,\omega\right) - \xi\right) & \leq \frac{1}{N}\hspace*{-20pt}\sum_{\substack{1\leq k\leq N\\ T_{\omega}^kx\notin S_{f(\cdot,\omega)}(\xi,\delta)}}\hspace*{-22pt}\log\left|f\left(T_{\omega}^kx,\omega\right) - \xi\right| \\
& \leq \frac{1}{N}\sum_{k=1}^N\log\rho\left(f\left(T_{\omega}^kx,\omega\right) - \xi\right) + \left(1 + B_1(f)\right)N^{-\frac{\sigma}{40}}
\end{align*}
for any $\xi\notin\cE_{\delta,\omega_0}$
 where $\delta=N^{-\frac{\sigma}{20}}$.
\end{corollary}

We now present a somewhat sharper version of Lemma~\ref{2:lem.capL} on large values
of certain logarithmic integrals.

\begin{lemma}
\label{lem:2.jensen}
Given $R > 0$, $\omega_0\in\TT^2$, $\eta > 0$, there exists
\[
\cL_{\omega_0,\eta,R}\subset(\omega_0 - \eta, \omega_0 + \eta)\times
\cJ_0(f), \quad \mes\cL_{\omega_0,\eta,R}\lesssim (1+B_0(f))\eta
\exp(-\sqrt{R}/2)
\] such that for $(\omega,\xi)\in\TT^2\times(\omega_0 - \eta,
\omega_0 + \eta)\times \cJ_0(f)\setminus\cL_{\omega_0,\eta,R}$ one
has
\begin{itemize}
\item[(1)]$\int_{\TT^2}\bigl|\log\left|f(x,\omega) - \xi\right|\bigr|^2\,dx\leq  R$
\item[(2)]
$\int_{S}\bigl|\log|f(x,\omega) -
\xi|\bigr|\,dx\leq R^{\frac{1}{2}}(\mes
S)^{\frac{1}{2}}$
\end{itemize}
where $S$ is an arbitrary measurable set in part (2).
Moreover, an analogous statement holds with $\omega=\omega_0$ fixed. In that case we only need to remove sets of~$\xi$.
\end{lemma}

\begin{proof}
The function $\Phi(y)=\exp(\sqrt{1+y})$ is convex on $y>0$. Then, by Jensen inequality
\begin{align}
& \mes \Big\{ (\omega,\xi)\::\: \int_{\TT^2} \Big|\log|f(x,\omega)-\xi|\Big|^2\,dx  >  R \Big\} \nn\\
& =\mes \Big\{ (\omega,\xi)\::\: \Phi\Big(\frac14\int_{\TT^2} \Big|\log|f(x,\omega)-\xi|\Big|^2\,dx \Big) > \Phi(\frac14 R) \Big\} \nn\\
& \le \mes \Big\{ (\omega,\xi)\::\: \int_{\TT^2} \Phi\Big(\frac14\Big|\log|f(x,\omega)-\xi|\Big|^2\Big)\,dx  > \Phi(\frac14 R) \Big\} \nn \\
& \le [\Phi(\frac14 R)]^{-1}  \int_{\TT^2}\int_\cJ\int_{\TT^2} \Phi\Big(\frac14\Big|\log|f(x,\omega)-\xi|\Big|^2\Big)\,dx\,d\xi\, d\omega \label{eq:Phi}
\end{align}
Note that $\Phi(y^2)\le \exp(y+1)$ for all $y\ge0$. Hence,
\begin{align*}
& \Phi\Big(\frac14\Big|\log|f(x,\omega)-\xi|\Big|^2\Big) \le  \exp\Big(\frac12 \Big|\log|f(x,\omega)-\xi|\Big|+1\Big) \\
& \lesssim (1+B_0(f))^{\frac12} |f(x,\omega)-\xi|^{-\frac12}
\end{align*}
Inserting this into \eqref{eq:Phi} yields
\begin{align*}
\mes \Big\{ (\omega,\xi)\::\: \int_{\TT^2} \Big|\log|f(x,\omega)-\xi|\Big|^2\,dx  >  R \Big\} &\lesssim e^{-\frac12 \sqrt{R}} (1+B_0(f))^{\frac12}  \int_{\TT^2}\int_\cJ\int_{\TT^2} |f(x,\omega)-\xi|^{-\frac12}\,dx d\xi\,d\omega \\
&\lesssim e^{-\frac12 \sqrt{R}} (1+B_0(f)) \eta,
\end{align*}
which proves (1). Finally,
claim (2) follows from (1) by Cauchy-Schwarz.
The final statement of the lemma follows by the same arguments but without averaging in~$\omega$.
\end{proof}

We can now formulate the analogue of Corollary~\ref{cor:2.11} for the case of functions which
depend on $\omega$, but with exceptional sets that do not depend on $\omega$ as long as $|\omega-\omega_0|$
is sufficiently small.

\begin{corollary}\label{cor:EM}
Given $N$, let $\delta=N^{-\frac{\sigma}{20}}$ and
$\eta=\frac{c}{2}(1+B_1(f))^{-1}N^{-1}$. For any $\omega_0\in\TT^2$
let $\cE_{\delta,\omega_0}$ be as in Lemma~\ref{lem:stableS}. There
exists $\cM_{\delta,\omega_0}\subset(\omega_0 - \eta,\omega_0 +
\eta)\times \cJ_0(f)$, $\mes\cM_{\delta,\omega_0}\leq
\exp(-N^{\sigma_1})$, where $\sigma_1>0$ is some small constant,
such that for any $(\omega,\xi)\in(\omega_0 - \eta,\omega_0 +
\eta)\times(\cJ_0(f)\setminus\cE_{\delta,\omega_0})\setminus\cM_{\delta,\omega_0}$
one has
\[
\left|\int_{\TT^2}\log\rho(f(x,\omega) - \xi)\,dx - \int_{\TT^2}\log|f(x,\omega) - \xi|\,dx\right|\lesssim \delta^{\frac{1}{8}}.
\]
\end{corollary}

\begin{proof}
Suppose $\xi\notin\cE_{\delta,\omega_0}$,
$(\omega,\xi)\notin\cL_{\omega_0,\eta,R}$ (see Lemma~\ref{lem:2.jensen}). By
Corollary~\ref{cor:2.DC},
\[\mes
S_{f(\cdot,\omega)}(\xi,\delta)\leq\delta^{\frac{1}{2}}\qquad \forall\;
\omega\in(\omega_0 - \eta,\omega_0 + \eta).
\]
Hence
\begin{align*}
&\left|\int_{\TT^2}\log\rho(f(x,\omega) - \xi)\,dx -
\int_{\TT^2}\log|f(x,\omega) - \xi|\,dx\right| \\&
\leq \int_{S_{f(\cdot,\omega)}(\xi,\delta)}\bigl|\log\rho(f(x,\omega) - \xi)\bigr|\,dx + \int_{S_{f(\cdot,\omega)}(\xi,\delta)}\bigl|\log|f(x,\omega) - \xi|\bigr|\,dx \\
& \leq \left[\mes S_{f(\cdot,\omega)}(\xi,\delta)\right]\left|\log\frac{\delta}{2}\right| +
R^{\frac{1}{2}}\left[\mes S_{f(\cdot,\omega)}(\xi,\delta)\right]^{\frac{1}{2}} \\
& \leq \delta^{\frac{1}{2}}\left|\log\frac{\delta}{2}\right| + R^{\frac{1}{2}}\delta^{\frac{1}{4}}
\end{align*}
Take $R = \delta^{-\frac{1}{4}}$, $\cM_{\delta,\omega_0} :=
\cL_{\omega_0,\eta,R}$.
\end{proof}

We are now ready to prove the analogue of Theorem~\ref{thm:A} for functions depending on~$\omega$.
The reader should take note of the fact that we first remove a large (i.e., of size $N^{-\kappa}$)
set of exceptional parameters $\xi$ which only depends on $\omega_0$ -- after that we proceed to
remove {\em exponentially small} sets in $(x,\omega,\xi)$.  In the following theorem we use the notion
of $(N,\gamma_1,\gamma_2)$--Diophantine $\omega$, cf.~Remark~\ref{rem:Nga1ga2} and Remark~\ref{rem:weak_dioph}.

\begin{theorem}
\label{thm:2.21}
Let $f(x,\omega)$ be $C^1$--smooth.  Let $T_\omega: \TT^2 \to \TT^2$ be a shift (a skew-shift).  Given large $N$ assume that $\omega$ is $(N,\gamma_1,\gamma_2)$--Diophantine $(\omega \in \TT_{c,\ve_1, N})$ for some small $\gamma_1, \gamma_2 > 0$ (for some small $\ve_1 > 0)$.
 Then there exists
$\cT(N)\subset\cJ_0(f)$, $\mes\cT(N) < N^{-\gamma}$, such that
\begin{multline*}
\mes\Biggl\{(x,\omega,\xi)\in\TT^2\times(\omega_0 - \eta,\omega_0 + \eta)\times\bigl(\cJ\setminus\cT(N)\bigr) : \Biggr. \\
\Biggl.\left|\frac{1}{N}\sum_{k=1}^N\log\left|f\left(T_{\omega}^kx,\omega\right) - \xi\right| - \la\log|f(\cdot,\omega) - \xi|\ra\right| > N^{-\gamma}\Biggr\}
\leq\exp\left(-N^{\gamma}\right)
\end{multline*}
provided $\eta=(1+B_1(f))^{-1}N^{-3}$. Here $\gamma>0$ is a small constant that depends on the
Diophantine condition.
\end{theorem}

\begin{proof} We proceed as in the proof of Theorem~\ref{thm:A}. Thus, let $\rho$ be as in Definition~\ref{def:rho} with $\delta$ to be specified later. Then
\begin{align*}
&\left|\frac{1}{N}\sum_{k=1}^N\log\left|f\left(T_\omega^kx,\omega\right) - \xi\right| - \la\log\left|f(\cdot,\omega) - \xi\right|\ra\right| \\
& \leq \left|\frac{1}{N}\sum_{k=1}^N\log\left|f\left(T_\omega^kx,\omega\right) - \xi\right| - \frac{1}{N}\sum_{k=1}^N\log\rho\left(f\left(T_\omega^kx,\omega\right) - \xi\right)\right| \\
& \hspace*{12pt} + \left|\frac{1}{N}\sum_{k=1}^N\log\rho\left(f\left(T_\omega^kx,\omega\right) - \xi\right) - \int_{\TT^2}\log\rho\left(f(y,\omega) - \xi\right)\,dy\right| \\
& \hspace*{12pt} + \left|\int_{\TT^2}\log\rho\left(f(y,\omega) - \xi\right)\,dy - \int_{\TT^2}\log\left|f(y,\omega) - \xi\right|\,dy\right|
\end{align*}
Let $B_{N,\xi,\omega} = \left\{ y\in\TT^2 : \min_{1\leq k\leq
N}\left|f\left(T_\omega^ky,\omega\right) - \xi\right| \leq e^{-N^\kappa}\right\}$ where $\kappa>0$ is small.
Moreover, let $\cE(N,\omega)$ be as in Lemma~\ref{lem:2.6} applied to $f(\cdot,\omega)$ and $\cE_{\delta, \omega_0}$ be as in Lemma~\ref{lem:stableS}.
If $\xi\notin\cE(N,\omega)\cup\cE_{\delta,\omega_0}$, and $x\notin B_{N,\xi,\omega}$, then
\begin{align*}
&\left|\frac{1}{N}\sum_{k=1}^N\log\left|f\left(T_\omega^kx,\omega\right) - \xi\right| -
\frac{1}{N}\sum_{k=1}^N\log\rho\left(f\left(T_\omega^kx,\omega\right) - \xi\right)\right| \\
& \leq \frac{1}{N}\#\left\{ 1\leq k\leq N : T_\omega^kx\in S_f(\xi,\delta)\right\}
 \left[\left|\log\min_{1\leq k\leq N}\left|f\left(T_\omega^kx,\omega\right) - \xi\right|\right| + \left|\log\frac{\delta}{2}\right|\right] \\
& \lesssim \left(1 + B_1(f)\right)\delta^{\frac{1}{2}}\left[N^\kappa + \left|\log\frac{\delta}{2}\right|\right]
\end{align*}
Moreover, $\mes B_{N,\xi,\omega} \leq \exp\left(-\frac{1}{4}N^\kappa\right)$.
Let $\varphi(y) = \log\rho\left(f(y,\omega) - \xi\right)$. By Lemma~\ref{lem:2.C4}, for any $\tau >
0$, there is $\psi\in C^4\left(\TT^2\right)$ such that
\begin{itemize}
\item[(1)]$\max_{y\in\TT^2}\left|\varphi(y) - \psi(y)\right|\lesssim B_1(\varphi)\tau\leq B_1(f)\delta^{-1}\tau$
\item[(2)]$B_4(\psi)\lesssim B_0(\varphi)\tau^{-4}\leq(|\log(\delta/2)|+1+B_0(f))\tau^{-4}$
\end{itemize}
Then with some $\sigma>0$,
\begin{align*}
\left|\frac{1}{N}\sum_{k=1}^N\varphi\left(T_\omega^kx\right) - \la\varphi\ra\right| & \leq \left|\frac{1}{N}\sum_{k=1}^N\varphi\left(T_\omega^kx\right)
- \frac{1}{N}\sum_{k=1}^N\psi\left(T_\omega^kx\right)\right| + \left|\frac{1}{N}\sum_{k=1}^N\psi\left(T_\omega^kx\right) -
\la\psi\ra\right| + \left|\la\psi\ra - \la\varphi\ra\right| \\
& \lesssim B_1(f)\delta^{-1}\tau + (|\log\delta/2|+1+B_0(f))\tau^{-4}N^{-\sigma}
\end{align*}
provided $N$ is sufficiently large (see Proposition~\ref{prop:2.largen}).
By Corollary~\ref{cor:EM} there exists $\cM_{\delta,\omega_0}$, with $\mes \cM_{\delta,\omega_0}\lesssim \exp(-N^{\sigma_1})$, such that for any
$(\omega,\xi)\notin \cM_{\delta,\omega_0}$ one has
\[
\left|\int_{\TT^2}\log\rho\left(f(y,\omega) - \xi\right)\,dy - \int_{\TT^2}\log\left|f(y,\omega) - \xi\right|\,dy\right| \leq \delta^{\frac{1}{8}}.
\]
Take $\delta^2 = \tau = N^{-\frac{\sigma}{10}}$, and let $\cT(N) := \cE_{\delta,\omega_0}$.
Then with some small $\gamma>0$, \[\mes\cT(N)\lesssim \delta^{\frac12} <
N^{-\gamma}.\]
Finally, we conclude from the preceding that   if $\xi\notin\cT(N)\cup\cE(N,\omega)$, $(\omega,\xi)\not\in \cM_{\delta,\omega_0}$,  and $x\notin
B_{N,\xi,\omega}$ then
\[
\left|\frac{1}{N}\sum_{k=1}^N\log\left|f\left(T_\omega^kx\right) - \xi\right| - \la\log\left|f(\cdot) - \xi\right|\ra\right|
 < N^{-\gamma}
\]
provided $\gamma$ was chosen sufficiently small.
\end{proof}

\section{Large Deviation Theorems in Frequencies and Elimination of resonances in a general setting}
\label{sec:3}

Let $f\in C^1\left(\TT^{2}\times\TT^2\right)$.  Let $T_\omega:
\TT^2 \to \TT^2$ be the shift (the skew-shift). We begin this
section with some simple statements concerning the introduction of
perturbations into the results of the previous section.

\begin{lemma}\label{lem:3.1}
For any $x,\ve_k,\omega\in\TT^{2}$, $\tau_k,\omega_1\in\TT^2$ $(x,\ve_k \in \TT^2$, $\tau_k, \omega, \omega_1 \in \TT$) one has
\begin{align*}
&\#\left\{ 1\leq k\leq N : \left|f\left(T_{\omega}^kx + \ve_k,\tau_k+\omega_1\right) - \xi\right| < \delta\right\} \\
&\leq \#\left\{ 1\leq k\leq N : \left|f\left(T_{\omega}^k x,\omega_1\right) - \xi\right| < \delta + B_1(f)\max_k\left\{|\ve_k| + |\tau_k|\right\}\right\}
\end{align*}
In particular, if $\ve := \max_k\left\{|\ve_k| + |\tau_k|\right\} <
\frac{\delta}{B_1(f)}$ then
\[
\#\left\{ 1\leq k\leq N : \left|f\left(T_\omega^kx + \ve_k,\omega_1+\tau_k\right)-\xi\right| < \delta\right\}
\leq\#\left\{ 1\leq k\leq N : T_\omega^kx\in S_{f(\cdot,\omega_1)}(\xi,2\delta)\right\}
\]
\end{lemma}

\begin{corollary}\label{cor:3.2} Let $N$ be large and assume that
$\omega_0$ is $(N,\gamma_1,\gamma_2)$-Diophantine for the case of
the shift, see Remark~\ref{rem:weak_dioph}
 (or $\omega_0 \in \TT_{c,\ve_1, N}$ for the skew-shift, see Remark~\ref{rem:A.skewn} in Appendix~A).
  Given $\delta \ge N^{-\frac{\sigma}{20}}$, assume that $\ve := \max_k\left\{|\ve_k| + |\tau_k|\right\}
  < \frac{\delta}{B_1(f)}$. Then, with $\omega_1$ fixed,
\begin{itemize}
\item[(1)] the estimate
\[\frac{1}{N}\#\left\{ 1\leq k\leq N : \left|f\left(T_\omega^kx + \ve_k,\omega_1+\tau_k\right) - \xi\right| < \delta\right\}
\lesssim\mes S_{f(\cdot,\omega_1)}(\xi,3\delta) + \left[1 +
B_1(f)\right]\delta^{\frac{1}{2}}\] holds for all $x\in\TT^2$,
provided $|\omega-\omega_0|<N^{-3}$, $\xi\in\cJ_0(f)$
\item[(2)]
there exists $\cE_{\omega_0,\omega_1,\delta}\subset\IR$,
$\mes\cE_{\omega_0,\omega_1,\delta}\lesssim\delta^{\frac{1}{2}}$
such that for any
$\xi\in\cJ_0(f)\setminus\cE_{\omega_0,\omega_1,\delta}$,
$|\omega-\omega_0|<N^{-3}$ one has
\[
\frac{1}{N}\#\left\{ 1\leq k\leq N \::\: \left|f\left(T_\omega^kx + \ve_k,\omega_1+\tau_k\right) - \xi\right|
< \delta\right\}\lesssim\left[1 + B_1(f)\right]\delta^{\frac{1}{2}}
\]
for all $x\in\TT^2$
($\cE_{\omega_0,\omega_1,\delta}$ does not depend on $\omega,\ve_k$, $\tau_k$).
Furthermore, if $|\omega_0-\omega_1|<\frac{\delta}{B_1(f)}$, then
$\cE_{\omega_0,\omega_1,\delta}$ can be chosen to depend only on $\omega_0,\delta$.
\end{itemize}
\end{corollary}

\begin{proof}
Recall that
\[
\frac{1}{N}\#\left\{ 1\leq k\leq N : T_\omega^k x \in S_{f(\cdot,\omega_1)}(\xi,2\delta)\right\}
\lesssim\mes S_{f(\cdot,\omega_1)}(\xi,3\delta) + \left[1 + B_1(f)\right]\delta^{\frac{1}{2}}
\]
for any $|\omega-\omega_0|<N^{-2}$ and any $x\in\TT^2$ due to Corollary~\ref{cor:2'.torus}.
Therefore, (1) follows from Lemma~\ref{lem:3.1}. Assertion (2) is a consequence of Lemma~\ref{lem:stableS}.
\end{proof}

\begin{remark}\label{rem:3.fat}
As we have noted in Remark~\ref{rem:2.8FAT}, the estimates of
Corollary~\ref{cor:2.xiset} can be stated in a slightly stronger
form which we need in our applications. Namely, the set $\cE_\delta$
in that corollary satisfies
\[ \mes[\cE_\delta](\rho)\less \delta^{\frac12},\quad \rho\le\delta \]
where $[\cE_\delta](\rho)=\{\xi\::\:\dist(\xi,\cE_\delta)\}\le\rho$. For the same reason the set $\cE_{\omega_0,\omega_1,\delta}$ in Corollary~\ref{cor:3.2} obeys
\[ \mes[\cE_{\omega_0,\omega_1,\delta}](\rho)\less \delta^{\frac12},\quad \rho\le\delta \]
\end{remark}

\begin{lemma}
\label{lem:3.ve}
Assume $\ve < \frac{1}{2}\left(\frac{\delta}{B_1(f)}\right)$. If $\left|f\left(T_\omega^kx + \ve_k,\omega_1+\tau_k\right) - \xi\right| > \delta$ then
\[
\Bigl|\log\left|f\left(T_\omega^kx + \ve_k,\omega_1+\tau_k\right) - \xi\right| - \log\left|f\left(T_\omega^kx,\omega_1\right) - \xi\right|\Bigr|\lesssim B_1(f)\,\ve\,\delta^{-1}
\]
uniformly in $\omega,\omega_1$.
\end{lemma}

\begin{proof}
\[
1 - \frac{B_1(f)\,\ve}{\delta}\leq\left|\frac{f\left(T_\omega^kx,\omega_1\right) - \xi}{f\left(T_\omega^kx + \ve_k,\omega_1+\tau_k\right) - \xi}\right|\leq 1 + \frac{B_1(f)\,\ve}{\delta}.
\]
By assumption, $B_1(f)\,\ve\,\delta^{-1} < \frac{1}{2}$ and the assertion follows.
\end{proof}

We shall also need the following analogue of \eqref{eq:oldstar}:
given $x\in\TT^2$, $|\ve_k|, |\tau_k| \ll 1$, let
\begin{equation}\label{eq:3.SETJ}
J_N(x,\xi,\delta) = \left\{ 1\leq k\leq N : \left|f\left(T_\omega^kx +
      \ve_k,\omega_1+\tau_k\right) - \xi\right| < \delta\right\}.
\end{equation}
 Then
\begin{equation}
\label{eq:star}
\begin{aligned}
\frac{1}{N}\hspace*{-8pt}\sum_{\substack{1\leq k\leq N\\ k\notin J_N(x,\xi,\delta)}}\hspace*{-11pt}\log\left|f\left(T_\omega^kx + \ve_k,\omega_1+\tau_k\right) - \xi\right| & +
\frac{1}{N}\bigl[\# J_N(x,\xi,\delta)\bigr]\log\left(\min_{1\leq k\leq N}\left|f\left(T_\omega^kx + \ve_k,\omega_1+\tau_k\right) - \xi\right|\right) \\
& \hspace*{9pt} \leq \frac{1}{N}\sum_{k=1}^N\log\left|f\left(T_\omega^kx + \ve_k,\omega_1+\tau_k\right) - \xi\right| \\
& \hspace*{9pt} \leq \frac{1}{N}\hspace*{-8pt}\sum_{\substack{1\leq k\leq N\\ k\notin J_N(x,\xi,\delta)}}\hspace*{-11pt}\log\left|f\left(T_\omega^kx + \ve_k,\omega_1+\tau_k\right) - \xi\right|. \\
\end{aligned}
\end{equation}

As in Lemma~\ref{lem:2.6}, Fubini's theorem immediately yields the
following statement (recall the definition of $\cJ_0(f)$
in~\eqref{eq:J_def}):

\begin{lemma}\label{lem:3.4}
Given $N\in\IN$, ${\ve}_k, \tau_k:\tor^2\to\tor^2$,
 $k =
1,2,\ldots,N$,  $\omega_1$ and any small $\kappa>0$ there exists
$\cE_{\omega_1}\left(N,\{{\ve}_k\},\{{\tau}_k\}\right)\subset \cJ_0(f)$, 
$\mes\cE_{\omega_1}\left(N,\{{\ve}_k\},\{{\tau}_k\}\right)\leq\exp\left(-\frac{1}{2}N^{\kappa}\right)$
such that for any
$\xi\in\cJ_0(f)\setminus\cE_{\omega_1}\left(N,\{{\ve}_k\},\{{\tau}_k\}\right)$
one has
\begin{equation}\label{eq:3.LOG}
\mes\left\{ x\in\TT^2 \::\: \log\left(\min_{1\leq k\leq N}\left|f\left(T_{\omega_1}^kx + {\ve}_k(x),\omega_1+{\tau}_k(x)\right) - \xi\right|\right) < - N^{\kappa}\right\} <
\exp\left(-\frac{1}{4}N^{\kappa}\right).
\end{equation}
Furthermore, just as in Remark~\ref{rem:3.fat} the set $\cE_{\omega_1}\left(N,\{{\ve}_k\},\{{\tau}_k\}\right)$
 satisfies
\begin{equation}
\label{eq:fett} \mes[\cE_{\omega_1}\left(N,\{{\ve}_k\},\{{\tau}_k\}\right)](\rho)\le \exp(-N^\kappa/4)
\end{equation}
for any $\rho\le \exp(-N^\kappa)$.
\end{lemma}

The following result is a perturbed version of Theorem~\ref{thm:A}. Note that in the statement of the following
theorem we introduce two different sets of $\xi$ which need to be removed. This is due to the fact that in later
applications we wish to sum over the perturbations $\ve_k$ and $\tau_k$.

\begin{prop}\label{prop:A'}
Let $f(x,\omega)$ be $C^1$--smooth.  Let $T_\omega: \TT^2 \to \TT^2$
be the shift (or the skew-shift).  Let $N$ be large.  Assume that
$\omega_0$ is $(N,\gamma_1,\gamma_2)$--Diophantine (or $\omega_0 \in
\TT_{c,\ve_1,N}$ for the skew-shift).  Let \[\ve_k(x,\omega),
\tau_k(x,\omega),\ (x,\omega) \in \TT^2 \times \TT^2\ ((x,\omega)\in
\TT^2\times \TT),\] obey
$\max_kB_0\left({\ve}_k\right)\lesssim\frac{N^{-1 }}{B_1(f)}$,
$\max_kB_0\left({\tau}_k\right)\lesssim\frac{N^{-1 }}{B_1(f)}$.
Moreover, let $|\omega_0-\omega_1|<[(1+B_1(f))N^3]^{-1}$. Then there
exist $\cE_{\omega_0}(N)$,
$\cE_{\omega_0,\omega_1}\left(N,\{{\ve}_k\},\{{\tau}_k\}\right)\subset\cJ_0(f)$,
such that
\[\mes\cE_{\omega_0}(N)\lesssim N^{-\kappa},\quad
\mes\cE_{\omega_0,\omega_1}\left(N,\{{\ve}_k\},\{{\tau}_k\}\right) <
\exp\left(-\frac{1}{2}N^{\kappa}\right),\]
and  so that for any
\[\xi\in\cJ_0(f)\setminus\left(\cE_{\omega_0}(N)\cup\cE_{\omega_0,\omega_1}(N,\{{\ve}_k\},\{{\tau}_k\})\right)\]
one has
\begin{align*}
&\mes\Big\{ x\in\TT^2\: :\: \Big|\frac{1}{N}\sum_{k=1}^N\log\Bigl|f\left(T_{\omega_1}^kx + {\ve}_k(x),\omega_1+{\tau}_k(x)\right) - \xi\Bigr|
- \la\log\left|f(\cdot,\omega_1) - \xi\right|\ra\Big| > N^{-\kappa}\Big\}\\
&\leq\exp\left(-N^{\kappa}\right)
\end{align*}
Here $\kappa>0$ is some small constant.
Moreover, the sets $\cE_{\omega_0}(N)$, $\cE_{\omega_0,\omega_1}(N,\{{\ve}_k\},\{{\tau}_k\})$ obey
\begin{equation}
\label{eq:fattened}
\mes [\cE_{\omega_0}(N)](\rho) \less N^{-\kappa}, \quad \mes [ \cE_{\omega_0,\omega_1}(N,\{{\ve}_k\},\{{\tau}_k\}) ](\rho) \less \exp(-N^\kappa/2)
\end{equation}
for any $\rho<\exp(-N^\kappa)$.
\end{prop}
\begin{proof}
We shall reduce this theorem to Theorem~\ref{thm:A} applied to the function
$f(\cdot,\omega_1)$ (see also Remark~\ref{rem:Nga1ga2}). Let $J_N(x,\xi,\delta)$ be as in~\eqref{eq:3.SETJ}. Set
$\delta=N^{-\frac{\sigma}{20}}$ with $0<\sigma<1$. Due to
Corollary~\ref{cor:3.2}, for all $x\in\TT^2$, $\xi \notin \cE_{\omega_0,\omega_1,\delta}$
\[
N^{-1}\# J_N(x,\xi,\delta) \lesssim (1+B_1(f)) N^{-\frac{\sigma}{40}}
\]
Let $\cE_{\omega_1}(N,\{\eps_k\},\{\tau_k\})$ be as in
Lemma~\ref{lem:3.4} and let $\cB_{\omega_1}(N,\xi)$ be the set
defined in~\eqref{eq:3.LOG}. Then for any
$\xi\in\cJ_0(f)\setminus\cE_{\omega_1}(N,\{\eps_k\},\{\tau_k\})$ and
any $x\in\TT^2\setminus \cB_{\omega_1}(N,\xi)$ one has due to
Lemma~\ref{lem:3.1}
\[ \Big|\frac{1}{N}\sum_{k=1}^N\log\bigl|f\left(T_{\omega_1}^kx +
    {\ve}_k(x),\omega_1+{\tau}_k(x)\right) - \xi\bigr|
-\frac{1}{N}\sum_{\substack{1\le k\le N\\k\not\in J_N(x,\xi,\delta)}}\log\bigl|f\left(T_{\omega_1}^kx + {\ve}_k(x),\omega_1+{\tau}_k(x)\right) - \xi\bigr|
\Big| \lesssim N^{-\frac{\sigma}{50}}
\]
By Lemma~\ref{lem:3.ve},
\[
\frac{1}{N}\sum_{\substack{1\le k\le N\\k\not\in J_N(x,\xi,\delta)}}\Big|\log\bigl|f\left(T_{\omega_1}^kx +
    {\ve}_k(x),\omega_1+{\tau}_k(x)\right) - \xi\bigr|
-\log\bigl|f\left(T_{\omega_1}^kx ,\omega_1\right) - \xi\bigr|
\Big| \lesssim N^{\frac{\sigma}{10}}B_1(f)\eps
\]
for any $\xi,x$ as above, where
$\eps:=\max_{x,k}(|\ve_k(x)|+|\tau_k(x)|)$. Let
$\cE_{\omega_1}(N,\{0\}_k,\{0\}_k)$ be defined as in
Lemma~\ref{lem:3.4} and let $\cF_{\omega_1}(N,\xi)$ be the set
defined in \eqref{eq:3.LOG} both times applied to $f(x,\omega_1)$
with $\eps_k=\tau_k=0$. Then for any
$\xi\in\cJ_0(f)\setminus\cE_{\omega_1}(N,\{0\}_k,\{0\}_k)$ and any
$x\not\in \cF_{\omega_1}(N,\xi)$ one has
\[
\Big|\frac{1}{N}\sum_{\substack{1\le k\le N\\k\not\in J_N(x,\xi,\delta)}}
  \log|f(T_{\omega_1}^kx,\omega_1)-\xi| - \frac{1}{N}\sum_{1\le k\le N} \log|f(T_{\omega_1}^kx,\omega_1)-\xi|
\Big| \lesssim N^{-\frac{\sigma}{50}}
\]
The main part of the theorem now follows from Theorem~\ref{thm:A} and Remark~\ref{rem:Nga1ga2} applied to $f(\cdot,\omega_1)$.
For the final statement \eqref{eq:fattened} we use Remark~\ref{rem:3.fat} and estimate \eqref{eq:fett}.
\end{proof}

\begin{remark}\label{rem:3.desect} Inspection of the proof of Proposition~\ref{prop:A'}, in view of Remark \ref{rem:2.ldtsect},
shows that the following version of the statement holds: With
$\ve_k, \tau_k$ as before,
\begin{equation*}
\begin{split}
&\mes\Bigl\{x_1 \in \TT: \Big |{1\over N}\sum^N_{k=1}\, \log
\big|f\bigl(T^k_{\omega_1}(x_1,x_2) + \ve_k(x_1, x_2),\omega_1
 + \tau_k(x_1,x_2)\bigr) - \xi\big| - \langle \log \big| f(\cdot,
\omega_1) - \xi\big| \rangle\Big| > N^{-\kappa}\Bigr\} \\&\le
\exp\bigl(-N^{\kappa}\bigr)
\end{split}
\end{equation*}
for all $x_2 \in \TT$, $|\omega_0 - \omega_1| < \bigl[\bigl(1+
B_1(f)N^3\bigr]^{-1}$, $\xi \in
 \cJ_0(f) \setminus\left(\cE_{\omega_0}(x_2, N)\cup
 \cE_{\omega_0,\omega_1}\bigl(x_2, N, \{\ve_k\}, \{\tau_k\}\bigr)\right)$,
 where \[\mes \bigl[\cE_{\omega_0} (x_2, N)\bigr](\rho) \le
 N^{-\kappa},\quad
 \mes\left[\cE_{\omega_0, \omega_1}\bigl(x_2, N, \{\ve_k\}, \{\tau_k\}\bigr)\right](\rho)
 \le \exp\bigl(-{1\over 2} N^{\kappa}\bigr),\] $\rho = \exp\bigl(-N^\kappa\bigr)$.
\end{remark}

Now we are going to apply Theorem~\ref{thm:2.21} to the evaluation of the measure of those frequencies
$\omega$ for which so-called {\em resonances} occur. In the general setting of Section~\ref{sec:2} we define
a resonance by means of the following inequality, where $\kappa>0$ is small and fixed:
\begin{equation}\label{eq:reson}
\Big|\frac{1}{N} \sum_{1\le k\le N}\log|f(T_\omega^k (T_\omega^{\oN}x_0),\omega)-\xi| - \langle \log|f(\cdot,\omega)-\xi|\rangle \Big| > N^{-\kappa}
\end{equation}
where $\oN\gg N$. The goal is to show that the measure of those
$(\omega,\xi)$ for which \eqref{eq:reson} occurs for some
$e^{N^\sigma}>\oN\gg N$ is small for any fixed $x_0\in \TT^2$ (here
$\sigma>0 $ is another small constant).

\begin{theorem}\label{thm:3.6} Fix $x_0\in\TT^2$ and $N$ large.
\begin{enumerate}
\item Let $T_\omega: \TT^2 \to \TT^2$ be a shift and let $\omega_0$ be
$(N,\gamma_1,\gamma_2)$-Diophantine with some choice of small
$\gamma_1,\gamma_2>0$, $|\omega_1 - \omega_0| < \bigl[1 +
B_1(f)N\bigr]^{-1}$. Given $\oN > B_1(f)\,N^2$,
 there exist sets
\[ \cE_{\omega_0}(N), \,\widetilde\cE_{\omega_0,\omega_1}(N,\oN)
\subset\cJ_0(f),\] with
\[ \mes[\cE_{\omega_0}(N)](\rho)\lesssim N^{-\kappa},\quad \mes[\widetilde\cE_{\omega_0,\omega_1}(N,\oN )](\rho)\lesssim \exp(-N^{\kappa}),\qquad \rho=\exp(-N^\kappa) \]
so that for
$\xi\notin \cE_{\omega_0}(N)\cup \widetilde\cE_{\omega_0,\omega_1}(N,\oN)$
one has
\begin{multline*}
\mes\Bigl\{\theta\in[0,1]^2 : \Bigl|\frac{1}{N}\sum_{k=1}^N\log\bigl|f\bigl(x_0 + \oN \omega_1 + \theta +
k\big(\omega_1 + {\theta}/{\oN }\big),\omega_1 + {\theta}/{\oN }\bigr) - \xi\bigr| - \la\log\left|f(\cdot,\omega_1) - \xi\right|\ra\Bigr| > N^{-\kappa}\Bigr\} \\
< \exp\left(-N^{\kappa}\right)
\end{multline*}
where $\kappa>0$ is some small constant. The constants (but not the
sets $\cE_{\omega_0}(N),\,\widetilde\cE_{\omega_0,\omega_1}(N,\oN)$)
are uniform in the choice of $x_0$.

\item Let $T_\omega: \TT^2 \to \TT^2$ be the skew-shift and let
$\omega_0 \in \TT_{c,\ve_1,N}$ with some small $\ve_1 > 0$,
$|\omega_1 - \omega_0| < \bigl(1 + B_1(f)\bigr)^{-1}N^{-3}$.  Given
$\oN > B_1(f) N^4$ there exist $\cE_{\omega_0}(N)$,
$\cE_{\omega_0,\omega_1}(N,\oN) \subset \cJ_0(f)$ with
\[\mes\bigl[\cE_{\omega_0}(N)\bigr](\rho) < N^{-\kappa},\quad
\mes\bigl[\cE_{\omega_0,\omega_1}(N,\oN)\bigr](\rho)<
\exp(-N^\kappa),\] $\rho = \exp(-N^\kappa)$, so that for any $\xi
\notin \cE_{\omega_0}(N)\cup \cE_{\omega_0, \omega_1}(N,\oN)$ one
has
\begin{equation*}
\begin{split}
&\mes\biggl\{\theta \in \TT :  \Big|{1\over N} \sum^N_{k=1}\, \log
\big| f\bigl(T^k_\omega(T^{\oN}_{\omega_1 +
\theta/{\ooN}}(x_0,\omega), \omega_1 + \theta/{\ooN})\bigr) -
\xi\big|
 - \langle\log\big|f(\cdot,\omega_0) - \xi\bigr|\rangle \Big| >
N^{-\kappa}\biggr\} \\
& < \exp\bigl(-N^\kappa\bigr)
\end{split}
\end{equation*}
Here ${\ooN} = \oN(\oN -1)/2$, $\omega := \omega_1 + \theta/{\ooN}$,
$T^{\oN}_{(\omega_1 + \theta/{\ooN})}(x_0) = \bigl(x_1^{(0)} + \oN
x_2^{(0)} + {\ooN}\omega_1 + \theta, x_2^{(0)} + \oN \omega_1 +
\oN\theta/{\ooN}\bigr)$, $(x_1^{(0)}, x_2^{(0)}) = x_0$.
\end{enumerate}
\end{theorem}

\begin{proof} (1)~
Define $g : \TT^2\times\TT^2\rightarrow\IR$ by $g(\theta,\omega) = f\left(x_0 + \oN \omega_1 + \theta, \omega_1 + \omega\right)$.
Set ${\ve}_k(\theta) := k\theta / \oN$, ${\tau}_k(\theta) := \theta / \oN$, $k = 1,2,\ldots,N$ then
\[
f\left(x_0 + \oN \omega_1 + \theta + k\Big(\omega_1 + \frac{\theta}{\oN }\Big),\omega_1 + \frac{\theta}{\oN }\right) =
g\left(\theta + k\omega_1 + {\ve}_k(\theta), {\tau}_k(\theta)\right).
\]
Note that, for any $\theta\in[0,1]$, $\xi\in\IR$, $\la\log\left|f(\cdot,\omega_1) - \xi\right|\ra = \la\log\left|g(\cdot,0) - \xi\right|\ra$.
Set \[\widetilde\cE_{\omega_0,\omega_1}(N,\oN):=\cE_{\omega_0,\omega_1}(N,\{{\ve}_k\},\{{\tau}_k\}).\] Then
\begin{align*}
&\mes\Bigl\{\theta\in[0,1]^2  : \Bigl|\frac{1}{N}\sum_{k=1}^N\log\bigl|f\bigl(x_0 + \oN \omega_1 + \theta +
k\bigl(\omega_1 + {\theta}/{\oN }\bigr),\omega_1 + {\theta}/{\oN }\bigr) - \xi\bigr| - \la\log\left|f(\cdot,\omega_1) - \xi\right|\ra\Bigr| > N^{-\kappa}\Bigr\} \\
& = \mes\Bigl\{\theta\in[0,1]^2 : \Bigl|\frac{1}{N}\sum_{k=1}^N\log\Bigl|g\bigl(\theta + k\omega_1 +
{\ve}_k(\theta),{\tau}_k(\theta)\bigr) - \xi\Bigr| - \la\log\left|g(\cdot,0) - \xi\right|\ra\Bigr| > N^{-\kappa}\Bigr\} \\
& \leq \exp\left(-N^{\kappa}\right)
\end{align*}
for
$\xi\notin\cE_{\omega_0}(N)\cup\cE_{\omega_0,\omega_1}\left(N,\{{\ve}_k\},\{{\tau}_k\}\right)$, where these sets are as
 in Proposition~\ref{prop:A'}. 

(2)~Note that with $\omega = \omega_1 + \theta/\ooN$ one has
$$
T^k_\omega\left(T^{\oN}_\omega(x_1^{(0)}, x_2^{(0)})\right) =
T^k_{\omega_1}(y_1^{(0)}+\theta, y_2^{(0)} + \oN\theta/\ooN) + \ve_k
$$
where
\begin{align*}
y_1^{(0)} & = x_1^{(0)} + \oN x_2^{(0)} + \ooN \omega_1,\qquad y_2^{(0)} = x_2^{(0)} + \oN \omega_1\ ,\\[6pt]
\ve_k & = T^k_{(\theta / \ooN)}(0,0)\ .
\end{align*}
Invoking now Remark~\ref{rem:3.desect} (instead of Proposition~\ref{prop:A'}) one obtains the statement.
\end{proof}

\begin{lemma}\label{lem:3.9} Let $f(x,\omega)$ be $C^1$--smooth.
 Let $T_\omega$ be the shift (or the skew-shift).
 Given large $N$ there exists a set $\cJ(N) \subset
 \left\{(m_1, m_2): 1 \le m_1, m_2 \le \oN_1\right\}$, $\oN_1 = N^2$ (resp.~$\cJ(N) \subset [1, \oN_1]$, $\oN_1 = N^4$)
 and subsets $\cE_m(N) \subset \cJ_0(f)$, for every $m = (m_1, m_2) \in [1, N^2]^2 \setminus \cJ(N)$
 (resp.~$m \in [1, \oN_1]\setminus\cJ(N))$, such that
\begin{enumerate}
\item Using the notations \eqref{eq:calP} (resp.~\eqref{eq:B.P1}),
$$
\mes\bigcup_{m \in \cJ(N)} \cP_m(\oN_1) \le N^{-\kappa}
$$
for some small $\kappa$

\item For any $m \notin \cJ(N)$
$$
\mes\bigl[\cE_m(N)\bigr] (\rho) \le N^{-\kappa}\ ,\quad \rho =
\exp\bigl(-N^\kappa\bigr)
$$

\item Let $\oN \ge \oN_1$ and let $\omega_j^{(\oN)}$ be as in
Lemma~\ref{lem:A.neun} (resp.~as in \eqref{eq:B.P1}).  If
$\omega_j^{(\oN)} \in \cP_m^{(\oN_1)}$ for some $m \notin \cJ(N)$
then there exists $\cE_{m,j} (N, \oN)\subset \bigl[-B_0(f),
B_0(f)\bigr]$ with $\mes\cE_{m,j}(N,\oN) \le
\exp\bigl(-N^\kappa\bigr)$ such that for any $\xi \in \bigl[-B_0(f),
B_0(f)\bigr] \setminus \bigl(\cE_m(N)\cup \cE_{m,j}(N, \oN)\bigr)$
one has
\begin{equation}
\begin{split}\label{eq:mes_est}
& \mes \biggl\{\omega \in \cP_j^{(\oN)}\;:\: \Big|{1\over
N}\sum^N_{k=1}\, \log \big| f\bigl(T^k_\omega(T_\omega^{\oN} x_0),
\omega\bigr) - \xi| -   \langle \log \big|f(\cdot, \omega_j^{(\oN)})
- \xi\big |\rangle\Big| > N^{-\gamma}\biggr\} \\
&\le \bigl(\mes \cP_j^{(\oN)}\bigr) \exp\bigl(-N^\gamma\bigr)
\end{split}
\end{equation}
where $\gamma > 0$ is some small constant.
\end{enumerate}
\end{lemma}

\begin{proof} The proof is basically the same for the shift and skew-shift.  So, assume that $T_\omega$ is the shift
with $\omega \in \TT^2$. By Corollary~\ref{cor:A.10} there exists
$\cJ(N)\subset\{1,\ldots,N^2\}$ such that
\begin{itemize}
\item $\mes \bigcup_{m\in\cJ(N)} \cP_m^{(N^2)} \le N^{-\kappa}$
\item each $\omega_j^{(N^2)}$ with $m\not\in \cJ(N)$ is $(N,\gamma_1,\gamma_2)$-Diophantine for some small $\gamma_1,\gamma_2>0$.
\end{itemize}
Let $m\not\in\cJ(N)$ and let $\cE_m(N)=\cE_{\omega_m^{(N^2)}} (N)$ be as in Theorem~\ref{thm:3.6} with $ \omega_m^{(N^2)}$ in place
of $\omega_0$ and any $\omega_j^{(\oN)}\in \cP_m^{(N^2)}$ in place of~$\omega_1$. Since
$|\omega_j^{(\oN)}-\omega_m^{(N^2)}|\lesssim N^{-2}$ this is legitimate.
Then by Theorem~\ref{thm:3.6} there exists $\cE_{m,j}(N,\oN):=\cE_{\omega_0,\omega_1}(N,\oN)$ such that for any
\[ \xi\in [-B_0(f),B_0(f)]\setminus(\cE_m(N)\cup\cE_{m,j}(N,\oN)) \]
one has
\begin{equation}
\nonumber
\begin{aligned}
 & \mes\Big\{\omega=\omega_j^{(\oN)}+\theta/\oN,\,\theta\in\TT^2 \::\: \Big|\frac{1}{N}\sum_{k=1}^N \log|f(x_0+\oN\omega_j^{(\oN)}+\theta+k\omega,\omega)-\xi|
-\langle |f(\cdot,\omega_j^{(\oN)})-\xi|\rangle
\Big| > N^{-\gamma} \Big\} \\
&\le \oN^{-2}\exp(-N^\gamma)
\end{aligned}
\end{equation}
and the lemma follows. For the case of the skew-shift we use
Corollary~\ref{cor:B.metn}.
\end{proof}

The following proposition is the main result of this section. It captures the mechanism needed
for the elimination of bad~$\omega$ in the sense of~\eqref{eq:reson}.
The exceptional set of $\xi$ which appears in the proposition
will be converted into an exceptional set of energies in Section~\ref{sec:green}.

\begin{prop}\label{prop:3.10}
Let $f(x,\omega)$ be $C^1$--smooth.  Let $T_\omega$ be the shift (or
the  skew-shift).  Fix $x_0 \in \TT^2$.  Given large $N$ there
exists $Q(N) \subset \TT^2$ (resp. $Q(N) \subset \TT$) and for each
$\omega \notin Q(N)$ a subset $\cE_\omega(N) \subset \bigl[-B_0(f),
B_0(f)\bigr]$ such that
\begin{itemize}
\item $\mes Q(N) < N^{-\kappa}$, $\mes\bigl[\cE_\omega(N)\bigr](\rho) < N^{-\kappa}$, $\rho =
\exp\bigl(-N^{2\kappa}\bigr)$

\item For each $\omega \notin Q(N)$, $\xi \in \bigl[-B_0(f), B_0(f)\bigr]
\setminus \cE_\omega(N)$ and any $N^2 < \oN <
\exp\bigl(N^\beta\bigr)\ $ (resp.~$N^4<\oN<\exp(N^\beta)$) one has
$$
\Big| {1\over N} \sum^N_{k=1} \log \big| f\bigl(T^k(T^{\oN} x_0),
\omega \bigr) - \xi \big| - \langle \log |f(\cdot, \omega) - \xi|
\rangle\Big | < N^{-\gamma}
$$

\item Each $\omega \in \TT^2 \setminus Q(N)$ is $(N, \gamma_1, \gamma_2)$--Diophantine (resp.~$\TT^2 \setminus Q(N) \subset \TT_{c,\ve, N}$).
\end{itemize}

Here $\kappa, \beta, \gamma, \gamma_1, \gamma_2, \ve$ are small constants.
\end{prop}

\begin{proof} The proof is the same for shift and skew-shift.  So, let $T_\omega$ be a shift.
 Using the notations \eqref{eq:calP}, for any $\omega\in\TT^2$, and any positive integer $N_0$ one has
\begin{align}
&\#\Big\{ m=(m_1,m_2)\::\: 1\le m_1,m_2\le N_0,\;\exists\, \oN\ge N,\, j=(j_1,j_2),\, 1\le j_1,j_2\le \oN, \nn \\
&\qquad\qquad\qquad |\omega-\omega_j^{(\oN)}|\le 1/\oN,\, |\omega_j^{(\oN)}-\omega_m^{(N_0)}|\le 1/N_0\Big\}\le 25 \label{eq:overshoot}
\end{align}
Using the notations of Corollary~\ref{cor:A.10} set
\[
{Q}'(N) = \TT^2 \setminus \bigcup_{N^2<\oN<\exp(N^\beta)} \bigcup \big\{ \cP_j^{(\oN)}\::\: j=(j_1,j_2),\,1\le j_1,j_2\le \oN,\, j/\oN\in \cup_{m\not\in\cJ(N)} \cP_m^{(N^2)}
\big\}
\]
Clearly, $\mes {Q}'(N)\le 4\mes \big(\cup_{m\in \cJ(N)}  \cP_m^{(N^2)}\big) $. Assume that $\omega\in\TT^2\setminus {Q}'(N)$ and let for instance
$m,\oN,j$ be such that $\omega\in\cP_j^{(\oN)}$, $j/\oN\in\cP_m^{(N^2)}$, and $m\not\in\cJ(N)$. Set
\begin{align}
\cE_\omega'(N) &:= \bigcup\{ \cE_m(N)\::\: \text{$m,\oN,j$ as above},\;\oN\le\exp(N^\beta)\} \label{eq:4.unibig}\\
\cE_\omega''(N) &:= \bigcup\{ \cE_{m,j}(N,\oN)\::\: \text{$m,\oN,j$ as above},\;\oN\le\exp(N^\beta)\} \label{eq:4.unismall}
\end{align}
By \eqref{eq:overshoot}, the set \eqref{eq:4.unibig} consists of a union of  at most $25$ sets. Hence
\[
\mes \cE_\omega'(N) \less N^{-\kappa}
\]
On the other hand, since $\mes \cE_{m,j}(N,\oN)\le \exp(-N^\kappa)$,
\[
\mes \cE_\omega''(N) \less \exp(-N^{\kappa} +4N^\beta) \less \exp(-N^\kappa/2)
\]
provided $\beta<\kappa$ and $N$ is large. Due to the estimate \eqref{eq:mes_est} and Fubini's theorem there exists
$Q''(N)\subset\TT^2\setminus Q'(N)$ with $\mes Q''(N) \le \exp(-N^\gamma/4)$ such that for any $\omega\in\TT^2\setminus(Q'(N)\cup Q''(N))$
any $\xi\not\in \cE_\omega'(N)\cup\cE_\omega''(N)$ and any $N^2<\oN<\exp(N^\beta)$ we have
\[
\Big| N^{-1}\sum_{k=1}^N \log|f(x_0+\oN\omega+k\omega,\omega)-\xi|-\la\log|f(\cdot,\omega)-\xi|\ra\Big| \le N^{-\gamma}
\]
Set $Q(N)=Q'(N)\cup Q''(N)$, $\cE_\omega(N)=\cE'_\omega(N)\cup \cE_\omega''(N)$. Then $\mes Q(N)< N^{-\kappa}$ and $\mes E_\omega(N)< N^{-\kappa}$.
Moreover, inspection of the sets $\cE_m(N)$, $\cE_{m,j}(N,\oN)$ in \eqref{eq:4.unibig} and \eqref{eq:4.unismall} shows that, due to Theorem~\ref{thm:3.6},
$\mes [\cE_\omega(N)](\rho)\less N^{-\kappa}$ provided $\rho<\exp(-N^\kappa)$.
\end{proof}

\section{Large Deviation Theorems and Elimination of Resonances for Dirichlet Determinants}
\label{sec:4}

Let $T=T_\omega:\TT^2\to\TT^2$ be the shift or the skew-shift and let $V(x)\in C^1(\TT^2)$ be a real-valued function. Consider the Schr\"odinger equation
\begin{equation}
\label{eq:4.E}
-\psi(n + 1) - \psi(n - 1) + V(T^nx)\psi(n) = E\psi(n),\qquad n\in\IZ
\end{equation}

Let $H_{[a,b]}(x,\omega)$ be the operator defined by~\eqref{eq:4.E} on the interval $[a,b]$ with Dirichlet boundary condition $\psi(a-1)=\psi(b+1)=0$.
Let $f_{[a,b]}(x,\omega,E)$ be the characteristic determinant of $H_{[a,b]}(x,\omega)$, i.e.,
\[ f_{[a,b]}(x,\omega,E)=\det[H_{[a,b]}(x,\omega)-E] \]
We refer to $f_{[a,b]}(x,\omega,E)$ as the Dirichlet determinant.
Let $E_1^{[a,b]}(x,\omega) < E_2^{[a,b]}(x,\omega) < \cdots < E_N^{[a,b]}(x,\omega)$ be the eigenvalues of $H_{[a,b]}(x,\omega)$
with corresponding normalized eigenfunctions $\psi_1^{[a,b]}(x,\omega),\ldots,\psi_N^{[a,b]}(x,\omega)$.
We reserve the notations $H_N(x,\omega), f_N(x,\omega), E_j^{(N)}(x,\omega)$, and $\psi_j^{(N)}(x,\omega)$
for $[a,b]=[1,N]$.

The following lemma is a simple consequence of Weyl's comparison theorem for the eigenvalues of Hermitian
matrices, see for example Appendix~C in~\cite{C} or~\cite{Bhat}.

\begin{lemma}\label{lem:4.1}
Assume $1=a_1<b_1<b_1+1=a_2<b_2<\ldots<a_n<b_n\le N$. Then for any
$x\in\TT^2$, $E\in\IC$ one has
\[
\Big|\log\left|f_N(x,\omega,E)\right| - \sum_{k=1}^{n}\log|f_{[a_k,b_k]}(x,\omega,E)|\Big|\lesssim (n+N-b_n)\log [(B_0(V)+1)\eta^{-1}] .
\]
where
\[ \eta= \dist(E,\{E_j^{(N)}(x,\omega)\::\:1\le j\le N\}\cup\{E_j^{[a_k,b_k]}(x,\omega)\::\:1\le k\le n, 1\le j\le \ell_k\})\]
with $\ell_k=b_k-a_k+1$.
\end{lemma}

\begin{proof}
See Appendix~C in~\cite{C}.
\end{proof}

\begin{theorem}\label{thm:B}
Let $T_\omega$ be the shift or skew-shift on $\TT^2$ and suppose $V\in C^1(\TT^2)$ is a real-valued function.
Let $N$ be large and assume that $\omega$ is $(N,\gamma_1,\gamma_2)$--Diophantine (resp. $\omega \in \TT_{c,\ve, N}$). Then there exists
$\cE_\omega(N)\subset [-B_0(V)-2,B_0(V)+2]$ with $\mes(\cE_\omega(N))<N^{-\kappa}$ such that for any
\[ E\in [-B_0(V)-2,B_0(V)+2]\setminus \cE_\omega(N) \]
one has
\begin{equation}
\label{eq:det_ldt}
\mes\Big\{x\in\TT^2\::\: \Big|\log|f_N(x,\omega,E)|-\langle \log|f_N(\cdot,\omega,E)|\rangle \Big| > N^{1-\kappa}\Big\} \le \exp(-N^\kappa)
\end{equation}
where $\kappa>0$ is a small constant depending on the Diophantine condition.
\end{theorem}
\begin{proof}
The proof is the same for shift and skew-shift.  So, assume $T_\omega$ is a shift, $\omega \in \TT^2$.  Let
$\ell\asymp N^\beta$ be an integer with some small $\beta>0$. Let $n=[N\ell^{-1}]$, $k=N-n\ell$. Then by Lemma~\ref{lem:4.1}
with $a_m=(m-1)\ell+1$, $b_m=m\ell$, $1\le m\le n$, $a_{n+1}=b_n+1$, $b_{n+1}=N$ one has
\begin{equation}
\label{eq:4.weyl}
\Big|\log|f_N(x,\omega,E)|-\sum_{m=1}^n \log|f_\ell(x+(m-1)\ell\omega,\omega,E)|\Big| \lesssim N^{1-\beta/2}
\end{equation}
provided $N $ is large and
\begin{equation}\label{eq:4.spec_dist}
\begin{aligned}
&\min_{1\le j\le N} |E_j^{(N)}(x,\omega)-E| \ge \exp(-N^{\frac{\beta}{2}}) \\
&\min_{j,m} |E_j^{[a_m,b_m]}(x,\omega)-E| \ge \exp(-N^{\frac{\beta}{2}})
\end{aligned}
\end{equation}
There exists $\tcE_\omega(N)$ with $\mes\tcE_\omega{(N)}< \exp(-N^{\frac{\beta}{2}}/4) $
such that for any $E\not\in \tcE_\omega{(N)}$ one has
\begin{equation}\label{eq:4.mes}
\mes \{ x\in\TT^2\::\: \eqref{eq:4.spec_dist} \text{\ \ fails\ }\} \les \exp(-N^{\frac{\beta}{2}}/4)
\end{equation}
Note that
\[
\begin{aligned}
\log|f_N(x,\omega,E)| &= \sum_{j=1}^N \log|E_j^{(N)}(x,\omega)-E| \\
\log|f_\ell(x,\omega,E)| &= \sum_{j=1}^\ell \log|E_j^{(\ell)}(x,\omega)-E|
\end{aligned}
\]
Let $S$ be the set in \eqref{eq:4.mes}. Due to Lemma~\ref{lem:2.jensen}
one has
\begin{equation}\label{eq:4.aver}
\begin{aligned}
\Big|\la \log|E_j^{(N)}(\cdot,\omega)-E| \ra - \int_{\TT^2\setminus S} \log|E_j^{(N)}(x,\omega)-E|\, dx \Big| &\le \exp(-N^{\frac{\beta}{2}}/10) \qquad \forall\;1\le j\le N\\
\Big|\la \log|E_j^{(\ell)}(\cdot,\omega)-E| \ra - \int_{\TT^2\setminus S} \log|E_j^{(\ell)}(x,\omega)-E|\, dx \Big| &\le \exp(-N^{\frac{\beta}{2}}/10) \qquad \forall\;1\le j\le \ell
\end{aligned}
\end{equation}
provided $E$ does not fall into some set of measure $\lesssim \exp(-N^{\frac{\beta}{2}}/4)$. We may assume that $\tcE_\omega(N)$ contains that set.
In particular, due to \eqref{eq:4.weyl},
\begin{equation}
\label{eq:4.rate}
\Big| N^{-1} \la \log |f_N(\cdot,\omega,E)|\ra - \ell^{-1} \la \log |f_\ell(\cdot,\omega,E)|\ra  \Big| \lesssim N^{-\frac{\beta}{2}}\asymp\ell^{-\frac12}
\end{equation}
Since $\omega$ is $(N,\gamma_1,\gamma_2)$--Diophantine, it follows that $\{\ell\omega\}$ is $(N,\gamma_1,\gamma_2/2)$--Diophantine
 provided $\beta<\gamma_2/2$. Recall that the functions $E_j^{(\ell)}(x,\omega)$ are $C^1$--smooth with \[
B_0(E_j^{(\ell)})\le 2+B_0(V),\qquad
B_1(E_j^{(\ell)})\lesssim \ell (1+B_1(V)).\] Therefore, Theorem~\ref{thm:A} applies to each average
\[ \frac{1}{n} \sum_{m=1}^n \log|E_j^{(\ell)}(T_\omega^{m\ell}x,\omega)-E| \quad \forall\; 1\le j\le\ell \]
Let $\cE_\omega(\ell,j,n)$ stand for the set $\cT(N)$ from Theorem~\ref{thm:A} applied to the function $f(x)=E_j^{(\ell)}(x,\omega)$ and
the shift~$T_{\ell\omega}$.
Then
\[
\mes \bigcup_j \cE_\omega(\ell,j,n) \le \ell n^{-\kappa} \le n^{-\kappa/2} \text{\ \ provided\ \ }\beta\le \frac{\kappa}{2}
\]
and for all $1\le j\le \ell$
\[
\mes\Big\{x\in\TT^2\::\: \Big|\frac{1}{n}\sum_{m=1}^n\log|E_j^{(\ell)}(T_{\ell\omega}^mx,\omega)-E|-\langle \log|E_j^{(\ell)}(\cdot,\omega)-E|\rangle \Big| >
N^{1-\beta/2}\Big\} \le \exp(-N^{\frac{\beta}{2}})
\]
for any $E\not\in \bigcup_{j=1}^\ell \cE_\omega(\ell,j,n)$.
\end{proof}

Let $M_{[a,b]}(x,\omega,E)$ be the monodromy matrices of equation~\eqref{eq:4.E}. We reserve the notation $M_N(x,\omega,E)$ for $[a,b]=[1,N]$.
Recall that
\begin{equation}\label{eq:4.mon}
M_{[a,b]}(x,\omega,E)= \left[\begin{array}{cc} f_{[a,b]}(x,\omega,E) & -f_{[a+1,b]}(x,\omega,E) \\
                                              f_{[a,b-1]}(x,\omega,E) & - f_{[a+1,b+1]}(x,\omega,E)
                           \end{array}
                      \right]
\end{equation}
Note that Lemma~\ref{lem:4.1} implies the following assertion.

\begin{lemma}\label{lem:4.3}
One has
\begin{itemize}
\item
for any intervals $[s_i,t_i]$, $i=1,2$,
\begin{equation}
\label{eq:4.weyl1}
\log|f_{[s_1,t_1]}(x,\omega,E)| \le \log |f_{[s_2,t_2]}(x,\omega,E)| + (|s_2-s_1|+|t_2-t_1|)\log[(1+B_0(V))\eta^{-1}]
\end{equation}
 where
\[ \eta=\min\Big(\frac12,\dist(E,\spec H_{[s_2,t_2]}(x,\omega))\Big) \]
\item  for the monodromies, one has
\begin{equation}
\label{eq:4.normon}
0\le \log\|M_{[a,b]}(x,\omega,E)\| - \log |f_{[a,b]}(x,\omega,E)| \lesssim \log[(1+B_0(V))\eta^{-1}]
\end{equation}
 where
\[ \eta=\min\Big(\frac12,\dist(E,\spec H_{[a,b]}(x,\omega))\Big) \]
\end{itemize}
\end{lemma}
\begin{proof}
Estimate \eqref{eq:4.weyl1} follows from Lemma~\ref{lem:4.1}.
Applying \eqref{eq:4.weyl1} entry-wise, we conclude that \eqref{eq:4.normon} follows
from~\eqref{eq:4.mon}.
\end{proof}

For any $x,\omega,E$ clearly
\[ 0\le \log\|M_N(x,\omega,E)\| \le N\log(1+B_0(V)) \]
Next, we can draw the following conclusion from Lemma~\ref{lem:4.3}. Recall that $L(\omega,E)$ is
the Lyapunov exponent.

\begin{prop}
\label{prop:4.4}
There exists $\cF_{\omega}(N)$ with $\mes \cF_\omega(N)\le\exp(-N^\kappa/4)$ such that for any
\[ E\in [-B_0(V)-2,B_0(V)+2]\setminus\cF_\omega(N)\]
one has
\begin{align}
\Big| N^{-1} \la \log\|M_N(\cdot,\omega,E)\|\ra - N^{-1}\la \log|f_N(\cdot,\omega,E)|\ra \Big| &\le N^{-\kappa} \\
\Big| N^{-1} \la \log\|M_N(\cdot,\omega,E)\|\ra - \ell^{-1}\la \log\|M_\ell(\cdot,\omega,E)\|\ra \Big| &\le \ell^{-\frac12}
\end{align}
for any $\ell\asymp N^{-\beta}$. Here $\kappa,\beta>0$. In
particular, given $\ell$, there exists $\overline{\cF}_\omega(\ell)$
with $\mes \overline{\cF}_\omega(\ell) \le
\exp\bigl(-\ell^\kappa\bigr)$ such that for any $E \in
\bigl[-B_0(V) -2, B_0(V) + 2\bigr] \setminus
\overline{\cF}_\omega(\ell)$ one has \be \label{eq:4.conrate}
\big|\ell^{-1}
\la\log|f_\ell(\cdot,\omega,E)|\ra-L(\omega,E)\big|\less
\ell^{-\frac12} \ee
\end{prop}

\begin{proof}
We shall use the notations from the proof of Theorem~\ref{thm:B}. Thus,
\[
\Big| N^{-1} \la \log\|M_N(\cdot,\omega,E)\|\ra - \ell^{-1}\la \log\|M_\ell(\cdot,\omega,E)\|\ra \Big| \le \ell^{-\frac12}\qquad 
\]
for all $E\not\in \tcE_\omega^{(N)}$ (see \eqref{eq:4.rate} in the proof).
Assume $E\not\in\tcE_\omega^{(N)}$. Let $S$ be the set in~\eqref{eq:4.mes} so that \eqref{eq:4.spec_dist} is valid whenever $x\not\in S$
($S$ depends on~$E$). Due to Lemma~\ref{lem:4.3} one has
\[
\big|\log\|M_N(x,\omega,E)\|-\log|f_N(x,\omega,E)|\big|\lesssim N^{1-\beta/2}
\]
provided $x\not\in S$. Due to \eqref{eq:4.aver},
\[
\big|N^{-1}\la\log|f_N(\cdot,\omega,E)|\ra-N^{-1}\int_{\TT^2\setminus S} \log|f_N(\cdot,\omega,E)|\, dx\big|\lesssim \exp(-N^{\beta/2}/20)
\]
Since
\[ 0\le \log\|M_N(x,\omega,E)\| \le N\log(1+B_0(V)) \]
for any $x$, see \eqref{eq:4.mon},
\[
\big|N^{-1}\la\log\|M_N(\cdot,\omega,E)\|\ra-N^{-1}\int_{\TT^2\setminus S} \log\|M_N(\cdot,\omega,E)\|\, dx\big|\lesssim \log(1+B_0(V))\mes S\lesssim \exp(-N^{\beta/2}/20)
\]
and (4.11), (4.12) follow.  Given $\ell$, set $\ell_t =
\left[\ell^{(\frac{2}{\beta})^t}\right]$, $N_t =\ell_{t+1}$, $t = 0,
1,\dots$, $\overline{\cF}_\omega(\ell) = \mybigcup_{t=0}^\infty
\widetilde{\cE}_\omega(N_t)$. Then \eqref{eq:4.conrate} is
valid for $E \notin \overline{\cF}_\omega(\ell)$.
\end{proof}

\begin{corollary}
\label{cor:4.ldtmon} Let $T_\omega$ be the shift or the skew-shift
on $\TT^2$ and suppose $V\in C^1(\TT^2)$ is a real-valued function.
Let $N$ be large and assume that $\omega$ is
$(N,\gamma_1,\gamma_2)$--Diophantine (resp.~$\omega \in
\TT_{c,\ve,N}$). Then there exista $\cE_\omega(N)\subset
[-B_0(V)-2,B_0(V)+2]$ with $\mes \cE_\omega(N)< N^{-\kappa}$ such
that for any
\[ E\in [-B_0(V)-2,B_0(V)+2]\setminus \cE_\omega(N) \]
one has
\begin{align}
\mes\big\{ x\in\TT^2\::\: \big|\log\|M_N(x,\omega,E)\|-\la\log\|M_N(\cdot,\omega,E)\|\ra\big|> N^{1-\kappa}\big\} &\le \exp(-N^\kappa)\,,\\
\sup_{x\in\TT^2}\log\|M_N(x,\omega,E)\| &\le \la\log\|M_N(\cdot,\omega,E)\|\ra + N^{1-\kappa}\,, \\
\mes\big\{ x\in\TT^2\::\: \big|\log\|M_N(x,\omega,E)\|-\la\log|f_N(\cdot,\omega,E)|\ra\big|> N^{1-\kappa}\big\} &\le \exp(-N^\kappa)\,,\\
 \big|\la\log\|M_N(\cdot,\omega,E)\|\ra - \la\log|f_N(\cdot,\omega,E)|\ra\big| &\le N^{1-\kappa}
\end{align}
\end{corollary}

\begin{prop}\label{prop:4.ldsp}
Let us use the notations of Theorem~\ref{thm:B}. Then for any
$E\not\in \cE_\omega(N)$ the following holds: if for some
$x_1\in\TT^2$, $\dist(E,\spec
H_N(x_1,\omega))>(1+B_1(V))\exp(-N^\kappa/2)$, then
\[
\log|f_N(x_1,\omega,E)| > \la\log|f_N(\cdot,\omega,E)|\ra -
2N^{1-\kappa}
\]
\end{prop}
\begin{proof}
Due to Theorem~\ref{thm:B} there exists $x\in\TT$ such that $|x-x_1|<\exp(-N^\kappa)$ and
\begin{equation}
\label{eq:4.ldtx}
\log|f_N(x,\omega,E)| > \la \log|f_N(\cdot,\omega,E)|\ra - N^{1-\kappa}
\end{equation}
Since $E_j^{(N)}(x,\omega)$ are $C^1$--smooth with $\cB_1(E_j^{(N)})\lesssim B_1(V)$, one obtains
\[
| E_j^{(N)}(x,\omega) - E_j^{(N)}(x_1,\omega) | \lesssim B_1(V)\exp(-N^\kappa)
\]
Hence,
\begin{equation}\label{eq:4.comxx1}
\sup_{1\le j\le
N}\frac{|E_j^{(N)}(x,\omega)-E|}{|E_j^{(N)}(x_1,\omega)-E|} \le 1 +
C\,[\dist(E,\spec H_N(x,\omega))]^{-1}\exp(-N^\kappa)
\end{equation}
The proposition follows from \eqref{eq:4.ldtx} and \eqref{eq:4.comxx1}.
\end{proof}

The following proposition -- which is a consequence of our main
elimination method of Proposition~\ref{prop:3.10} -- shows that we
can insure that the large deviation theorem holds for a fixed phase
$x_0$ as long as we shift it by an amount $\oN\omega$ with $\oN\gg
N$; of course this requires the removal of a small set of
frequencies $\omega$ and energies~$E$ depending on~$x_0$.

\begin{prop}\label{prop:4.7}
Let $V\in C^1(\TT^2)$ and fix $x_0\in\TT^2$. Given large $N$, there exist a set $\cP(N)$ and for each $\omega\notin\cP(N)$
a subset $\cR_\omega(N)\subset[-B_0(V)-2,B_0(V)+2]$ so that
\begin{itemize}
\item with $\rho=\exp(-N^\kappa)$, we have
\[ \mes\cP(N)< N^{-\kappa}, \mes [\cR_\omega(N)](\rho)< N^{-\kappa}\]
\item for each $\omega\notin\cP(N)$, $E\in[-B_0(V)-2,B_0(V)+2]\setminus\cR_\omega(N)$, and any
$N^3<\oN<\exp(N^\beta)$ there is the bound
\[
\big| N^{-1}\log|f_N(x_0+\oN\omega,\omega,E)| - N^{-1}\la \log|f_N(\cdot,\omega,E)|\ra \big| \less N^{-\gamma}
\]
\end{itemize}
Here $\kappa,\beta,\gamma>0$ are small constants.
\end{prop}

\begin{proof}
We consider the case of the shift $T_\omega$, $\omega \in \TT^2$ and
shall use the notation from the proof of Theorem~\ref{thm:B}. By
Proposition~\ref{prop:3.10} applied to
$f(x,\omega)=E_j^{(\ell)}(x,\omega)$ there exists a set
$Q_j(N)\subset\TT^2$ and a subset $\cE_{\omega,j}(N)\subset
[-B_0(V)-2,B_0(V)+2]$ with $\mes Q_j(N)<N^{-\kappa}$ and
$\mes[\cE_{\omega,j}(N)](\rho)<N^{-\kappa}$, $\rho=\exp(-N^\kappa)$
so that for each $\omega\in\TT^2\setminus Q_j(N)$ and any
\[ E\in [-B_0(V)-2,B_0(V)+2]\setminus \cE_{\omega,j}(N) \]
as well as $N^2\le\oN\le\exp(N^\beta)$ one has
\[ \Big|\frac{1}{n}\sum_{k=1}^n \log|E_j^{(\ell)}(x_0+\oN\omega+k\omega,\omega)-E|
-\la\log|E_j^{(\ell)}(\cdot,\omega)-E|\ra \Big| \less N^{-\gamma}
\]
Hence,
\[
\Big|\frac{1}{n}\sum_{k=1}^n \log|f_\ell(x_0+\oN\omega+k\omega,\omega)-E|
-\la\log|f_\ell(\cdot,\omega)-E|\ra \Big| \less \ell N^{-\gamma} \less N^{-\gamma/2}
\]
for $\omega\in\TT^2\setminus \bigcup_{j=1}^\ell Q_j(N)$, $E\in [-B_0(V)-2,B_0(V)+2]\setminus \bigcup_{j=1}^\ell
\cE_{\omega,j}(N)$. As in the proof of Theorem~\ref{thm:B},
\[
\big|\log|f_N(x_0+\oN\omega,\omega,E)|-\sum_{k=1}^n \log|f_\ell(x_0+\oN\omega+k\ell\omega,\omega,E)|\big|\less N^{1-\gamma/2}
\]
provided
\begin{equation}\label{eq:4.19}
\begin{aligned}
&\min_{1\le j\le N}\quad\min_{N^2\le\oN\le\exp(N^\beta)} |E_j^{(N)}(x_0+\oN\omega,\omega)-E| \ge \exp(-N^{\frac{\gamma}{2}}) \\
&\min_{\substack{1\le j\le\ell\\0\le k\le n}}\quad  \min_{N^2\le\oN\le\exp(N^\beta)} |E_j^{(\ell)}(x_0+\oN\omega+k\ell\omega,\omega)-E| \ge \exp(-N^{\frac{\gamma}{2}})
\end{aligned}
\end{equation}
Given $\omega$, let $\cE'_\omega(N)$ be the set of
$E\in[-B_0(V)-2,B_0(V)+2]$ such that \eqref{eq:4.19} fails. Clearly,
$\mes [\cE'_\omega(N)](\rho)<\exp(-N^{\gamma/2}/2)$,
$\rho=\exp(-N^\gamma)$ provided $\beta\ll\gamma$.  Finally,
define
\[
 \cP(N)=\bigcup_{j=1}^\ell Q_j(N), \quad \cR_\omega(N) = \cE'_\omega(N)\cup \bigcup_{j=1}^\ell \cE_{\omega,j}(N)
\]
and the proposition is proved.
\end{proof}

\section{Estimates on the Green function and the proof of Theorem \ref{th:1.1}}
\label{sec:green}

By Cramer's rule,
\begin{equation}\label{eq:cramer}
\begin{aligned}
G_{[a,b]}(x,\omega,E)(m,n) &:= (H_{[a,b]}(x,\omega)-E)^{-1}(m,n) \\
&= \frac{f_{[a,m-1]}(x,\omega,E)f_{[n+1,b]}(x,\omega,E)}{f_{[a,b]}(x,\omega,E)}
\end{aligned}
\end{equation}
for all $a\le m\le n\le b$. To evaluate the Green function $G_{[a,b]}(x,\omega,E)$ we need to obtain appropriate
estimates on the Dirichlet determinants in \eqref{eq:cramer} which are uniform in $x,m,n$. To derive such estimates
we need to modify the proof of Theorem~\ref{thm:B} slightly; in fact, we refer the reader to Remark~\ref{rem:2.18} for these matters.
We shall use the notations from the proof of Theorem~\ref{thm:B}. Note that for any $x,\omega,E$ and $1\le N'\le N$
\[
\log\|M_{N'}(x,\omega,E)\| \le \sum_{m=1}^{n'} \log\|M_\ell(T_{\ell\omega}^{m-1}x,\omega,E)\|+\ell\log[2B_0(V)+4]
\]
where $n'=[\ell^{-1}N']$. This is because
\[
M_{N'}(x,\omega,E) = M_{[n'\ell,N']}(x,\omega,E)\prod_{m=n'}^1 M_{[a_m,b_m]}(x,\omega,E)
\]
where $a_m,b_m$ are as in the proof of Theorem~\ref{thm:B}.
Assume that $\omega$ is $(N,\gamma_1,\gamma_2)$--Diophantine.
Recall that due to Remark~\ref{rem:2.18}
one then has
\begin{equation}
\label{eq:5.supj}\sup_{n^{\frac12}\le n'\le n}\sup_{\#\cB\le (n')^{1-2\kappa}}\sup_{x\in\TT^2}\frac{1}{n'}
\sum_{m\in[1,n']\setminus\cB}\log\left|E_j^{(\ell)}(x+(m-1)\ell\omega,\omega) - E\right| \le \la\log|E_j^{(\ell)}(\cdot,\omega) - E|\ra+ n^{-\kappa}
\end{equation}
for any
$E\in[-B_0(V)-2,B_0(V)+2]\setminus\widetilde\cE_\omega(\ell,j,n)$
where $\mes \widetilde\cE_\omega(\ell,j,n)\le 2n^{-\frac{\kappa}{2}}$.
Since $\ell\le N^\beta$ and taking $\beta\ll\kappa$ one has the following assertion

\begin{lemma}
\label{lem:5.1} Let $f\in C^1(\TT^2)$ and assume that  $\omega$ is
$(n,\gamma_1,\gamma_2)$--Diophantine with some large $n$
(resp.~$\omega \in \TT_{c,\ve, N}$ for the case of the skew-shift).
There exists a set $\cE'_\omega(n)$, $\mes
\cE'_\omega(n)<n^{-\frac{\kappa}{4}}$ so that
\begin{equation}
\label{eq:5.unichez}\sup_{n^{\frac12}\le n'\le n}\;\;\sup_{\#\cB\le (n')^{1-2\kappa}}\sup_{x\in\TT^2}\frac{1}{n'}\;\;
\sum_{m\in[1,n']\setminus\cB}\log\left|f_{\ell}(T_\omega^{(m-1)\ell},\omega,E)\right| \le
\la\log|f_{\ell}(\cdot,\omega,E)|\ra+ n^{-\kappa}
\end{equation}
Similar estimates hold for $f_{[a,\ell+b]}$ with $|a|,|b|\le 1$. Moreover, the average on the right-hand side
of~\eqref{eq:5.unichez} can be kept the same for all $f_{[a,\ell+b]}$, $|a|,|b|\le1$.
\end{lemma}

\begin{proof}
Adding up \eqref{eq:5.supj} over $1\le j\le\ell$ one obtains \eqref{eq:5.unichez}. The same arguments are valid
for $f_{[a,\ell+b]}$. The last part follows from Lemma~\ref{lem:4.3}.
\end{proof}

To proceed we need to compare the following two sums:
\[
\sum_{m=1}^{n'} \log\|M_\ell(x+(m-1)\ell\omega,\omega,E)\|
\]
versus
\[
\sum_{m=1}^{n'} \log|f_\ell(x+(m-1)\ell\omega,\omega,E)|
\]

\begin{lemma}\label{lem:5.2} Under the assumptions of the previous lemma
there exists $\widetilde\cE_\omega(\ell,n)\subset [-B_0(V)-2,B_0(V)+2]$ with $\mes\widetilde\cE_\omega(\ell,n)
\le n^{-\kappa_1}$ such that for any $E\in [-B_0(V)-2,B_0(V)+2]\setminus\widetilde\cE_\omega(\ell,n)$, $n^{\frac12}\le n'\le n$,
and any $x\in\TT^2$ one has
\[
\#\big\{ 1\le m\le n'\::\:
\log\|M_\ell(T_\omega^{m\ell},\omega,E)\|>\log|f_\ell(T_\omega^{m\ell},\omega,E)|+
\log[(B_0(V)+1)n] \big\} \le (n')^{1-\kappa_2}
\]
Here $0<\mu\ll\kappa_1$, $\kappa_2,\kappa\ll1$, $\ell\asymp n^{2\mu}$.
\end{lemma}

\begin{proof} We consider the case of a shift $T_\omega$, $\omega \in \TT^2$.
Due to Lemma~\ref{lem:4.3}
\[
\log\|M_\ell(x+m\ell\omega,\omega,E)\|\le \log|f_\ell(x+m\ell\omega,\omega,E)|+ \log[(B_0(V)+1)n]
\]
unless
\begin{equation}
\label{eq:5.minlj}
\min_j |E-E_j^{(\ell)}(x+m\ell\omega,\omega)|< n^{-\kappa}
\end{equation}
Due to Corollary~\ref{cor:2.xiset}, one can find $\widetilde\cE_\omega(\ell,j,n)$ with $\mes \widetilde\cE_\omega(\ell,j,n)
\le n^{-\kappa/2}$ such that for any $E\in [-B_0(V)-2,B_0(V)+2]\setminus\widetilde\cE_\omega(\ell,j,n)$ and any $x\in\TT^2$
one has
\[
\#\{ 1\le m\le n'\::\: |E-E_j^{(\ell)}(x+m\ell\omega,\omega)|< n^{-\kappa}\} \le (n')^{1-\kappa/2}
\]
Let $\widetilde\cE_\omega(\ell,n)=\bigcup_j \widetilde\cE_\omega(\ell,j,n)$. Then $\mes \widetilde\cE_\omega(\ell,n)
\le n^{-\kappa/3}$ and for any $E\in [-B_0(V)-2,B_0(V)+2]\setminus\widetilde\cE_\omega(\ell,n)$, $x\in\TT^2$ we have
\[
\#\{ 1\le m\le n'\::\: \eqref{eq:5.minlj} \text{\ \ fails }\} \le \ell (n')^{1-\kappa/2} \le (n')^{1-\kappa/3}
\]
as desired.
\end{proof}

We can now prove the following uniform upper bound for $\log\|M_N(x,\omega,E)\|$.

\begin{prop}\label{prop:5.3}
Let $V\in C^1(\TT^2)$, and assume $N$ is large. Let $T_\omega$ be a
shift or a skew-shift. Suppose that $\omega$ is
$(N,\gamma_1,\gamma_2)$-Diophantine (resp. $\omega \in
\TT_{c,\ve,N}$). There exists
$\cE_\omega'(N)\subset[-B_0(V)-2,B_0(V)+2]$ with $\mes
\cE'_\omega(N)< N^{-\kappa}$ such that
\[\sup_{N^{\frac12}\le N'\le N}\sup_{x\in\TT^2} \frac{1}{N'}
\log\|M_{N'}(x,\omega,E)\| \le
\frac{1}{N}\la\log\|M_{N}(\cdot,\omega,E)\|\ra+ N^{-\kappa}
\]
for all $E\in [-B_0(V)-2,B_0(V)+2]\setminus \cE_\omega'(N)$.
\end{prop}

\begin{proof}
We shall use the notations from the proof of Theorem~\ref{thm:B}. Let
$E\in [-B_0(V)-2,B_0(V)+2]\setminus \widetilde\cE_\omega(\ell,n)$
where $\ell\asymp N^\mu$, $n=[N\ell^{-1}]$, and $\widetilde\cE_\omega(\ell,n)$
is defined in Lemma~\ref{lem:5.2}. Then for any $n^{\frac12}\le n'\le n$, $x\in\TT^2$ one has
\begin{equation}
\label{eq:5.REL}
\begin{aligned}
\sum_{m=1}^{n'} \log\|M_\ell(\mapt,\omega,E)\|&\le \sum_{m\in[1,n']\setminus\cB_{n'}(x,\omega,E)}
\log|f_\ell(\mapt,\omega,E)|\\
&+ n'\log[(B_0(V)+1)n] + (n')^{1-\kappa_1} \sup_{y} \log\|M_\ell(y,\omega,E)\|
\end{aligned}
\end{equation}
where
\[\cB_{n'}(x,\omega,E) = \{1\le m\le n'\::\: \log\|M_\ell(\mapt,\omega,E)\|
> \log|f_\ell(\mapt,\omega,E)|+ \log[(B_0(V)+1)n]\}
\]
and $\#\cB_{n'}(x,\omega,E)\le (n')^{1-\kappa}$. Combining \eqref{eq:5.REL} with Lemma~\ref{lem:5.1} yields
\[
\frac{1}{n'\ell} \sum_{m=1}^{n'} \log\|M_\ell(\mapt,\omega,E)\|\le \frac{1}{\ell} \la\log|f_\ell(\cdot,\omega,E)|\ra
+ \frac{1}{\ell}\log[(B_0(V)+1)n] + \ell^{-1}n^{-\kappa/2}
\]
Finally, by Proposition~\ref{prop:4.4},
\[
\ell^{-1}\la\log|f_\ell(\cdot,\omega,E)|\ra \le N^{-1}\la\log|f_N(\cdot,\omega,E)|\ra + \ell^{-1/2}
\]
The proposition is proved.
\end{proof}

Combining Propositions~\ref{prop:4.7} and~\ref{prop:5.3} yields the following.

\begin{prop}\label{prop:5.4}
Let $V\in C^1(\TT^2)$ and fix $x_0\in\TT^2$. Let $T_\omega$ be the
shift or the skew-shift. Assume that $L(\omega,E)\ge L_0>0$ for any
$\omega$ and $E\in (E_1,E_2)$. Given large $N$ there exists a set
$\cP(N)$ and for each $\omega\notin\cP(N)$ a subset
$\cK_\omega(N)\subset \bigl[-B_0(V) - 2, B_0(V) +2\bigr]$ so that
\begin{itemize}
\item[(a)] with $\rho=\exp(-N^\kappa)$ one has
\[ \mes\cP(N) < N^{-\kappa},\quad \mes[\cK_\omega(N)](\rho)<N^{-\kappa} \]
\item[(b)] for each $\omega\notin\cP(N)$, $E\in  [-B_0(V)-2,B_0(V)+2]\setminus \cK_\omega(N)\cap (E_1,E_2)$
and $N^3\le\oN\le\exp(N^\kappa)$ one has
\begin{equation}\label{eq:5.GRDEC}
|G_N(x_0+\oN\omega,\omega,E)(m,n)|\le \exp(-L_0|m-n|/2)
\end{equation}
for any $|m-n|>N/2$
\item[(c)]  with $N_0=N^3$, $N_0\ll N_1\le \exp(N^\beta)$, for any $\omega\notin\cP(N)$,
and any \[ E\in
 [-B_0(V)-2,B_0(V)+2]\setminus [\cK_\omega(N)](\rho)\cap (E_1,E_2)\] we have
 $\dist(\spec H_{[N_0,N_1]}(x_0,\omega),E)>\exp(-N^\kappa)$
and
\[
|G_{[N_0,N_1]}(x_0,\omega,E)|(m,n) \le \exp(-L_0|m-n|/3)
\]
for any $|m-n|>N/2$.
Here, $\beta,\kappa$ are as in Proposition~\ref{prop:4.7}
\end{itemize}
\end{prop}

\begin{proof}
Let $\cP(N)$, $\cE_\omega(N)\subset [-B_0(V)-2,B_0(V)+2]$, $\omega\in\TT^2\setminus\cP(N)$
be as in Proposition~\ref{prop:4.7}. Then every
\[
E\in [-B_0(V)-2,B_0(V)+2]\cap (E_1,E_2)\setminus \cE_\omega(N)
\]
satisfies, for any $\oN$ as above,
\begin{equation}\label{eq:5.LOWB}
\begin{aligned}
|f_N(T_\omega^{\oN}(x_0),\omega,E)| &> \exp\big(\la\log|f_N(\cdot,\omega,E)|\ra-N^{1-\gamma}\big) \\
& > \exp(NL(\omega,E)-N^{1-\gamma}) > \exp(NL(\omega,E)/2)
\end{aligned}
\end{equation}
provided $N$ is large. Here we used Proposition~\ref{prop:4.4}. Due to Proposition~\ref{prop:3.10} each
$\omega\notin\cP(N)$ is $(N,\gamma_1,\gamma_2)$--Diophantine (resp. $\omega \in T_{c,\ve, N}$).
Therefore, by Proposition~\ref{prop:5.3}
there exists $\widetilde\cE_\omega(N)$, $\mes\widetilde\cE_\omega(N)<N^{-\kappa}$ such that for any
\[ E\in  [-B_0(V)-2,B_0(V)+2]\cap (E_1,E_2)\setminus \widetilde\cE_\omega(N) \]
any $x\in\TT^2$ and any interval $[s,t]$ with $N^{\frac12}<t-s\le N$
one has
\begin{equation}
\label{eq:5.UPPB}
|f_{[s,t]}(x,\omega,E)|\le \exp
\big((t-s)(N^{-1}\la\log|f_N(\cdot,\omega,E)|\ra) + (t-s)^{1-\kappa}\big) \le \exp\big((t-s)L(\omega,E)+2(t-s)^{1-\kappa}\big)
\end{equation}
Assume $\omega\notin\cP(N)$, and
\[
 E\in  [-B_0(V)-2,B_0(V)+2]\cap (E_1,E_2)\setminus (\widetilde\cE_\omega(N)\cup\cE_\omega(N))
\]
Then \eqref{eq:cramer}, \eqref{eq:5.LOWB}, \eqref{eq:5.UPPB} imply \eqref{eq:5.GRDEC}
when $a=1<N^{\frac12}\le m<n<N-N^{\frac12}$ and $n-m>N^{1-\kappa}$. To prove~\eqref{eq:5.GRDEC} when
$1\le m\le N^{\frac12}$ or $N-N^{\frac12}<n<N$ one can use the trivial upper bound
\[
|f_{[s,t]}(\omega,\omega,E)|\les \exp(2(s-t)(2+B_0(V)))
\]
Thus (b) holds. We invoke now the following general fact which is
valid for general discrete Schr\"odinger equations: if for some $E$
the estimate
\begin{equation}
\label{eq:5.gencr}
|G_{[a',a'+N]}(x,\w,E)(m,n)|\le \exp(-L|m-n|)
\end{equation}
holds for all $a'\in[a,b]$, $|m-n|>N/2$, then $E\notin\spec
H_{[a,b]}(x,\w)$ provided $N>(\log(b-a))^2$, $b-a>R_0(L)$ where $R_0(L)$
is a suitable constant (see for example Appendix~C in \cite{C}).
>From this and (b) we conclude that $E\not\in \spec
H_{[N_0,N_1]}(x_0,\omega)$ for any $\omega\notin\cP(N)$,
\[
 E\in  [-B_0(V)-2,B_0(V)+2]\cap (E_1,E_2)\setminus \cK_\omega(N)
\]
That means, in particular, if
\[
 E\in  [-B_0(V)-2,B_0(V)+2]\cap (E_1,E_2)\setminus [\cK_\omega(N)](\rho)
\]
then $(E-\rho,E+\rho)\cap \spec
H_{[N_0,N_1]}(x_0,\omega)=\emptyset$.

Finally, to complete the proof of (c) we refer to yet another
general fact about Schr\"odinger equations: if for some $E$ the
estimate~\eqref{eq:5.gencr} holds for all $a'\in[a,b]$ and
$\dist(E,\spec H_{[a,b]}(x,\w))>\exp(-N^\kappa)$, then
\[
|G_{[a,b]}(x,\w,E)(m,n)|\le \exp(-L_1|m-n|/2)
\]
for all $|m-n|>N/2$, provided $N>(\log(b-a))^2$, $b-a>R_0(L_1)$.
\end{proof}

\begin{remark}\label{rem:5.5} A similar statement is valid for the Green
functions $G_N(x_0 - \oN\omega, \omega, E)(m, n)$, $N^3 \le \oN \le
\exp\bigl(N^\alpha\bigr)$ and $G_{[-N_1, -N_0]}(x_0, \omega, E)(m,
n)$.  We will use the same notations $\cP(N), \cK_\omega(N)$ for the
exceptional sets of $\omega \in \TT^2$ and $E \in [-B_0(V)-2, B_0(V)
+2]$ needed to guarantee the estimates for these Green
functions as for the Green function in Proposition \ref{prop:5.4}.
\end{remark}

\begin{corollary}\label{cor:5.6} Assume that some function $\psi(n)$, $n \in \IZ$ obeys
\be \label{eq:5.9} -\psi(n+1) - \psi(n-1) + V\bigl(T^n_\omega
x_0\bigr) \psi(n) = E\psi(n)\ \text{for}\ -N_1 \le n \le N_1 \ee for
some $\omega \notin \cP(N)$, $E \in \bigl[-B_0(V) - 2, B_0(V) +
2\bigr] \setminus\cK_\omega(N)$ where $\cP(N), \cK_\omega(N)$ are as
in Proposition~\ref{prop:5.4} and $\exp\bigl(N^\kappa\bigr)\ge N_1$.
Assume in addition that
$$
\max_{|n|\le N_1}\ |\psi(n)| \le 1
$$
then for any $N^3 < |n| < N_1$ holds
$$
|\psi(n)| \lesssim \exp\left(-L_0 \min \bigl(|n| - N^3, N_1 - |n|\bigr)/2\right)\ .
$$
If $\psi(n)$ in addition satisfies the Dirichlet boundary conditions
on $[-N_1 +1, N_1 - 1]$, i.e., $\psi(-N_1) = \psi(N_1) = 0$, then
$$
|\psi(n)| \le \exp\left(-L_0\bigl(|n| - N^3\bigr)/2\right)
$$
for any $N^3 < |n| \le N_1$.
\end{corollary}

\begin{proof} Let $N^3 < n < N_1$.  Then
$$
\psi(n) = \sum_{m \in \{N^3, N_1\}} G_{[N^3, N_1]}(x_0, \omega,
E)(n, m) \psi(m)\ .
$$
Both estimates follow now from Proposition~\ref{prop:5.4} (see also Remark~\ref{rem:5.5}).
\end{proof}

\begin{proof}[Proof of Theorem \ref{th:1.1}] We shall use the
notations from the Proposition~\ref{prop:5.4}.  Fix $x_0$. Given
$N$, define $N_t := \left[\exp\bigl(N_{t-1}^{\kappa\over
9}\bigr)\right]$, $t = 1, 2,\dots$, $N_0 = N_1$, $\overline{\cP}(N)
= \mybigcup_{t \ge 0} \cP(N_t)$, $\overline{\cK}_\omega(N) =
\mybigcup_{t\ge 0}\bigl[\cK_\omega(N_t)\bigr](\rho_{N_t})$, for
$\omega \notin \overline{\cP}(N)$ where $\cP(N_t),\cK_\omega(N_t)$
are as in Remark~\ref{rem:5.5}, $\rho_{N_t} =
\exp\bigl(-(N_t)^\kappa\bigr)$. Then
\begin{align}\label{eq:5.mes}
\mes \overline{\cP}(N) & \le \sum\, N_t^{-\kappa} \lesssim N^{-\kappa}\\
\mes \overline{\cK}_\omega(N) & \le \sum\, N_t^{-\kappa} \lesssim N^{-\kappa}\nn
\end{align}
If $\psi(n)$ obeys
\be \label{eq:5.poly}
\begin{aligned}
& -\psi(n+1) - \psi(n-1) + V\bigl(T^n_\omega x_0\bigr) \psi(n) = E\psi(n)\qquad n \in \IZ^1\\
& |\psi(n)| \le |n|^2
\end{aligned}
\ee with $\omega \notin \overline{\cP}(N)$ and $E \in
\bigl[-B_0(V)-2, B_0(V)+2\bigr] \setminus \overline\cK_\omega(N)$,
then due to Corollary~\ref{cor:5.6}, \be \label{eq:5.dec}
|\psi(n)| \le \min_{s = t-1, t, t+1} N_s^2\exp\left(-L_0
\min\left(\bigl(|n|- N^3_s, N_{s+1} - |n|\bigr)\right)\right) \le
\exp(-L_0|n|/4) \ee where $N_t \le n < N_{t+1}$.
Theorem~\ref{th:1.1} follows from \eqref{eq:5.dec} and
\eqref{eq:5.mes}.
\end{proof}

\section{Skew shifted $C^1$--potentials at large disorder}\label{sec:6}

Consider \be \label{eq:6.1} H(x,y,\lambda, \omega)\psi(n):= -
\psi(n-1)- \psi(n+1) + \lambda V\bigl(T^n(x,y)\bigr) \psi(n),\ n \in
\IZ^1\ , \ee where $V(x,y), (x,y) \in \TT^2$ is a real valued
$C1$--function, $T = T_\omega: \TT^2 \to \TT$ is the skew--shift
$T_\omega(x,y) = (x+y, y + \omega)$, $\lambda$ is a parameter.  Let
$f_N(x,y, \lambda, \omega, E)$ be the characteristic determinant of
the operator $H_N(x,y,\lambda,\omega)$ which is the restriction
of $H(x,y,\lambda,\omega)$ on $[1, N]$ with Dirichlet boundary
conditions, i.e., \be\label{eq:6.2} f_N(x,y,\lambda,\omega,E) =
\det\begin{vmatrix}
\lambda V_1 - E & -1 & 0 & 0 &\cdots & 0\\
-1 & \lambda V_2 - E & -1 & 0 & \cdots &\\
0 & -1 & \lambda V_3 - E & -1 & \cdots & 0\\
\hdotsfor{5} & -1\\
0 & 0 & \cdots & 0 & -1 & \lambda V_N-E\end{vmatrix}
\ee
where $V_j = V\bigl(T^j(x,y)\bigr)$.  Recall that the monodromy matrices are as follows
\be \label{eq:6.mon}
M_N(x,y,\lambda,\omega, E) =
\begin{bmatrix}
f_N(x,y,\lambda,\omega, E) & - f_{N-1}\bigl(T(x,y),\lambda,\omega, E\bigr)\\
f_{N-1}(x,y,\lambda,\omega, E) & - f_{N-2}\bigl(T(x,y),\lambda,\omega, E\bigr)
\end{bmatrix}
\ee Consider also the following diagonal matrix \be\label{eq:6.diag}
D_N(x,y,\lambda,\omega,E) = \diag(\lambda V_1 - E,\dots, \lambda V_N
- E)(x,y)\ee

\begin{lemma}\label{lem:6.1}
There exists $\lambda_0 = \lambda_0\bigl(B_0(V)\bigr)$ such that for
$|\lambda| \ge \lambda_0$ and with $N \asymp \lambda^{1/2}$ the
following assertion holds: there exists $\cE_{\omega,\lambda}(N)
\subset \bigl[-\lambda B_0(V) - 2, \lambda B_0(V) + 2\bigr]$ with
$\mes\cE_{\omega,\lambda}(N) < N^{-\kappa}$ such that for any $E \in
\bigl[-B_0(V)-2, B_0(V)+2\bigr]\setminus
\cE_{\omega,\lambda}(N)$ one has
$$
 {1\over N} \langle \log
\big|f_N(\cdot,\lambda,\omega,E)\big|\rangle > {1\over 2} \log
|\lambda|
$$
and
$$
\mes\left\{(x,y) \in \TT^2 : \Bigl|\frac{1}{N}\log\big
|f_N(x,y,\lambda,\omega, E)\big| - \frac{1}{N} \langle \log
\big|f_N(\cdot,\lambda,\omega,E)\big|\rangle\Big| >
N^{-\kappa}\right\} \le N^{-\kappa}
$$
\end{lemma}

\begin{proof} The matrix in the right-hand side of \eqref{eq:6.2} can be written in the form $D_N(x,y,\lambda,\omega,E) + B_N$, where $D_N$ is given by \eqref{eq:6.diag}.  Clearly $\|B_N\| = 2$ and
\be \label{eq:6.detd}
 {1\over N}\log|\det D_N(x,y,\lambda,\omega,E)| = \log |\lambda| + N^{-1} \sum\limits^N_{m=1} \log |V
\bigl(T^m(x,y)\bigr) - E/\lambda|\ . \ee By Lemma~\ref{lem:2.jensen}
and Theorem~\ref{thm:A}, there exists $L_{\lambda,\omega}(N) \subset
[-B_0(V)-2/\lambda, B_0(V)+2/\lambda]$ with $\mes L_{\omega,\lambda}(N) < N^{-\kappa}$ such
that for any $E/\lambda \in [-B_0(V)-2/\lambda, B_0(V)+2/\lambda]\setminus
L_{\omega,\lambda}(N)$ holds
 \begin{align}
 & \mes\bigl\{(x,y)\in\TT^2: \min_{1\le m\le N} | V\bigl(T^m(x,y)\bigr) - E/\lambda| \le |\lambda|^{-1/2}\bigr\}\le |\lambda|^{-1/4}\label{eq:6.minv}\\[6pt]
 & \mes \Bigl\{(x,y) \in \TT^2: \Big|N^{-1} \sum^N_{m=1} \log \big| V\bigl( T^m(x,y)\bigr) - E/\lambda\big| -\label{eq:6.ldv}\\
  &\qquad \langle \log \big| V(\cdot) - E/\lambda\big| \rangle \Big| > N^{-\kappa}\Bigr\} \le \exp(-N^{\kappa}) = \exp(-\lambda^{-\kappa_1})\nn \\[6pt]
 & \Big| \langle \log \big| V(\cdot) - E/\lambda\big|\rangle - \int_{\TT^2\setminus B_{\omega,\lambda}(N)} \log \big|V(x) - E/\lambda\big| dx\Big| \le N^{-\kappa}\label{eq:6.appr}
 \end{align}
 Then
 \be \label{eq:6.dinv}
 \big\|D_N(x,y,\lambda,\omega,E)^{-1}\big\| \le \lambda^{-1/2}\ \text{for any}\ (x,y) \in \TT^2 \setminus B_{\omega,\lambda}(N)
 \ee
 where $B_{\omega,\lambda}(N)$ is the set in \eqref{eq:6.minv}, $\mes B_{\omega,\lambda}(N) < \lambda^{-1/4}$
 \be \label{eq:6.logd}
 {1\over N} \langle \log |\det D_N(\cdot, \lambda,\omega, E)|\rangle > {1\over 2}\log|\lambda|\ .
 \ee
Note that \eqref{eq:6.dinv} implies
\begin{gather*}
\| D_N^{-1} B_N\| < 2|\lambda|^{-1/2}\\[6pt]
(1+2|\lambda|^{-1/2})^N > |\det(1 + D_N^{-1} B)| > (1 - 2 |\lambda|^{-1/2})^N > \exp(-4N|\lambda|^{-1/2}) \gtrsim 1\\[6pt]
\log|f_N| = \log|\det(D_N + B_N)| = \log|\det D_N| + O(|\lambda|^{-1/2})
\end{gather*}
Let $E_1^{(N)}(x,y) < \dots < E_N^{(N)}(x,y,\omega)$ be the eigenvalues $H_N(x,y,\lambda,\omega)$.  Let $S$ be the union of the sets in \eqref{eq:6.minv} and \eqref{eq:6.ldv}.  Then due to Lemma~\ref{lem:2.jensen}, there exists $\widetilde L_{\omega,\lambda}(j, N)$ with $\mes \widetilde{L}_{\omega,\lambda}(j, N) \le \exp\bigl(-\lambda^{\kappa_1}\bigr)$ such that
$$
\big| \langle\log | \lambda^{-1} E_j^{(N)}(\cdot,\omega) - \lambda^{-1} E| \rangle - \int_{\TT^2\setminus S} \log |\lambda^{-1} E_j^{(N)}(x,y,\omega) - \lambda^{-1} E|dx\, dy \big| \le N^{-\kappa}
$$
provided $\lambda^{-1} E \in [-B_0, B_0] \setminus
\bigl(L_{\omega,\lambda}(N)\cup \widetilde{L}_{\omega,\lambda}(j,
N)\bigr)$, $j = 1, 2,\dots, N$.  Hence, \be \label{eq:6.adf} \left|
\langle \log |f_N(\cdot,\lambda,\omega, E)| \rangle -
\int_{\TT^2\setminus S} \log |f_N(x,y,\lambda,\omega, E)|dx\, dy
\right| \le N^{1-\kappa} \ee whenever $\lambda^{-1} E \in [-B_0,
B_0] \setminus \cE_{\omega,\lambda}(N)$, where
$\cE_{\omega,\lambda}(N) = L_{\omega,\lambda}(N)\cup
\Bigl(\mybigcup_j \widetilde L_{\omega,\lambda}(j, N)\Bigr)$. The
lemma now follows  from \eqref{eq:6.ldv}, \eqref{eq:6.appr}, and
\eqref{eq:6.adf}.
\end{proof}

\begin{proof}[Proof of Theorem \ref{thm:1.2}]
Let $\lambda_0$ be as in Lemma~\ref{lem:6.1},
$|\lambda| \ge \lambda_0$, $\ell = \bigl[\lambda^{1/2}\bigr]$ and
let $\cE_{\omega,\lambda}(\ell)$ be the set in the statement of Lemma~\ref{lem:6.1}, then
$$
{1\over \ell} \langle\log |f_\ell(\cdot,\lambda,\omega, E)|\rangle > {1\over 2} \log |\lambda|
$$
for any $E \in \bigl[-B_0(V) - 2, B_0(V) +2\bigr]\setminus \cE_{\omega,\lambda}(\ell)$.
By Proposition~\ref{prop:4.4} there exists $\overline{\cF}_{\omega,\lambda}(\ell)$ with
 $\mes\overline{\cF}_{\omega,\lambda}(\ell) \le \exp\bigl(-\ell^\kappa\bigr)$ such that for any
$$
\big| \ell^{-1} \langle \log |f_\ell(\cdot,\lambda,\omega,E)| \rangle - L(\lambda,\omega,E)\big| < \ell^{-1/2}\ .
$$
Thus $L(\lambda,\omega, E) > 1/3 \log |\lambda|$ for any $E \in
\bigl[-\lambda B_0(V) - 2, \lambda B_0(V) + 2\bigr] \setminus
\cE_{\omega,\lambda}$, where $\cE_{\omega,\lambda} =
\cE_{\omega,\lambda}(\ell) \cup
\overline{\cF}_{\omega,\lambda}(\ell)$.  That proves the first part
in Theorem~\ref{thm:1.2}. The second part now follows from
Theorem~\ref{th:1.1}.
\end{proof}

\setcounter{section}{1}
\renewcommand{\thesection}{\Alph{section}}
\section*{Appendix A: Quantitative Ergodic Theorem}
\label{sec:A}
\setcounter{theorem}{0}

In this section, let $\psi\in C^4\left(\TT^2\right)$. We first derive a quantitative ergodic
theorem for the shift $T_\omega :
\TT^2\rightarrow\TT^2$ which is defined as $T_\omega(x_1,x_2) = (x_1 + \omega_1, x_2
+ \omega_2)$. We recall from Definition~\ref{def:Dioph} that $\omega$ is said to satisfy a
{\em Diophantine condition} provided
\begin{equation}
\label{eq:diophA} \| k_1\omega_1 + k_2\omega_2\| > c\left(|k_1| +
|k_2|\right)^{-A}\text{\ \ for all\ }(k_1,k_2)\in\IZ^2\setminus\{
0\}
\end{equation}
where $c > 0$, $A > 2$ are fixed.

\begin{prop}
\label{prop:A.ineq} Suppose $\tom=(\tom_1,\tom_2)$ satisfies the Diophantine assumption
\eqref{eq:diophA}. Then for all $\omega=(\omega_1,\omega_2)$ with $|\tom-\omega|<\frac{c}{2}N^{-1}$
\begin{equation}\label{eq:psi_av}
\left|\frac{1}{N}\sum_{m=1}^{N}\psi\left(T_\omega^mx\right) - \la\psi\ra\right|\lesssim B_4(\psi)N^{-\sigma}
\end{equation}
for all $x\in\TT^2$, $N\geq N_0(c,A)$.  Here, $\sigma={1\over 1+A}$.
\end{prop}

\begin{proof} The Fourier series
\[
\psi\left(T_\omega^mx\right)  = \sum_{k_1,k_2}c(k_1,k_2)e\bigl(k_1(x_1 +
m\omega_1) + k_2(x_2 + m\omega_2)\bigr) \] where \[c(k_1,k_2)  =
\int_{\TT^2}\psi(x)e(-k_1x_1 - k_2x_2)\,dx\]
 converges absolutely. Indeed, integrating by parts yields
\begin{align*}
\bigl|c(k_1,k_2)\bigr| & \lesssim B_4(\psi)\left(1 + |k_1|^2\right)^{-1}\left(1 + |k_2|^2\right)^{-1} \\
c(0,0) & = \int_{\TT^2}\psi(x)dx = \la\psi\ra \\
\frac{1}{N}\sum_{m=1}^N\psi\left(T_\omega^mx\right) - \la\psi\ra & =
\hspace*{-10pt}\sum_{(k_1,k_2)\neq 0}\hspace*{-8pt}c(k_1,k_2)e(k_1x_1 + k_2x_2)\frac{1}{N}
\sum_{m=1}^Ne\bigl(m(k_1\omega_1 + k_2\omega_2)\bigr).
\end{align*}
Recall the simple bound \[\left|\frac{1}{N}\sum_{m=1}^Ne(m\theta)\right|  \leq
\frac{2}{1 + N\|\theta\|}\]
Moreover,  for all $|k|\le N^\sigma$ with $\sigma=\frac{1}{1+A}$,
\begin{align*}
 \|k_1\omega_1 + k_2\omega_2\| &\ge \|k_1\tom_1 + k_2\tom_2\| - (|k_1|+|k_2|)|\omega-\tom| \\
&\ge c|k|^{-A} - |k|\frac{c}{2N} > \frac{c}{2}|k|^{-A}
\end{align*}
>From this we conclude that
\begin{align*}
\left|\frac{1}{N}\sum_{m=1}^N\psi\left(T_\omega^mx\right) - \la\psi\ra\right| & \leq  \hspace*{-8pt}\sum_{(k_1,k_2)\neq 0}
\hspace*{-8pt}\bigl|c(k_1,k_2)\bigr|\left|\frac{1}{N}\sum_{m=1}^Ne\bigl(m(k_1\omega_1 + k_2\omega_2)\bigr)\right| \\
& \lesssim \hspace*{-8pt}\sum_{(k_1,k_2)\neq 0}\hspace*{-8pt}B_4(\psi)\left(1 + |k_1|^2\right)^{-1}\left(1 + |k_2|^2\right)^{-1}
\bigl(1 + N\|k_1\omega_1 + k_2\omega_2\|\bigr)^{-1} \\
& \leq \hspace*{-8pt}\sum_{(k_1,k_2)\neq 0}\hspace*{-8pt}B_4(\psi)\left(1 + |k_1|^2\right)^{-1}\left(1 + |k_2|^2\right)^{-1}\bigl(1 + Nc(|k_1| + |k_2|)^{-A}\bigr)^{-1} \\
& \leq \hspace*{-7pt}\sum_{\substack{|k_1|\leq N^{\sigma}\\ |k_2|\leq N^{\sigma}}}
\hspace*{-9pt}B_4(\psi)\left(1 + |k_1|^2\right)^{-1}\left(1 + |k_2|^2\right)^{-1}\left(1 + cN^{\frac{1}{1+A}}\right)^{-1} \\
& \hspace*{17pt} + \hspace*{-18pt}
\sum_{\substack{|k_1| > N^{\sigma}\\ \text{OR } |k_2| > N^{\sigma}}}
\hspace*{-12pt}B_4(\psi)\left(1 + |k_1|^2\right)^{-1}\left(1 + |k_2|^2\right)^{-1} \\
& \lesssim B_4(\psi)N^{-\sigma}.
\end{align*}
as claimed.
\end{proof}

\begin{remark}\label{rem:weak_dioph}
Inspection of the previous proof reveals that Proposition~\ref{prop:A.ineq}
only requires the weaker condition
\begin{equation}\label{eq:NDioph} \|k\cdot\omega\|=\|k_1\omega_1+k_2\omega_2\|\ge N^{-\gamma_1} \quad \forall \; 1\le |k|\le N^{\gamma_2}
\end{equation}
were $\gamma_1,\gamma_2>0$ are any small numbers (but fixed; one then also needs to change the power
on the right-hand side of~\eqref{eq:psi_av} according to the choice of these constants).
Note that this is stable under perturbations $|\omega-\tilde\omega|<N^{-1}$. We refer to $\omega$ as in~\eqref{eq:NDioph}
as $(N,\gamma_1,\gamma_2)$-Diophantine.  Proposition \ref{prop:A.ineq} holds for $\w$ being $(N,\gamma_1,\gamma_2)$-Diophantine, with $\sigma={1\over 2}\min{\gamma_1,\gamma_2}$.
\end{remark}

\begin{lemma}\label{lem:A.sieb}
Let $0<\mu<\frac12$. Then
\begin{itemize}
\item for any $k\in\IZ^2$ and $\theta\in\TT$, $\mes\left\{\omega\in\TT^2 \::\: \|k\cdot\omega+\theta\| < \mu\right\} = 2\mu$
\item for any $\theta\in\TT$, $N_0\ge1$,  \[\mes\left\{\omega\in\TT^2 \::\: \min_{1\leq |k_1|,\,
|k_2|\leq N_0}\|k\cdot\omega+\theta\| < \mu\right\}\lesssim N_0\mu\]
\end{itemize}
\end{lemma}

\begin{lemma}\label{lem:A.acht}
For any $\omega,\omega_0\in\TT^2$,
$ \left|\min_{1\leq |k|\leq N_0}\|k\cdot\omega\| - \min_{1\leq |k|\leq N_0}\|k\cdot\omega_0\|\right|\les N_0|\omega - \omega_0| $
\end{lemma}

\begin{lemma}\label{lem:A.neun}
Given $\oN \in\IN$, $j=(j_1,j_2)\in [1,\oN]^2\subset \IZ^2$ set $\omega_j^{(\oN)} = (\frac{j_1}{\mathstrut{\oN} },\frac{j_2}{\mathstrut{\oN}})$. Then for any $N_0\in\IN$,
$\mu>0$ one has
\begin{equation}\label{eq:A.J}
\#\left\{ 1\leq j\leq\oN  \::\: \min_{1\leq |k|\leq N_0}\|k\cdot\omega_j^{(\oN)}\| < \mu\right\} \lesssim  \mu N_0^2\oN^2 + N_0^3 \oN
\end{equation}
Denote the set on the left-hand side by $J(\oN,N_0,\mu)$. Then
\[
\min_{1\le|k|\le N_0} \|k\cdot\omega\| > \mu/2
\]
for any $|\omega-\omega_j^{(\oN)}|<\frac{\mu}{2N_0}$ with $j\not\in J(\oN,N_0,\mu)$.
\end{lemma}

\begin{proof}
If $j=(j_1,j_2)\in J(\oN,N_0,\mu)$, then due to Lemma~\ref{lem:A.acht} one has
\[ \min_{1\le|k|\le N_0} \|k\cdot \omega\| \le \mu + \frac{N_0}{\oN} \]
for any
\begin{equation}\label{eq:calP}
\omega\in \cP_j^{(\oN)} = (j_1/\oN,(j_1+1)/\oN)\times (j_2/\oN,(j_2+1)/\oN).
\end{equation}
Hence,
\begin{equation}\label{eq:A.P}
\mes\Big[ \bigcup_{j\in J(\oN,N_0,\mu)} \cP_j^{(\oN)}  \Big] \le \mes\Big\{\omega\in\TT^2\::\: \min_{1\le|k|\le N_0} \|k\cdot \omega\| \le \mu + \frac{N_0}{\oN} \Big\} \lesssim N_0^2 (\mu+N_0/\oN)
\end{equation}
due to Lemma~\ref{lem:A.sieb}. The bound \eqref{eq:A.J} now follows from \eqref{eq:A.P}.
The second assertion follows from Lemma \ref{lem:A.acht}.
\end{proof}

\begin{corollary}\label{cor:A.10}
Given $N>0$ there exists $\cJ(N)\subset \{(j_1,j_2)\::\: 1\le j_1,j_2\le N^2\}$ such that (with $\cP_j^{(N^2)}$ as in~\eqref{eq:calP})
\begin{itemize}
\item[(1)] for any $\oN\ge N^2$ and any $\omega_j^{(\oN)}=(j_1/\oN,j_2/\oN)\in \bigcup_{k\not\in \cJ(N)} \cP_k^{(N^2)}$ one has
\[ \min_{1\le|k|\le N^{\kappa/4}} \|k \omega_j^{(\oN)}\|\ge N^{-\kappa} \]
\item[(2)] $\mes \Big[ \bigcup_{k\in\cJ(N)} \cP_k^{(N^2)} \Big] \le N^{-\kappa/2} $
\end{itemize}
Here $\kappa>0$ is a small constant.
\end{corollary}

\begin{proof} Using the notations of Lemma~\ref{lem:A.neun} set $\cJ(N)= J(N^2,N^{\frac{\kappa}{4}}, N^{-\kappa})$. If $\omega_j^{(\oN)}\in\bigcup_{k\not\in \cJ(N)} \cP^{(N^2)}_k$,
then $|\omega_j^{(\oN)}-\omega_k^{(N^2)}|\lesssim N^{-2}$ for some $k\not\in J(N)$. Therefore, (1), (2) follow from Lemma~\ref{lem:A.neun}.
\end{proof}

Let $\psi\in C^4\left(\TT^2\right)$. Let $\omega_1,\omega_2$ satisfy
\[
\|k\omega_i\|\geq c|k|^{-A}
\]
for any $k\in\IZ\setminus\{ 0\}$, where $c > 0$, $A > 2$.

To deal with the skew-shift, we need the following well-known estimate.

\begin{lemma}
\label{lem:A.pq}
Let $\alpha\in(0,1)$. If
\[
\left|\alpha - \frac{p}{q}\right|\leq\frac{1}{q^2}
\]
where $p,q\in\IN$, $(p,q) = 1$, then for any $N\in\IN$, one has
\[
\sum_{k=1}^N\min\left(N,\|k\alpha\|^{-1}\right)\lesssim\left(q + N +
N^2q^{-1}\right)\max\{ 1,\log q\}.
\]
\end{lemma}

\begin{proof}
This is Lemma~4.1~in~[Nat].
\end{proof}

Suppose $\alpha$ satisfies
\begin{equation}
\label{eq:A.kalpha}
\|k\alpha\| > c|k|^{-(1 + \ve)}\text{ for }k\in\IZ\setminus\{ 0\}
\end{equation}
where $0 < c < 1$ and $\eps>0$.
Let $\left\{\frac{p_s}{q_s}\right\}_{s=1}^{\infty}$ be the convergents of $\alpha$. Recall that
\[
\left|\alpha -
\frac{p_s}{q_s}\right|\leq\frac{1}{q_sq_{s+1}}\leq\frac{1}{q_s^2}
\quad\forall\, s\ge1.
\]
Note
\[
cq_{s - 1}^{-(1 + \ve)} < \|q_{s - 1}\alpha\| = |q_{s  - 1}\alpha - p_{s - 1}|\leq q_s^{-1}.
\]
Hence
\[
q_s < c^{-1}q_{s - 1}^{1 + \ve}.
\]
If $q_{s - 1} < N\leq q_s$, then $q_s < c^{-1}N^{1 + \ve}$. Combining this with Lemma~\ref{lem:A.pq}, one has

\begin{lemma}
\label{lem:A.omega}
Suppose $\omega$ satisfies (\ref{eq:A.kalpha}). Then for any $N\geq N_0(c,\ve)$
\[
\sum_{k=1}^N\min\left(N,\|k\alpha\|^{-1}\right)\lesssim c^{-1}N^{1 + \ve}\log N.
\]
\end{lemma}

We also need the following version of Weyl's inequality (a more
general version can be found in \cite{Nath}.)

\begin{lemma}
\label{lem:A.omega2} Suppose $\alpha$ satisfies (\ref{eq:A.kalpha}).
Let $S = \sum_{k=1}^Ne\left(k^2\alpha + k\beta\right)$ where
$N\in\IN$, $\beta\in\IR$. Then $|S|\lesssim N^{\frac{1}{2} + \ve}$
for $N\geq N_0(c,\ve)$.
\end{lemma}

\begin{proof}
\begin{align*}
|S|^2 & = \sum_{k=1}^Ne\left(k^2\alpha + k\beta\right)\sum_{\ell=1}^Ne\left(-\ell^2\alpha - \ell\beta\right) \\
& = \sum_{k=1}^N\sum_{\ell=1}^{N}\hspace*{-0pt}e\Bigl(\left(k^2 - \ell^2\right)\alpha + (k - \ell)\beta\Bigr) \\
& = \sum_{k=1}^N\sum_{m=k-N}^{k-1}\hspace*{-5pt}e\left(m(2k - m)\alpha + m\beta\right) \\
& = \hspace*{-8pt}\sum_{m=1-N}^{-1}\hspace*{-2pt}\sum_{k=1}^{N+m}e\bigl(m(2k - m)\alpha + m\beta\bigr) + \sum_{k=1}^Ne(0)
+ \sum_{m=1}^{N-1}\sum_{k=m+1}^N\hspace*{-4pt}e\bigl(m(2k - m)\alpha + m\beta\bigr) \\
& \leq \sum_{m=1}^{N-1}\left|\sum_{k=1}^{N-m}e\bigl(k(-2m\alpha)\bigr)\right| + N + \sum_{m=1}^{N-1}\left|\sum_{k=m+1}^Ne\bigl(k(2m\alpha)\bigr)\right| \\
& \lesssim N + \sum_{m=1}^N\min\left(N,\|m(2\alpha)\|^{-1}\right)
\end{align*}
\[
\|m(2\alpha)\| > c|2m|^{-(1 + \ve)} = \left(c2^{-(1 + \ve)}\right)|m|^{-(1 + \ve)}
\]
So $2\alpha$ satisfies~(\ref{eq:A.kalpha}) if $c$ is replaced by $c2^{-(1 + \ve)}$. By Lemma~\ref{lem:A.omega},
\[
|S|^2\lesssim N + c^{-1}2^{1 + \ve}N^{1 + \ve}\log N\lesssim N^{1 + 2\ve}.
\]
\end{proof}

\begin{remark} \label{rem:A.skewn}
We require $N\geq N_0(c,\ve)$ only to make sure that
\begin{itemize}
\item $\log c^{-1}N^{1 + \ve}\lesssim\log N$
\item $c^{-1}N^{1 + \ve}\log N\lesssim N^{1 + 2\ve}$.
\end{itemize}
Hence we can choose $N_0\asymp c^{-\frac{2}{\ve}}$.  Moreover, for
the purpose of Lemma \ref{lem:A.omega2}, there is no need in
condition \eqref{eq:A.kalpha}  for all $k \in \IZ \setminus \{0\}$.
Indeed, assume that \be\label{eq:A.skewn} \|k\omega\| \ge
cN^{-(1+\ve)}\quad \text{for all}\quad 0 < |k| \le N\ . \ee If
$q_{s-1} \le N < q_s$, then $\|q_{s-1}\omega \| \ge cN^{-(1+\ve)}$.
Since $\|q_{s-1}\omega\| \le q_s^{-1}$, that implies $q_s \le c^{-1}
N^{(1+\ve)}$ thus, $\big|\omega - \frac{p_s}{q_s}\big| \le
q_s^{-2}$, and $N \le q_s \le c^{-1} N^{(1+\ve)}$.  Therefore,
Lemma~\ref{lem:A.omega} holds and Lemma~\ref{lem:A.omega2} follows.
We denote by $\TT_{c,\ve,N}$ the set of all $\omega$ for which
\eqref{eq:A.skewn} is valid (see Appendix B).
\end{remark}

Now suppose $T : \TT^2\rightarrow\TT^2$ is a skew-shift $T(x_1,x_2)
= (x_1 + x_2,x_2 + \omega)$ with $\alpha = \frac{\omega}{2}$
satisfying~(\ref{eq:A.kalpha}).

\begin{prop}\label{prop:A.skewineq}
\[
\left|\frac{1}{N}\sum_{m=1}^N\psi\left(T^mx\right) - \la\psi\ra\right|\lesssim B_4(\psi) N^{-\frac{\ve}{2}}
\]
for all $x\in\TT^2$, $N\geq N_1(c,\ve)$.
\end{prop}

\begin{proof}
\[
T^mx = \left(x_1 + mx_2 + \frac{m(m - 1)}{2}\omega, x_2 + m\omega\right)
\]
\begin{align*}
\frac{1}{N}\sum_{m=1}^N\psi\left(T^mx\right) - \la\psi\ra & = \hspace*{-10pt}\sum_{(k_1,k_2)\neq 0}\hspace*{-10pt}c(k_1,k_2)e(k_1x_1 + k_2x_2)\frac{1}{N}\sum_{m=1}^Ne\biggl(m^2k_1\frac{\omega}{2} + m\left(k_1x_2 - k_1\frac{\omega}{2} + k_2\omega\right)\biggr) \\
\left|\frac{1}{N}\sum_{m=1}^N\psi\left(T^mx\right) - \la\psi\ra\right| & \leq
\sum_{k_2\neq 0}\bigl|c(0,k_2)\bigr|\frac{1}{N}\left|\sum_{m=1}^Ne(mk_2\omega)\right| \\
& \hspace*{17pt} + \sum_{\substack{k_1\neq 0\\ k_2\in\IZ}}\bigl|c(k_1,k_2)\bigr|\frac{1}{N}\left|\sum_{m=1}^Ne\biggl(m^2k_1\frac{\omega}{2} + m\left(k_1x_2 - k_1\frac{\omega}{2} + k_1\omega\right)\biggr)\right| \\
& \lesssim \sum_{k_2\neq 0} B_4(\psi)\left(1 + |k_2|^2\right)^{-1}\bigl(1 + N\|k_2\omega\|\bigr)^{-1}+\sum_{k_2\in\IZ}B_4(\psi)\left(1 + |k_2|^2\right)^{-1} \\
& \hspace*{21pt} \cdot \sum_{k_1\neq 0}\left(1 + |k_1|^2\right)^{-1}\frac{1}{N}\left|\sum_{m=1}^Ne\biggl(m^2k_1\frac{\omega}{2} + m\left(k_1x_2 - k_1\frac{\omega}{2} + k_2\omega\right)\biggr)\right|
\end{align*}

As in the proof of Proposition~(\ref{prop:A.ineq}),
\[
\sum_{k_2\neq 0}B_4(\psi)\left(1 + |k_2|^2\right)^{-1}\left(1 + N\|k_2\omega\|\right)^{-1}\lesssim B_4(\psi)N^{-\frac{1}{2(1 + \ve)}}
\]
provided $N\geq N_1(c,\ve)$.

We want to apply Lemma~\ref{lem:A.omega2} to estimate the second sum. Note that $\alpha = k_1\frac{\omega}{2}$ satisfies~(\ref{eq:A.kalpha}) if $c$ is replaced by $c|k_1|^{-(1+\ve)}$. For $N\geq\left[c|k_1|^{-(1+\ve)}\right]^{-\frac{2}{\ve}}$, $\left|\sum_{m=1}^Ne\left(m^2k_1\frac{\omega}{2} + m\beta\right)\right|\lesssim N^{\frac{1}{2} + \ve}$.
\begin{align*}
\sum_{k_1\neq 0}\left(1 + |k_1|^2\right)^{-1}\frac{1}{N}\left|\sum_{m=1}^Ne\left(m^2k_1\frac{\omega}{2} + m\beta\right)\right| & \lesssim \hspace*{-10pt}\sum_{|k_1| < c^{-1}N^{\frac{\ve}{2}}}\hspace*{-12pt}\left(1 + |k_1|^2\right)^{-1}N^{-\frac{1}{2} + \ve} + \hspace*{-10pt}\sum_{|k_1|\geq c^{-1}N^{\frac{\ve}{2}}}\hspace*{-12pt}\left(1 + |k_1|^2\right)^{-1} \\
& \lesssim cN^{-\frac{\ve}{2}}
\end{align*}
provided $N\geq N(c,\ve)$.
\end{proof}

\addtocounter{section}{1}
\setcounter{theorem}{0}
\section*{Appendix B: Metric  estimates for
 the typical rational approximation rate}

Given $\omega \in (0,1)$ denote by $[a_1, a_2,
\ldots]=[a_1(\omega), a_2(\omega), \ldots ]$ its
continued fraction, $a_j \in \ZZ, a_j \geq 1$.
Take arbitrary integers $k_j \geq 1, j=1,2,
\ldots, n$ and put

\bigskip

\begin{equation*}
\label{eq:1.a}
\cE \left(\begin{array}{cccc}
1, & 2, & \ldots, & n\\k_1, & k_2, & \ldots, & k_n\\
\end{array}\right )= \{ \omega \in (0,1):
a_j(\omega)=k_j, 1 \leq j \leq n\}
\end{equation*}

\bigskip

Denote by $\omega_s = p_s / q_s$ the convergents
for $\omega=[a_1, a_2, \ldots ]$.  Recall that

\bigskip

\begin{equation*}
\label{eq:2.a}
\omega_s=\frac{p_s}{q_s}=\frac{1}{a_1 + \frac{1}{
\begin{array}{cc} a_2\ + & \\ & \ddots + \frac{1}{a_s} \end{array} }}
\end{equation*}

\bigskip

\begin{equation*}
\begin{gathered}
\label{eq:3.a}
q_s=a_s q_{s-1} + q_{s+2}, \quad q_0=1,
q_{-1}=0\\
p_s=a_s p_{s-1} + p_{s-2}, \quad p_0=0,
p_{-1}=1
\end{gathered}
\end{equation*}

\bigskip

The set $\cE \left(\begin{array}{cccc}
1, & 2, & \ldots, & n\\k_1, & k_2, & \ldots, & k_n\\
\end{array}\right )$
consists of an interval,

\bigskip

\begin{equation*}
\label{eq:4.a}
\cE \left(\begin{array}{cccc}
1, & 2, & \ldots, & n\\k_1, & k_2, & \ldots, & k_n\\
\end{array}\right )
=\Big[\frac{p_n}
{q_n}, \frac{p_n+p_{n-1}}{q_n+q_{n-1}} \Big],
\qquad {p_s\over q_s}=[k_1, \ldots, k_s]
\end{equation*}

\bigskip

In particular,

\bigskip

\begin{align*}
\label{eq:5.a}
& \Bigg|\cE \left(\begin{array}{cccc}
1, & \ldots, & n, & n+1\\k_1, & \ldots, & k_n, & k\\
\end{array}\right )\Bigg| \\
& \qquad = k^{-2}(1+(kq_n)^{-1} q_{n-1})^{-1}
(1+k^{-1}+(k q_n)^{-1} q_{n-1}) (1+ q_n^{-1}
q_{n-1})
\Bigg |
\cE \left(\begin{array}{cccc}
1, & 2, & \ldots, & n\\k_1, & k_2, & \ldots, & k_n\\
\end{array}\right )
\Bigg | \ ,
\end{align*}

\bigskip

\begin{equation*}
\label{eq:6.a}
{1\over 3k^2} \Bigg|\cE \left(\begin{array}{cccc}
1, & 2, & \ldots, & n\\k_1, & k_2, & \ldots, & k_n\\
\end{array}\right )
\Bigg | \leq \Bigg |
\cE \left(\begin{array}{cccc}
1, & \ldots, & n, & n+1 \\k_1, & \ldots, & k_n, & k\\
\end{array}\right ) \Bigg|
\leq {2\over k^2} \Bigg |
\cE \left(\begin{array}{cccc}
1, & 2, & \ldots, & n\\k_1, & k_2, & \ldots, & k_n\\
\end{array}\right )
\Bigg |
\end{equation*}

\bigskip Hence ${1\over 3k^2} \leq {\rm mes} \{ \omega :
a_{n+1}(\omega)=k \} \leq {2\over k^2}$ and

\begin{equation}
\label{eq:8.a}
\mathop\prod\limits_{j=1}^r {c\over k_j} \leq
{\rm mes} \{ \omega : a_{n_{j}}(\omega) \geq
k_j, j=1,2 \ldots, r \} \leq
\mathop\prod\limits_{j=1}^r {C\over k_j}
\end{equation}

\bigskip

\begin{prop}
\label{1.a}
\begin{equation*}
\label{eq:9.a}
{\rm mes} \{ \omega \in (0,1):
\max\limits_{t \leq s} \ a_t \geq K\}
\leq CsK^{-1}
\end{equation*}

\bigskip

\begin{equation*}
\label{eq:10.a}
{\rm mes} \{ \omega \in (0,1): \omega_s= p_s
q_s^{-1}, \log q_s > s_0 K\ \text{for some}\
 s \leq s_0 \} \leq C s_0 e^{-K}
\end{equation*}

\end{prop}

\bigskip

\begin{proof}
Due to (\ref{eq:8.a})

\begin{equation*}
\label{eq:11.a}
{\rm mes} \{ \omega \in (0,1):
\max\limits_{t \leq s}\ a_t \geq K \}
\leq \mathop\sum\limits_{1 \leq t \leq s} \
{\rm mes} \{ \omega: a_t(\omega) \geq K \} \leq Cs
K^{-1}
\end{equation*}
Note that $q_t \leq (a_t +1) q_{t-1}$.
If $\log q_s > K$ for some $s \leq s_0$,
then $\max\limits_{1 \leq t \leq s_0}\ \log
(a_t (\omega)+1) \geq K$.  Therefore

\begin{equation*}
\label{eq:14.a}
{\rm mes} \{ \omega \in (0,1):
\omega_s=p_sq_s^{-1}, \log q_s \geq K\ \text{
for some}\ s \leq s_0 \} \leq C s_0 e^{-K}
\end{equation*}

\end{proof}

\bigskip

Note also that (\ref{eq:8.a})implies for arbitrary $k(n)$

\begin{equation*}
\label{eq:15.a}
{\rm mes} \{ \omega : a_n(\omega) \leq k(n),
n=n_0+1, \ldots n_0+m \} \leq
\mathop\prod\limits_{n=n_0+1}^{n+m} \Big( 1 -
\frac{c}{k(n)} \Big)
\end{equation*}

In particular

\begin{equation*}
\label{eq:16.a}
{\rm mes} \{ \omega : a_n(\omega) \leq k(n),
n=n_0+1, n_0+m \} \leq \exp \Big( -
\frac{cm}{k(n_0)}\Big )
\end{equation*}

\bigskip

\begin{equation*}
\begin{gathered}
\label{eq:17.a}
{\rm mes} \bigl\{ \omega : a_n(\omega) \leq n_0 \log
n_0(\log \log n_0)(\log \log \log n_0),
n=n_0+1, \ldots n_0+m \bigr\}\\[16pt]
\leq {\rm exp}(-C(\log \log n_0))\\
\end{gathered}
\end{equation*}
provided $m> n_0(\log n_0) (\log \log n_0)^2$.

\bigskip

\begin{lemma}\label{lem:B.2}
Let $\omega_s = \frac{p_s}{q_s}$ be the
convergents for $\omega$.  Take $1 \leq m <
q_s, s \geq 1$.  Then
$\Vert m \omega \Vert \geq \frac{a_{s+1}}{q_{s+1}}$.
\end{lemma}

\begin{proof}
By induction in $s$.  If $s=1$, then can assume
that $a_1=q_s \geq 2$.  Since $\frac{a_2}{q_2}
= \omega_2 < \omega < \omega_1 =
\frac{1}{a_1}$, and $1 \leq m \leq a_1 -1$ one
has $\frac{a_2}{q_2} \leq m \omega < 1 -
\frac{1}{a_1}$.  Thus $\Vert m \omega \Vert
\geq {\rm min} (\frac{a_2}{q_2}, \frac{1}{a_1})
\geq \frac{a_2}{q_2}$, which is $(1)$ for
$s=1$.  Let $s \geq 2$ and set $m=a q_{s-1} +
b$, $0 \leq b < q_{s-1}$, $ a < \frac{q_s}{q_{s-1}}
\leq a_s +1$.  Thus $a \leq a_s$.

\ \newline {\bf Case 1}
$\quad a=0$.  Then $m < q_{s-1}$ and by the reductive
hypothesis,

$$
\Vert m \omega \Vert \geq \frac{a_s}{a_s} \geq
\frac{a_{s+1}}{q_{s+1}}
$$
The final inequality holds since $q_{s+1}
\geq a_{s+1} q_s$ and $a_s \geq 1$.

\ \newline {\bf Case 2}
$\quad 1 \leq a < a_s, b \not= 0$.  Then

$$
\Vert m \omega \Vert \geq \Vert b \omega \Vert
- \Vert a q_{s-1} \omega \Vert \geq \frac
{a_s}{q_s} - a \Vert q_{s-1} \omega \Vert
$$

$$
\geq \frac{a_s}{q_s} - a \frac{1}{q_s} \geq
\frac{1}{q_s} \geq \frac{a_{s+1}}{q_{s+1}}
$$

\ \newline {\bf Case 3}
$\quad m=a q_{s-1}, 1 \leq a \leq a_s$.  (i.e. $b=0$)  One has

\begin{equation*}
(-1)^{s-1} \omega_{s+1} < (-1)^{s-1} \omega <
(-1)^{s-1} \omega_s\ \text{for}\ s=1,2, \ldots
\end{equation*}
Hence, for $s \geq 2$, and $a$ as above,

\begin{equation*}
\underbrace{(-1)^{s-1} a q_{s-1} \omega_{s-1}}_{\in
\ZZ} < (-1)^{s-1}a q_{s-1} \omega_{s+1} < (-1)^{s-1}
a q_{s-1} \omega < (-1)^{s-1} a q_{s-1} \omega_s.
\end{equation*}
where

\begin{equation*}
|aq_{s-1} \omega_{s-1} - aq_{s-1} \omega_s|=
\frac {a}{q_s} \leq \frac{a_s}{q_s} <1.
\end{equation*}
Therefore,

\begin{equation*}
\Vert aq_{s-1} \omega \Vert \geq {\rm min} (|
aq_{s-1} \omega_{s+1} - aq_{s-1} \omega_{s-1}|,
1-|aq_{s-1} \omega_s -aq_{s-1} \omega_{s-1}|)
\end{equation*}
On one hand

$$
|aq_{s-1} \omega_{s+1} - aq_{s-1}
\omega_{s-1}|= aq_{s-1} \Bigg | \frac{p_{s+1}}{q_{s+1}}
- \frac{p_{s-1}}{q_{s-1}}\Bigg |
= a\frac{a_{s+1}}{q_{s+1}} \geq \frac{a_{s+1}}{q_{s+1}},
$$
and on the other hand,

\begin{gather*}
1 - |a q_{s-1} \omega_s - a q_{s-1}
\omega_{s-1}|=1-a q_{s-1}
\frac{1}{q_{s}q_{s-1}}=1 - \frac{a}{q_s}\\[16pt]
\geq \frac{q_s - a_s}{q_s} =
\frac{a_s(q_{s-1}-1) + q_{s-2}}{q_s} \geq \frac
{1}{q_s} \geq \frac{a_{s+1}}{q_{s+1}}
\end{gather*}

\ \newline {\bf Case 4}
$\quad m=a_s q_{s-1} + b < q_s, 1 \leq b <a_{s-2}$  Then

\begin{gather*}
\Vert m \omega \Vert \geq \Vert (q_s - m)
\omega \Vert - \Vert q_s \omega \Vert = \Vert
\overbrace{(q_{s-2}-b)}^{<q_{s-2}}\omega
\Vert - \Vert q_s \omega \Vert\\[16pt]
\geq \frac{a_{s-1}}{q_{s-1}} -
\frac{1}{q_{s-1}} \geq \frac{1}{q_{s-1}} -
\frac{1}{q_{s+1}}=
\frac{a_{s+1}q_s}{q_{s-1}q_{s+1}} \geq \frac
{a_{s+1}}{q_{s+1}}.
\end{gather*}

\end{proof}

Given $N \in \IN$, $\ve > 0$ let
$$
\cJ(N,\ve) = \left\{\omega \in [0,1]: q_{s-1}(\omega) \le N < q_s(\omega)
\quad\text{for some $s$ and $a_s(\omega) > N^\ve$}\right\}
$$

\begin{lemma}\label{lem:B.elem}
$\mes\cJ(N,\ve) \le N^{-\ve/2}$, provided $N > N_0(\ve)$.
\end{lemma}

\begin{proof} Recall that $q_{s-1}(\omega) \ge 2^{{s-1\over 2}}$.  Therefore $q_{s-1}(\omega) \le N$ implies $s \lesssim \log N$.  Hence.
\begin{align*}
\mes \cJ(N, \ve) & \le \sum_{s\lesssim\log N}\, \mes\left\{\omega \in [0,1]: a_s(\omega) > N^\ve\right\}\\
& \lesssim N^{-\ve} \log N\
\end{align*}
due to Proposition~\ref{1.a}
\end{proof}

\begin{defi}\label{def:B.typ} Given $N \in \IN$, and $0 < c,\ve \ll 1$.  Let
$$
\TT_{c,\ve, N} = \left\{\omega \in \TT: \|k\omega\| \ge cN^{-(1+\ve)}\ \quad \text{for any}\quad 0 < |k| \le N\right\}\ .
$$
\end{defi}

\begin{lemma}\label{lem:B.3}
$\TT \setminus \TT_{c,\ve, N} \subset \cJ(N, \ve)$, provided $N > N_0(\ve)$.  In particular, $\mes\bigl(\TT\setminus\TT_{c,\ve, N}\bigr) \le N^{-\ve/2}$.
\end{lemma}

\begin{proof} Assume that $\|k\omega\| < cN^{-(1+\ve)}$ for some $\omega \in \TT$, $0 < k \le N$.  Find $s$ so that $q_{s-1}(\omega) \le N < q_s(\omega)$.  Then $k < q_s(\omega)$ and by Lemma~\ref{lem:B.2}, $\|k\omega\| \ge {1\over 2 q_s(\omega)}$.  Hence, $q_s(\omega) \ge (2c)^{-1} N^{1+\ve}$.  Therefore
$$
a_s(\omega) +1 \ge q_s(\omega)/q_{s-1}(\omega) \ge (2c)^{-1} N^{1+\ve}\, \big /\, N \ge c^{-1} N^\ve > N^\ve + 1\ .
$$
\end{proof}

Given $\oN$ set $\omega_j^{(\oN)} = j/\oN$
\be\label{eq:B.P1}
\cP_j^{(\oN)} = \left(\omega_j^{(\oN)}, \omega_{j+1}^{\oN}\right)\ ,\quad j = 1, 2,\dots
\ee

\begin{corollary}\label{cor:B.metn}
Let $N, \oN \in \IN$, $\oN > N^3$, $0 < \ve, c\ll 1$.  Then
\be\label{eq:B.estn}
\#\left\{1 \le j \le \oN: \min_{1 \le |k| \le N} \big\|k\omega_j^{(\oN)}\big\| < c N^{-(1+\ve)}\right\} \lesssim \oN\, N^{-\ve/2}\ .
\ee
Denote the set on the left-hand side of \eqref{eq:B.estn} by $\cJ(\oN, N, c, \ve)$, then
$$
\min_{1 \le |k| \le N}\, \|k\omega\| \ge {c\over 2} N^{-(1+\ve)}
$$
for any $|\omega - \omega_j^{(\oN)}| < 1/ \oN$, with $j \notin \cJ(\oN, N, c, \ve)$.
\end{corollary}

\begin{proof} If $j \in \cJ(\oN, N, c, \ve)$ then by Lemma~\ref{lem:A.acht}
$$
\min_{1 \le k \le N}\, \|k\omega\| \le N/\oN + cN^{-(1+\ve)} \le 2cN^{-(1+\ve)}
$$
for any $\omega\in \cP_j^{(\oN)}$.  Hence,
\be\label{eq:B.unimes}
\mes\Biggl(\bigcup_{j\in \cJ(\oN, N, c,\ve)}\, \cP_j^{(\oN)}\Biggr) \le \mes
\left(\TT \setminus \TT_{2c,\ve, N}\right) \le N^{-\ve/2}\ .
\ee
Estimate \eqref{eq:B.estn} follows from \eqref{eq:B.unimes}.  The second assertion follows from Lemma~\ref{lem:A.acht}.
\end{proof}

\addtocounter{section}{1}
\setcounter{theorem}{0}
\section*{Appendix C: $C^\alpha$-smooth potentials}

Here we discuss the modifications needed for the case of
$C^\alpha$-smooth potential, $0\le\alpha<1$, in
Theorems~\ref{th:1.1} and \ref{thm:1.2}. Let $\varphi(x),  x=(x_1,
x_2)\in\TT^2$  be a $C^\alpha$-smooth function, i.e.,
\begin{equation}\label{c.bal}
B_\alpha(\varphi):=\sup_{x\not=y} |x-y|^{-\alpha}
| \varphi(x)-\varphi(y)| < +\infty
\end{equation}
where $0<\alpha\le 1$. Given $\tau>0$, let $h_\tau(x)\in C^4 (\IR)$
be as in Definition~\ref{def:moll}, i.e., $h_\tau (x)$ is
$1$-periodic,
\begin{itemize}
\item $h_\tau\ge 0$
\item ${\rm{supp}}\ h_\tau \subset \bigcup_{k\in \cZ}
     [k-\tau,\ k+\tau]$
\item $\int^1_0 h_\tau (y)\, dy=1$
\item $\max_{y\in\IR} \Big|\left(\frac{\displaystyle d}{\displaystyle
dy}\right)^{m} h_\tau (y)\Big|
\lesssim \tau^{-(m+1)}\quad$ for $m\le 4$
\end{itemize}

Set $\tilde h_\tau (x_1, x_1)=h_\tau (x_1)h_\tau (x_2)$.

\begin{lemma} \label{c.1}
Define \label{c.mol}
\begin{equation}
\psi (x_1, x_2):=\int_{\tor^2}
\varphi(y_1, y_2) \tilde h_\tau (x_1 - y_1, x_2-y_2)\, dy
\end{equation}
Then $\psi \in C^4(\tor^2)$ satisfies
\begin{itemize}
\item[(1)] $\max_{x\in\tor^2} | \varphi(x)-\psi(x)| \le
B_\alpha (\varphi) \tau^\alpha$
\item[(2)] $B_4 (\psi)\lesssim B_0 (\varphi)\tau^{-4}$
\end{itemize}
\end{lemma}

\begin{proof} Note that
\begin{align*}
\psi (x_1, x_2)&=\int_{\|y_1-x_1\| \le \tau,\ \|y_2-x_2\| \le \tau}\
       \varphi(x_1-y_1,\ x_2-y_2)\tilde h_\tau(y_1, y_2)\, dy_1 dy_2\\
\varphi(x_1, x_2) &=\int_{\|y_1 -x_1\| \le \tau, \|y_2-x_2\| \le\tau}\
\varphi(x_1, x_2)\tilde h_\tau (y_1, y_2)\, dy_1 dy_2\\
\end{align*}

Therefore, (1) follows from (\ref{c.bal}).  Part (2)
follows just from the definition (\ref{c.mol}).
\end{proof}

Similarly to the $C^1$ case (see Corollary~\ref{cor:2.cprime}) we
proceed with the following.

\begin{corollary} \label{c.count}
Let $f\in C^\alpha(\tor^2),\ T_\omega (x)=x+\omega$. Let
$N\in \IN$ be large and let $\omega$ be $(N, \gamma_1, \gamma_2)$-Diophantine.
Then for any $\xi \in\IR$, $0<\delta <1$ one has
$${1\over N}\# \{1\le k\le N:T_\omega^k x\in S_f (\xi, \delta)\}
\lesssim \mes S_f (\xi, 2\delta) + (1+B_\alpha
(f))\delta^{1\over2}$$ provided $N>\delta^{-20\over \alpha\sigma}$.
Here $\sigma={1\over 2} \min (\gamma_1, \gamma_2)>0$ is a constant,
$S_f (\xi, \delta')=\{x\in\tor^2: |f(x)-\xi|<\delta'\}$
\end{corollary}

\begin{proof}
Using the notations of Lemma \ref{lem:2.list}, we have to estimate
${1\over N} \sum^n_{k=1} \chi_{\delta}(f(T^k_x)-\xi)$. Note that
$\varphi(x)=\chi_\delta (f(x)-\xi)$ is $C^\alpha$-smooth with
$B_\alpha(\varphi) \lesssim \delta^{-1} B_\alpha(f)$. Define
$\psi(x)$ as in Lemma \ref{c.1}, then, due to (\ref{c.mol}),
$|\langle\psi\rangle-\langle\varphi\rangle| \le B_\alpha(\varphi)
\tau^{\alpha}\le \delta^{-1} B_\alpha (f)\tau^{\alpha}$ and \[\max_x
\Big|{1\over N} \sum^n_{k=1} \varphi (T^m x) -{1\over N} \psi (T^m
x)\Big| \le \delta^{-1} B_\alpha (f) \tau^{\alpha}.\] As in
Corollary~\ref{cor:2.cprime}, one obtains
\begin{align*}
{1\over N}\# \{1\le k\le N:T^k_\omega x \in S_f(\xi, \delta)\}
& \le {1\over N} \sum^n_{k=1} \varphi (T^k_w x)\\
&\le {1\over N} \sum^n_{k=1} \psi (T^k_\omega x) +\delta^{-1}
B_\alpha (f) \tau^{\alpha} \\
& \le \langle\psi\rangle +B_4 (\psi) N^{-\sigma}+ \delta^{-1} B_\alpha(f) \tau^\alpha \\
& \le \langle\varphi\rangle +B_0(f)
\tau^{-5} N^{-\sigma}+2\delta^{-1} B_\alpha (f) \tau^{\alpha}
\end{align*}
due to Proposition \ref{prop:A.ineq} and Remark \ref{rem:weak_dioph}. The assertion
follows if we take here $\tau=\delta^{3\over 2\alpha}$.
\end{proof}

One can see that the rest of the auxiliarly assertions needed for
the proof of Theorem~\ref{thm:A} do not rely on the smoothness of
the function $f(x)$, and therefore does not require any
modifications. Thus, Theorem~\ref{thm:A}  as well as the remarks
after it hold for $f\in C^\alpha (\tor^2)$.

The modifications needed for the rest of the Section~\ref{sec:2} and
whole Section~\ref{sec:3} consist only in stronger restrictions on
the interval in which $\omega$ runs. For instance, the assertions in
Theorem~\ref{thm:3.6} are valid, provided
$|\omega_1-\omega_0|<(1+B_\alpha (f))^{-1} N^{-{4\over \alpha}}$.

Section 4, 5, 6 rely only on the application of Theorem~\ref{thm:A}
and Theorem~\ref{thm:3.6} to $\log|E^{(N)}_j (x, w)-E|$, where
$E_1^{(N)} (x, \omega) <\cdots < E^{(N)}_N (x, \omega)$ are the
eigenvalues of $H_N (x, \omega)$. The only fact needed for the
validity of these applications is as follows:

\begin{lemma}\label{c.rel}
Suppose $V(x)\in C^\alpha (\tor^2)$, let $T=T_\omega$ be the shift
(or the skew-shift).  Let $E_1^{(N)} (x, \omega) <E^{(N)}_2 (x,
\omega)< \cdots < E_N^{(N)} (x, \omega)$ be the eigenvalues of
$H_N(x , \omega)$.  Then the functions $E_j^{(N)}$ are
$C^\alpha$--smooth and $B_\alpha \bigl(E_j^{(N)} \bigr) \lesssim N^2
B_\alpha (V) , j=1, 2,\cdots, N$.
\end{lemma}

\begin{proof}
Note that
$$\|H_N (x, \omega)-H_N (\tilde x, \tilde \omega)\|
\lesssim B_\alpha (V) N^2 (|x-\tilde x|^\alpha
+|\w-\tilde \w|^\alpha)$$

Recall that due to the minimax principle, if $A_i,\ i=1,2$, are
Hermitian operators in $\IC^N$, and $\ E_{1, i}<E_{2, i}, \cdots$
are the eigenvalues of $A_i$, $i=1,2$ , then
$$|E_{j, 1}, - E_{j, 2}| \le
\| A_1 - A_2\|$$ for $j= 1,2,...,N$, and the lemma follows.
\end{proof}


\begin{thebibliography}{GolMolPasxxx}


\bibitem[Bha]{Bhat} Bhatia, R. {\em Perturbation bounds for matrix eigenvalues.} Pitman research notes in mathematics series~162, Longman, 1987.


\bibitem[Bje]{Bje} Bjerkl\"ov, K. {\em Positive Lyapunov exponents for continuous quasiperiodic Schr\"odinger equations.}
  J.\ Math.\ Phys.~47,  no.~2,  (2006)

\bibitem[Bou]{Bou} Bourgain, J. {\em Green's function estimates for lattice Schr\"odinger
operators and applications.} Annals of Mathematics Studies, 158.
Princeton University Press, Princeton, NJ, 2005.

\bibitem[BouGol]{BG} Bourgain, J., Goldstein, M.
{\em On nonperturbative localization with quasi-periodic potential.}
Ann.\ of Math.~(2) 152 (2000), no.~3, 835--879.

\bibitem[BouGolSch]{BGS} Bourgain, J., Goldstein, M., Schlag, W.
{\em Anderson localization for Schr\"odinger operators on~$\IZ$ with
potentials given by the skew-shift.} Comm.\ Math.\ Phys.\ 220
(2001), no.~3, 583--621.

\bibitem[Cha]{C} Chan, J. {\em Method of variations of potential of quasi-periodic Schr\"odinger equation}, preprint 2005, to appear in GAFA.










\bibitem[GolKle] {GK} Goldstein, M., Klein, S. {\em Anderson localization
for random potentials with fast decaying correlations. } In preparation.

\bibitem[GolSch1]{GS1} Goldstein, M., Schlag, W.
{\em H\"older continuity of the integrated density of states
 for quasiperiodic Schr\"odinger equations and averages of shifts of subharmonic functions.}
Ann.\ of Math.~(2) 154 (2001), no.~1, 155--203.


\bibitem[GolSch2]{GS2} Goldstein, M., Schlag, W.
{\em Fine properties of the integrated density of states and a
quantitative separation property of the Dirichlet eigenvalues},
Preprint 2005.




\bibitem[Jit]{Jit} Jitomirskaya, S. {\em Metal--insulator transition for the almost Mathieu operator.}
Annals of Math.~150, no.~3 (1999).

\bibitem[Khin]{Kh} Khinchin, A.\ Ya.\ {Continued Fractions}, Dover, 1992.

\bibitem[Kle]{Kl} Klein, S. {\em Anderson localization for the discrete one-dimensional quasi-periodic Schr\"odinger operator
with potential defined by a Gevrey-class function.}  J.\ Funct.\
Anal.~218  (2005),  no.~2, 255--292.







\bibitem[Nat]{Nath} Nathanson, M. {Additive Number Theory I},
Springer.










\end{thebibliography}
\end{document}